\documentclass[draft]{amsart}
\usepackage{amsfonts, amsmath, amscd,amssymb,fullpage}
\usepackage{amssymb, amsbsy, amsthm,  amstext, amsopn, verbatim,
multicol}

\numberwithin{equation}{subsection}
\newtheorem*{thm}{Theorem}
\newtheorem*{prop}{Proposition}
\newtheorem*{lem}{Lemma}
\newtheorem*{sublemma}{Sublemma}

\newtheorem*{cor}{Corollary}

\newtheorem*{hypothesis}{Hypothesis}
\newtheorem*{question}{Question}
\newtheorem*{conjecture}{Conjecture}
\newtheorem*{question1}{Question 1}
\newtheorem*{question2}{Question 2}

\newtheorem*{problem}{Problem}

\newcommand{\reg}[1]{#1^{\text{reg}}}
\newcommand{\PP}{\mathcal{P}}
\newcommand{\BB}{\mathcal{B}}

\newcommand{\ep}{\epsilon}

\newcommand{\TT}{\mathbb{T}}

\newcommand{\h}{\mathfrak{h}}

\newcommand{\LL}{\mathcal{L}}
\newcommand{\CC}{\mathbb{C}}
\newcommand{\OO}{\mathcal{O}}

\newcommand{\hcdc}{{{}_c\hskip -0.8pt\mathcal{HC}_d}}
\newcommand{\NN}{\mathbb{N}}
\newcommand{\Z}{\mathbb{Z}}
\newcommand{\Q}{\mathbb{Q}}

\newcommand{\C}{\mathbb{C}}
\newcommand{\V}{\mathcal{V}}
 \newcommand{\ttt}{\textsf}

\newcommand{\irr}[1]{\ttt{Irrep}(#1)}
\newcommand{\hr}{\mathfrak{h}^{\text{reg}}}
\newcommand{\cxy}{\C [\h\oplus \h^*]}
\newcommand{\cx}{\C [\h]}
\newcommand{\cy}{\C [\h^*]}

\newcommand{\EE}{\mathbf{E}}

\newcommand{\WW}{{W}}
\newcommand{\UU}{U}

\newcommand{\AAA}{\mathbb{A}}
\newcommand{\JJJ}{\mathbb{J}}

\newcommand{\ddd}{S}

\newcommand{\twoslash}{{{/} \hskip -3pt {/}}}
\newcommand{\bmd}{\ttt{\text{-}bimod}}
\newcommand{\lgr}{\ttt{\text{-}grmod}}
 \newcommand{\lGr}{\ttt{\text{-}Grmod}}
 \newcommand{\lqgr}{\text{-}\ttt{qgr}}
 \newcommand{\lQgr}{\text{-}\ttt{Qgr}}
 \newcommand{\ltors}{\text{-}\ttt{tors}}
 \newcommand{\lTors}{\text{-}\ttt{Tors}}
 \newcommand{\lmod}{\text{-}\ttt{mod}}
 
 \newcommand{\fMd}{\text{-}\ttt{Filt}}
 \newcommand{\fmd}{\text{-}\ttt{filt}}
 \newcommand{\fbmd}{\text{-}\ttt{bifilt}}
\DeclareMathOperator{\odeg}{\ttt{ord}\text{-deg}}
 \DeclareMathOperator{\tdeg}{\ttt{tot}\text{-deg}}
\DeclareMathOperator{\Edeg}{\EE\text{-deg}}

\newcommand{\md}{\text{-mod}}

\DeclareMathOperator{\triv}{\ttt{triv}}
\DeclareMathOperator{\sign}{\ttt{sign}}
\DeclareMathOperator{\bideg}{\text{bideg}}

\DeclareMathOperator{\hin}{Hilb^n\CC^2}
\DeclareMathOperator{\hi}{Hilb(n)}
\DeclareMathOperator{\hio}{Z_n}
\DeclareMathOperator{\tgr}{\ttt{tgr}}
\DeclareMathOperator{\ogr}{\ttt{ogr}}
\DeclareMathOperator{\ord}{\ttt{ord}}
\DeclareMathOperator{\tot}{\ttt{tot}}
 
\DeclareMathOperator{\Hom}{Hom}

\DeclareMathOperator{\prj}{Proj}
\DeclareMathOperator{\spc}{Spec}
\DeclareMathOperator{\supp}{Supp}
\DeclareMathOperator{\mult}{mult}

\DeclareMathOperator{\coh}{\ttt{Coh} }
\DeclareMathOperator{\Qcoh}{\ttt{Qcoh} }

\DeclareMathOperator{\gr}{gr}

\DeclareMathOperator{\qgr}{\ttt{qgr}}  %{\ttt{Coh}}
\DeclareMathOperator{\QGr}{\ttt{Qgr}}  %{\ttt{Qcoh}}
 
 \DeclareMathOperator{\minp}{Min}
 
 \DeclareMathOperator{\ann}{ann}
\DeclareMathOperator{\Chro}{Char_0}
 \DeclareMathOperator{\Chr}{\mathbf{Char}}
  \DeclareMathOperator{\Ch}{\mathbf{Ch}}
\DeclareMathOperator{\Chw}{\mathbf{rCh}}

\begin{document}

\title[Rational Cherednik algebras and Hilbert schemes II]{Rational
Cherednik algebras and Hilbert schemes II:  representations and sheaves}
   \author{I. Gordon} \address{Department of Mathematics,
 Glasgow  University, Glasgow G12 8QW, Scotland}
\email{ig@maths.gla.ac.uk}
 \author{J. T. Stafford}
\address{Department of Mathematics, University of Michigan, Ann Arbor,
MI 48109-1109, USA.}
\email{jts@umich.edu}
   \thanks{The first author was supported by the Nuffield Foundation
    Grant NAL/00625/A and the Leverhulme trust. He would like to thank
    the University of Washington and the University of California at
    Santa Barbara for their hospitality while parts of this paper
    were written. The second author was supported in part by the
NSF   grant DMS-0245320.}

  \begin{abstract}  Let $H_c$ be the rational Cherednik algebra of type
 $A_{n-1}$ with spherical subalgebra $U_c=eH_ce$. Then $U_c$ is filtered by
 order of differential operators with associated graded ring $\gr
 U_c=\C[\h\oplus\h^*]^\WW$, where $\WW$ is the $n$-th symmetric group. Using the
 $\Z$-algebra construction from   \cite{GS} it is also possible to  associate
 to a filtered $H_c$- or $U_c$-module $M$ a coherent sheaf
 $\widehat{\Phi}(M)$ on the Hilbert scheme $\hi$.

  Using this technique, we study the representation theory of $U_c$ and $H_c$,
  and relate it to $\hi$ and to   the resolution of singularities $\tau:\hi\to
  \h\oplus\h^*/\WW$.
For example, we prove:
\begin{itemize}
\item If $c=1/n$, so that
$L_c(\triv)$ is the unique one-dimensional simple $H_c$-module,
then $\widehat{\Phi}  (eL_c(\triv)) \cong
 \OO_{\hio}$, where $\hio=\tau^{-1}(0)$
  is the {\it punctual Hilbert scheme}.
\item If $c=1/n+k$ for $k\in \NN$ then, under a canonical
filtration on the finite dimensional module $L_c(\triv)$, $\gr
eL_{c}(\triv)$ has a natural bigraded structure which coincides
with that on $\mathrm{H}^0(\hio, \LL^k)$, where
$\LL\cong\OO_{\hi}(1)$; this confirms  conjectures of Berest,
Etingof and Ginzburg. \item Under mild restrictions on $c$, the
characteristic cycle of $\widehat{\Phi}  (e\Delta_c(\mu))$ equals
$\sum_{\lambda}K_{\mu\lambda}[Z_\lambda]$, where
$K_{\mu\lambda}$ are Kostka numbers  and the $Z_{\lambda }$
are (known) irreducible components of $\tau^{-1}(\h/\WW)$.
 \end{itemize}
  \end{abstract}
   \maketitle
\vskip -12pt \vfill
\vbox{\tableofcontents}
\clearpage

\section{Introduction} \subsection{}\label{intro-1.1}
Let $c\in \C$. We assume throughout the paper that $c\not\in
\frac{1}{2}+\Z$ and, for simplicity, we will also assume in this
introduction that $c\notin \Q_{\leq 0}$. The aim of this paper is
to study the representation theory of  the rational Cherednik
algebra of type $A_{n-1}$, denoted $H_c$, and its spherical
subalgebra $U_c = eH_ce$ and to relate $H_c$- and $U_c$-modules to
sheaves over the Hilbert scheme of points in the plane.

Glimpses  of such a connection  appear in \cite{gordc} and
\cite{BEGfd} where it is shown that  finite dimensional  $H_c$-
and $U_c$-modules deform the sections of some  remarkable sheaves
on   $\hi$. A more formal relationship is provided by the main
theorem from  the companion paper \cite{GS}. In order to describe
that result recall that both $H_c$ and $U_c$ can be filtered by
degree of differential operators with associated graded rings $\gr
H_c\cong \cxy\ast\WW$ and $\gr U_c\cong \C[\h\oplus \h^*]^\WW$,
where  $\h\subset \C^{n}$ is the reflection representation of the
$n^{\mathrm{th}}$ symmetric group $\WW$; one can therefore regard
$U_c$ as a deformation of   $ \C[\h\oplus \h^*]^\WW$. The theorem
below provides a second way of passing to associated graded
objects which shows that $U_c$ can also be regarded as a
noncommutative deformation of  a homogeneous coordinate ring of
the Hilbert scheme $\hi$ (this scheme provides a crepant
resolution $\tau:\hi\to  \h\oplus \h^*/\WW$).

\subsection{Theorem}\label{intro-1.2} \cite[Theorem~1.4]{GS} {\it
 There exists a filtered $\Z$-algebra $B$ such that
\begin{enumerate}
\item
$U_c\md$, the category of finitely generated $U_c$-modules, is equivalent
to $B\lqgr$,
 the quotient category of finitely generated graded $B$-modules
modulo those of finite length; \item $\gr B$, the associated
graded ring of $B$, is isomorphic to (the $\Z$-algebra associated
with) the homogeneous coordinate ring $\bigoplus_{k \geq 0}
\mathrm{H}^0(\hi, \LL^k)$ for  a certain ample line bundle $\LL$
on $\hi$.
\end{enumerate}
}

\medskip
The construction of $B$ is described in more detail in
 Subsection~\ref{shift-defn}. In brief
$B=\bigoplus_{i\geq j\geq 0}B_{ij}$, where $B_{ij}$ is the
$(U_{c+i},\,U_{c+j})$-bimodule that provides a Morita equivalence
between $U_{c+i}$ and $U_{c+j}$, and multiplication in $B$ is given in matrix
fashion: $B_{ij}B_{jk}\subseteq B_{ik}$, but $B_{ij}B_{\ell k}=0$ if
$j\not=\ell$. One should regard $B$ as being a little like a Rees ring, but in a
situation where one may not have any proper ideals. Each $B_{ij}$ is
filtered by order of differential operators, which induces the desired
filtration $B=\bigcup \ord^n B$  on $B$.

\subsection{}\label{intro-1.3}
Theorem~\ref{intro-1.2} provides the following
recipe  for passing from a left
$U_c$-module $M$ with a good filtration $\Lambda$ to a coherent
sheaf $\Phi(M)=\Phi_{\Lambda}(M)$ on $\hi$.
Each  $M(k)=B_{k0}\otimes_{U_c} M$ has an  induced
  {\it tensor product filtration} $\Lambda$
defined by $\Lambda^n M(k) =\sum_{i\leq n}\ord^i B_{k0}\otimes \Lambda^{n-i}M$.
This   induces a filtration on
$\widetilde{M}=B\otimes M$, again written $\Lambda$, whose  associated graded
object $\gr_{\Lambda} \widetilde{M}$ is a finitely generated  $\gr B$-module.
Finally, by the theorem,
  $\gr_{\Lambda} \widetilde{M}$  corresponds to a coherent
sheaf $\Phi(M)$ on $\hi$.

This also works for  any filtered $H_c$-module $(N,\Gamma)$:  use the  induced
tensor product filtration $\Lambda $ on $eH_c \otimes_{H_c} N$ to construct a
sheaf $\widehat{\Phi}(N)=\Phi_{\Lambda}(eH_c \otimes N)$ on $\hi$. Under our
hypotheses on $c$,  tensoring with $eH_c $ gives a Morita  equivalence between
$H_c$ and $U_c$ \cite[Theorem~3.3]{GS}, so  $\widehat{\Phi}$ is really the same
construction  as $\Phi$.

The basic technique of this paper is to exploit the functors $\Phi$ and
$\widehat{\Phi}$ to better understand the representation theory
of $U_c$ and $H_c$.

\subsection{}\label{intro-1.11}  This works particularly well  for
$\OO_c$,  the   category of  finitely generated $H_c$-modules on
which $\C[\h^*]$ acts locally nilpotently. The basic building
blocks here are the {\it standard modules} $\Delta_c(\mu)$, where
$\mu\in \irr{\WW}$ is an irreducible representation of  $\WW$.
These are analogues of Verma modules over a simple complex Lie
algebra; in particular $\Delta_c(\mu)\cong \C[\h]\otimes \mu$ as
 left $\C[\h]$-modules.  It follows from this
that,  for any good filtration $\Lambda$ on $e\Delta_c(\mu)$,
 the  associated variety  $
\mathcal{V}(\mathrm{ann}_{\gr U_c}(\gr_\Lambda e\Delta_c(\mu))) =
\h/\WW\subset \h\oplus\h^*/\WW$ is independent of $\mu$.

An interesting phenomenon occurs if we follow the procedure of
\eqref{intro-1.3} since  the characteristic variety
$\mathcal{V} (\mathrm{ann}(\gr_\Lambda \widetilde{e\Delta_c(\mu)}))$ of
$\widetilde{e\Delta_c(\mu)}$ is   a subvariety of  $\hi$ that
 {\it does} depend upon $\mu$. We will state the result
using a finer invariant called  the {\it
 characteristic cycle}, $\Ch(\Delta_c(\mu))$, which counts the irreducible
components of  the characteristic variety  of $\widehat{\Phi}(\Delta_c(\mu))$ in $\hi$
with multiplicities (see \eqref{carcycle-subsec} for the formal definition).

\begin{thm} {\rm (Theorem~\ref{stand-cc})}
Let $\Delta_c(\mu)$ be the standard
$H_c$-module corresponding to $\mu\in\irr\WW$. Then
$$ \Ch (\Delta_c(\mu))   = \sum_{\lambda} K_{\mu
\lambda} [Z_{\lambda}],$$
where $K_{\mu\lambda}$ are   Kostka numbers as defined in \eqref{kostka-defn}
and the $Z_{\lambda}$'s are the irreducible components of
$Z=\tau^{-1}(\h/\WW)$ described in \eqref{Z-variety}.
\end{thm}

When $n=2$ and $c=1/2$  an easy direct proof of this and other theorems is given in
(\ref{s2-eg}--\ref{s2-eg3}) which the reader may like to use as an introduction
to the results of this paper.

\subsection{}\label{intro-1.45}  The above theorem illustrates a
fact that will appear repeatedly in this paper:  important modules $M$
over $H_c$ or $U_c$ correspond naturally to important sheaves and
combinatorial data  on $\hi$.

One reason why Theorem~\ref{intro-1.2} and the functors $\Phi$ and $\widehat{\Phi}$ provide
a bridge between Hilbert schemes and  Cherednik algebras is that the
construction of $B$  carries key elements from both theories. For
instance, the shift functors $(B_{j+1,j}\otimes -)$  are fundamental to the
theory of Cherednik algebras; they are also the analogue of the endo-functor $(\LL\otimes -)$
in $\coh(\hi)$: given a   $U_c$-module
$M$ with a good filtration $\Lambda$, we show that
$\LL\otimes\Phi_{\Lambda} (M) = \Phi_{\Gamma}(B_{c+1,c} \otimes M)$ for
the appropriate filtration $\Gamma$ (see Lemma~\ref{tensorline}).
As Proposition~\ref{1d-intro} and Theorem~\ref{intro-1.13}
illustrate, this is important since it enables one to
study whole families of modules simultaneously.

\subsection{}\label{intro-1.15} A consequence of Theorem~\ref{intro-1.11}
is that the Grothendieck group $K(\OO_c)$ is naturally isomorphic to the top
degree  Borel-Moore  homology group $H_{2n-2}(Z,\C)$
(see Corollary~\ref{borel}).  If
\cite[Remark~5.17]{GGOR} is confirmed,\footnote{This has been announced
recently by Rouquier.}   this can be considerably refined:
the characteristic cycles of the simple
modules in $\OO_c$ give a canonical basis of $H_{2n-2}(Z,\C)$
which corresponds to   Leclerc and Thibon's
lower canonical basis of the ring of symmetric functions.
 The details can be found in
Subsection~\ref{justify-1.15}.

\subsection{}\label{intro-1.5} The category $\OO_c$ also includes the  finite
dimensional $H_c$-modules: they occur precisely when $c=r/n$ for some $r\in
\NN$ with $(r,n)=1$. In such a case $H_c$ has a  unique finite dimensional
simple  module $L_c(\triv)$ (see \eqref{subsec-3.7}  for the notation) and
every finite dimensional $H_c$-module is a direct sum of copies of
$L_c(\triv)$. We note  that $L_{1/n}(\triv)$
is one-dimensional and, by
 \cite[Proposition~3.16]{GS}, that
$
 eL_{k+ 1/n}(\triv)\cong
B_{k0}\otimes_{U_{1/n}}eL_{1/n}(\triv),$ for $ k\in \NN.$

\begin{prop}\label{1d-intro} {\rm (Proposition~\ref{1d-prop})}
Give $L_{1/n}(\triv)$ the trivial filtration
$\Lambda^0=L_{1/n}(\triv)$. Then $\widehat{\Phi}(L_{1/n}(\triv))
\cong \OO_{\hio},$ where $\hio = \tau^{-1}({\bf 0})$ is the
punctual Hilbert scheme.
\end{prop}

\subsection{}\label{intro-1.13} Let $c=k+1/n$ and give
$eL_c(\triv)=B_{k0}\otimes eL_{1/n}(\triv)$ the tensor product
filtration induced from the trivial filtration on
$eL_{1/n}(\triv)$. Then $\gr eL_c(\triv)$ and (under an analogous
filtration)
  $\gr L_c(\triv)$ are  bigraded: one component of this
 grading is induced from the Euler grading on differential operators;
 the other comes from the fact that we have taken an associated graded
 module. Similarly, $\widehat{\Phi}(L_c(\triv))$ is a $(\C^*)^2$-equivariant
 sheaf on $\hi$.
Using Proposition~\ref{1d-intro} we are able to identify this sheaf,
and we can then describe the  corresponding bigraded characters of
$\gr L_c(\triv)$ and $\gr eL_c(\triv)$. Modulo a technical  change in the
filtrations,
this  gives the following result,  the first
part of which  confirms conjectures of Berest, Etingof and Ginzburg,
\cite[Conjectures~7.2 and 7.3]{BEGfd}.

 \begin{thm} {\rm (Theorem~\ref{beg-conj} and \eqref{bigrfd})}
 Let $c=k + 1/n$ for some $k\in \NN$. Then there are
  natural filtrations $\Lambda$
 on $eL_{c}(\triv)$ and $L_c(\triv)$   such that
\begin{enumerate}
\item there are bigraded isomorphisms
 $\gr_\Lambda L_c(\triv) \cong \mathrm{H}^0(\hio , \PP\otimes \LL^{k-1})$
 and $\gr_\Lambda eL_c(\triv) \cong \mathrm{H}^0(\hio,   \LL^k)$;
\item under this bigrading, $\gr_\Lambda eL_c(\triv)$ has Poincar\'e series
 equal to the $(t,q)$-Catalan-like number defined by
 \cite[(3.26)]{garhai2}.
\end{enumerate}
\end{thm}

\subsection{}\label{intro-1.6} A further example of an  $H_c$-module
corresponding to an important sheaf is given by the module $H_c$ itself.  Under
the order filtration $\ord$ of differential operators, Theorem~\ref{cohpro}
shows that  $\widehat{\Phi}_{\ord} (H_c) \cong \PP$, the {\it Procesi bundle} on
$\hi$.  This is a remarkable rank $n!$ bundle that arises from Haiman's
isospectral Hilbert scheme. See \eqref{PP-defn} or \cite{hai3} for more
details.

\subsection{}\label{intro-1.14}   There is a second  construction that produces
 a coherent sheaf on $\hi$ from an $H_c$-module $N$ with a good
filtration $\Gamma$, and which provides a potential connection between
our $\Z$-algebra construction and the Fourier-Mukai equivalences of Bridgeland-King-Reid.
The construction  is described as follows. Theorems of Bridgeland-King-Reid and   Haiman provide an equivalence
of categories $D(\cxy
\ast\WW) \xrightarrow{\sim} D(\hi)$ from the bounded derived category of
finitely generated $\cxy\ast\WW$-modules to that of
coherent sheaves on $\hi$. Applying this equivalence to the
finitely generated $\cxy\ast \WW$-module $\gr_{\Gamma} N$,
regarded as a complex concentrated in degree zero,
 gives a bounded complex $\Psi (N)=\Psi_{\Gamma} (N)$  of coherent sheaves
  on $\hi$.

In all the examples we have calculated, {\it the complex $\Psi (N)$ is
quasi-isomorphic to $\widehat{\Phi}(N)$}. We ask whether this is true in general; we
   are only able to show that there is
a surjective mapping from the  zeroth cohomology sheaf $\Psi^0(N)$ onto
$\widehat{\Phi}(N)$, see Proposition~\ref{compBKR}.
If this question had a positive answer it would imply
that  all higher
cohomology of $\PP^*\otimes \widehat{\Phi}(N)$ must vanish. A vanishing theorem
in this spirit,  which has consequences for symmetric function theory,
has been proposed in \cite[Conjecture~3.2]{hai1}, see \eqref{cohvan}.

\subsection{}\label{intro-1.16}
 Our final application is an analogue of the  Conze embedding \cite{con} for $U_c$.
 The Dunkl-Cherednik map \cite[Proposition~4.5]{EG} provides an inclusion
 $U_c \hookrightarrow e(D(\hr)\ast \WW) e \cong D(\hr)^\WW$ and it is
 straightforward to see that this actually embeds $U_c$  into
 the Weyl algebra $D(\h /\WW)$. Our analogue of the Conze embedding gives a
 precise description of the $U_c$-module structure of this algebra:

\begin{prop} {\rm (Proposition~\ref{subsec-8.4})}
Under the natural identifications,
  $D(\h/\WW)=\bigcup_{k\geq 0} B^*_{k0}$,  where the
$B_{k0}^{\ast} = \Hom_{U_{c }}(B_{k0}, U_{c})$ are projective left
$U_c$-modules. Thus $D(\h/\WW)$ is  a flat      left
 $\UU_c$-module.
\end{prop}

Once again, this corresponds to a fact about Hilbert schemes: in
this case it is the result from \cite[Section~2]{haidis} that
shows that the   subset $U_{(1^n)}$ of $\hi$ is isomorphic to the
affine space $\mathbb{A}^{2n-2}$.

\subsection{Plan}\label{intro-1.19} The paper is organised as follows. In
Section~\ref{zalg} we study  filtered $\Z$-algebras and their modules and show
that classical facts concerning characteristic varieties and characteristic
cycles for filtered noetherian algebras do generalise to   $\Z$-algebras.
Section~\ref{rat-intro} provides some background material while
Section~\ref{cohsh} sees the construction of  the functors $\Phi$ and
$\widehat{\Phi}$ and a discussion of their relationship  to the equivalences of
Bridgeland-King-Reid. In Section~\ref{var-sect} we  study the effect of $\Phi$
and $\widehat{\Phi}$ on the finite dimensional  simple $U_c$- and $H_c$-modules, in
particular using this to determine  their bigraded characters. Our results on
characteristic cycles for  $\OO_c$ are presented in Section~\ref{ses-sect},
whilst analogous results for Harish-Chandra modules are proved in
Section~\ref{bimodule-sect} and
the  Conze embedding is constructed
in  Section~\ref{conze}.

%%%%%%%%%%%%%%%%%%%%%%%%%%%%%%%%%%%%%%%%%%%%

 \section{Graded and filtered modules for $\Z$-algebras} \label{zalg}

\subsection{}\label{Z-alg-defn}  Throughout this paper a {\it $\Z$-algebra}
 will mean a {\it lower triangular $\Z$-algebra}. By definition, this is a
 (non-unital)  algebra $B=\bigoplus_{i\geq j\geq 0} B_{ij}$, where
 multiplication is defined in matrix fashion:
 $B_{ij}B_{jk}\subseteq B_{ik}$ for $i\geq j\geq k\geq 0$
 but $B_{ij}B_{\ell k}=0$ if $j\not=\ell$.
 Although $B$ cannot have a unit element, we do require that each
 subalgebra $B_{ii}$  has a unit element $1_i$ such that
 $1_ib_{ij}=b_{ij}=b_{ij}1_j$, for all  $b_{ij}\in B_{ij}$.

Let $B$ be a $\Z$-algebra. We define the category $B\lGr $  to be
the category of $\NN$-graded left  $B$-modules $M=\bigoplus_{i\in
\NN} M_i$ such that $B_{ij}M_j\subseteq M_i$ for all $i\geq j$
and $B_{ij}M_k=0$ if $k\not=j$. Homomorphisms are defined to be
graded homomorphisms of degree zero. The subcategory of noetherian
graded left $B$-modules will be denoted $B\lgr $. In all  examples
considered in this paper $B\lgr $ will consist precisely of the
finitely generated graded left  $B$-modules.

A module $M\in B\lGr $ is {\it bounded} if $M_n = 0$ for all but
finitely many $n\in \Z$
 and {\it  torsion} if it is a direct limit of bounded modules. We denote the
  category  of torsion modules
  by $B\lTors $. The corresponding subcategory of $B\lgr $
is denoted $B\ltors $; thus  $B\ltors $
 just consists of bounded, noetherian modules.
 Finally, in order to obtain geometric categories, we construct
the quotient categories $B\lQgr = B\lGr /B\lTors $ and $B\lqgr  =
B\lgr /B\ltors $.\label{qgr-defn} We write $\pi(M)$ for the image
in $B\lQgr $ of $M\in B\lGr$.

\subsection{} \label{zalgex1} We will be interested in two types of examples
of $\Z$-algebras.

First, suppose that
$S= \bigoplus_{n\geq 0} S_n$ is a noetherian
$\NN$-graded  algebra. To $S$ we
  associate a $\Z$-algebra
   $\widehat{S}=\bigoplus_{i\geq j\geq 0}\widehat{S}_{ij}$
   by setting $\widehat{S}_{ij} = S_{i-j}$ with multiplication
induced from  that in $S$. Define  categories $S\lGr,\dots,S\lqgr$
in the usual manner. In particular,
 $S\lGr$ denotes the category
of $\Z$-graded $S$-modules, from which the other definitions
follow as in the last paragraph. As explained in \cite[(5.3)]{GS}
we have equivalences $S\lQgr\simeq \widehat{S}\lQgr$ and $ S\lqgr\simeq
\widehat{S}\lqgr$.
 If $S$ is  commutative and generated as an
$S_0$-algebra by $S_1$,  then $S\lqgr$ and hence $\widehat{S}\lqgr$
 are equivalent to
$\coh\prj(S)$, the category of coherent sheaves on
$\prj(S)$   \cite[Proposition~3.3.5]{EGA}.

Second, suppose that  $\{R_n : n\in \NN\}$ are Morita equivalent
noetherian algebras, with the equivalence induced from the
progenerative $(R_{n+1},R_{n})$-bimodules $P_n$. Define
$R_{ij} = P_{i-1}\otimes_{R_{i-1}}\otimes\cdots
\otimes_{R_{j+2}}P_{j+1}\otimes_{R_{j+1}}P_j$ and
$R_{jj}=R_j$, for $i>j\geq 0$. The corresponding $\Z$-algebra
$R_\Z=\bigoplus_{i\geq j\geq 0} R_{ij}$ will be called the {\it
Morita $\Z$-algebra associated to the data $\{R_n,P_n : n\in
\NN\}$}. We again have a useful description of $R_\Z \lqgr$,

\begin{lem} {\rm \cite[Lemma~5.5]{GS}} Let $R_{\Z}$ be a Morita $\Z$-algebra
as above. Then
\begin{enumerate}
  \item[(1)] each finitely generated graded left $R_\Z$-module is noetherian;
  \item[(2)] the
 association $\phi: M \mapsto \bigoplus_{n\in \NN}R_{n0}\otimes_{R_0}M$
 induces an equivalence of categories
 $R_0\lmod\xrightarrow{\sim}R_\Z\lqgr$.
 \end{enumerate}
\end{lem}

\subsection{Associated graded techniques}\label{filter-discuss} A number
standard constructions for unital noetherian rings extend routinely to
$\Z$-algebras: those that will be useful in this paper are the notions of  associated graded
modules and characteristic cycles.

In order for filtrations of a  $\Z$-algebra $R$ to work properly,
the filtrations   on the various components $R_{ij}$ need to be
compatible. The abstract requirements for this are quite lengthy
and so we will assume that we are in the following special
situation.
\begin{enumerate}
\item {\it $R=R_\Z = \bigoplus_{i\geq j\geq 0} R_{ij}$ is the
Morita $\Z$-algebra associated to the data
 $\{R_n,P_n : n \in \NN\}$.}
\item {\it Each $R_i$ is a noetherian domain, and each $P_i$ has
rank one (on one side and therefore the other). As a result, we
may identify each $R_{n+1}=\text{End}_{R_n}(P_n)$ with a subring
of the quotient division ring, say $D$, of $R_0$.}
\item {\it
There is a ring $E\subset D$ containing each $R_i$ and $P_i$ such
that $E$ has an ascending $\NN$-filtration $F^kE$ with a
commutative associated graded ring $\gr_FE$. }
\end{enumerate}
The point of these assumptions is that we can take the induced
filtrations $\{F^k R_{ij}\}$ on the $R_{ij}$ and hence on $R$
itself. It is then routine to check that this filtration is
compatible with the ring structure of $R$ and so we obtain an {\it
associated graded $\Z$-algebra}
$$\gr_F R = \bigoplus_{i\geq j\geq 0}\gr_F R_{ij},
\qquad\text{where}\qquad \gr_F R_{ij} = \bigoplus_{n\in
\NN}\frac{F^nR_{ij}}{F^{n-1}R_{ij}}.$$ The final assumption
that we will make is:
\begin{enumerate}
\item[(4)] {\it $\gr_F R\cong \widehat{S}$, the $\Z$-algebra
associated to some finitely generated commutative graded domain $S$.}
\end{enumerate}

Conditions (1)--(3) occur reasonably frequently. One case is given
by taking the $R_i=U_{c+i}$ (see \eqref{balg-defn}), but one
can also take appropriate factor rings of the enveloping
$U(\mathfrak{g})$ of a semisimple Lie algebra $\mathfrak{g}$,
where the progenerators arise from translation functors (see, for example,
\cite{JS}). Condition (4) is  more stringent,
 although Theorem~\ref{intro-1.2}(2) shows that it does
apply when $R_i=U_{c+i}$.

\subsection{}\label{graded-section}
Keep hypotheses (1)--(3) of \eqref{filter-discuss} and let
$M=\bigcup \Lambda^iM$ be a filtered left  $R_0$-module.
Each module $M(i)=R_{i0}\otimes_{R_0}M$ is then filtered by the
 {\it tensor product filtration}
 \begin{equation}\label{indfil}
 \Lambda^kM(i)
= \sum_{\ell\in\Z} F^\ell R_{i0}\otimes \Lambda^{k-\ell}M.
\end{equation}
We therefore obtain an associated graded $\gr_FR_i$-module
$\gr_\Lambda M(i) =\bigoplus \Lambda^kM(i)/\Lambda^{k-1}
M(i)$. Summing over all $i$ gives  an associated graded
$\gr_FR$-module $\gr_\Lambda \widetilde{M} = \bigoplus
\gr_\Lambda M(i)$ for the $R$-module
$\widetilde{M}=\bigoplus_{i\in \Z}M(i)$. When
Hypothesis~\ref{filter-discuss}(4) holds we can and will regard
$\gr_\Lambda M$ as a graded $S$-module.

We say that $\Lambda$ is a {\it good filtration}\label{good-defn}
 on $M$ if $\gr_\Lambda M$
is a finitely generated $\gr_FR_0$-module. Similarly,
a filtration $\Lambda$ on $\widetilde{M}$ is
 {\it good} if $\gr_\Lambda\widetilde{M}$
  is a finitely generated $\gr_FR$-module. If $m\in \Lambda^iM\smallsetminus
  \Lambda^{i-1}M$, we write $\sigma(m)=[m+\Lambda^{i-1}M]\in \gr_\Lambda M$ for
  the {\it principal symbol} of $m$\label{princ-symbol}.

\subsection{Lemma}\label{graded-lemma} {\it Keep  Hypotheses (1)--(3) of
\eqref{filter-discuss} and let $M$ be a left  $R_0$-module.
 \begin{enumerate}
 \item If $\Lambda$ is a good filtration on $M$ then the induced
 filtration $\Lambda$ on $\widetilde{M}$ is good.

  \item Given two good filtrations $\Lambda$ and $\Gamma$ on $M$, there exists
  a positive integer $t$ such that, for all  $n$ and $k$,
  $$\Lambda^nM(k)\  \subseteq\  \Gamma^{n+t}M(k)\  \subseteq\
  \Lambda^{n+2t}M(k). $$

  \item Assume that Hypothesis \eqref{filter-discuss}(4)
  also holds and regard
  $\gr_\Lambda \widetilde{M}$ as an object in $S\lgr$.
   Define the \emph{associated radical ideal}   $N=N_{\widetilde{M}}$
   of $\widetilde{M}$ to be the radical ideal $N=\sqrt{I}$ of
  $I=\mathrm{ann}_S\big(\mathrm{gr}_{\Lambda}\widetilde{M}\big) .$
  Then $N$ is independent of the choice of good filtration $\Lambda$.
  \end{enumerate}
   }

 \begin{proof}
 (1) Since the filtration $\Lambda$ on $M$ is good, we may pick
 generators $m_i$ for $M$ with
 $m_i\in \Lambda^{j_i}M\smallsetminus \Lambda^{j_i-1}M$ such that the
 principal symbols $\sigma(m_i)$ generate $\mathrm{gr}_\Lambda M$.
    Pick a free $R_0$-module $G=\bigoplus R_0 g_i$ and filter
    $G$ by putting $g_i$ in degree $j_i$; thus
$G$ is filtered by    $ \Phi^nG=\bigoplus (F^{n-j_i}R_0) g_i$. Then the
natural surjection $\alpha : G\to M$ given by $g_i\mapsto m_i$ is
\emph{filtered surjective} in the sense that, for each $n$, we
have $\alpha(\Phi^nG)= \Lambda^nM$.
 Equivalently, $\alpha$ induces a surjection
$\gr \alpha : \mathrm{gr}_{\Phi} G\to \mathrm{gr}_\Lambda M$.

Now consider the induced map $G(k)=R_{k0}\otimes_{R_0}G \to M(k)$.
The tensor product  filtration \eqref{indfil} on  $G(k)$ induces
a surjection $\alpha_k : \Phi^nG(k)\twoheadrightarrow
 \Lambda^nM(k)$ for all $n$ and $k$ and hence a surjection of graded
 groups $\mathrm{gr}\,  \alpha_k : \mathrm{gr}_{\Phi}G(k)\to
 \mathrm{gr}_{\Lambda }M(k)$.
 Since  $G$ is free, we have,  for each $n$ and $k$,
 $$\displaystyle
   \Phi^nG(k)= \displaystyle\sum_{j} F^{n-j}R_{k0}\otimes \Phi^{j}G =
    \sum_{j,t} F^{n-j}R_{k0}\otimes (F^{j-j_t}R_0)g_t
   = \sum_{t} (F^{n-j_t}R_{k0}) \otimes g_t.
  $$
 The filtration $\Phi G(k)$ is therefore the natural filtration
 on the (weighted) direct sum of copies of $R_{k0}$ and so
  the associated graded module
 $\gr_\Phi\widetilde{G}$ is  just the weighted
 sum $E=\bigoplus_t S[j_t]$.
 Since $gr_{\Lambda}\widetilde{M}$ is a
 homomorphic image of $E$,  it is  finitely generated.

(2) By  \cite[Corollary~6.12]{kralen} the desired equation holds
for $k=0$. By \eqref{indfil} it then  holds for each $k$, with the
same value of $t$.

(3) By using intermediary filtrations, as in
\cite[Lemma~8.6.11]{MR},
 we may assume in part (2)  that $t=1$ for $k=0$ and hence for all $k$.
 As in the proof of \cite[Lemma~8.6.12]{MR}, we then
 obtain short exact sequences
 \begin{equation}\label{trick1}
 0\to Z\to \gr_{\Lambda}(\widetilde{M}) \to Y\to 0\qquad\text{and}\qquad
 0\to Y\to \gr_{\Gamma}(\widetilde{M}) \to Z[1]\to 0,
 \end{equation}
 where $Y=\bigoplus_{k}\bigoplus_{n} \Lambda^{n+1}M(k)/\Gamma^nM(k)$
 and $Z=\bigoplus_{k}\bigoplus_{n}\Gamma^{n}M(k)/\Lambda^nM(k)$.
If $J=\ann_S(\gr_\Gamma\widetilde{M})$, then it follows that
 $I^2\subseteq J$ and $J^2\subseteq I$.
 \end{proof}

\subsection{Definition}\label{charvar-defn}
 Keep the hypotheses of Lemma~\ref{graded-lemma}(3). Then $N_{\widetilde{M}}$
 is a graded ideal of $S$ and
  we  define the \emph{characteristic variety}
$\Chr(M)$ of $M$ to be the projective
subvariety  $\V(N_{\widetilde{M}})$ of $ \mathrm{Proj}(S)$
defined by $N_{\widetilde{M}}$. The lemma ensures that this subvariety is
independent of the choice of good filtration $\Lambda$, and so it is
an invariant of $M$.
A second variety associated to $M$ is
$$\Chro(M)=\V(\mathrm{ann}_{\gr R_0}(\gr_\Lambda M))\subseteq
\mathrm{Spec}(\gr R_0).$$  In order to avoid any possible
confusion with $\Chr M$ we will always call $\Chro M$ the {\it
associated variety of $M$}, as it is sometimes termed in the
literature.

\subsection{Characteristic cycles}\label{carcycle-subsec}
In the classical situation of filtered noetherian rings it is
often useful to refine the characteristic variety of a module to
the characteristic cycle; see, for example, \cite[Remark~5.7]{Bj2}.
As we show next, this also works for
$\Z$-algebras.

Keep the hypotheses and notation of Lemma~\ref{graded-lemma},
and write  $\minp \widetilde{M}$ for the
set of prime ideals  minimal over the associated
 radical ideal $N_{\widetilde{M}}$.
 If $P\in \minp  \widetilde{M}$, we
define $n_{\widetilde{M},P}$ to be the length of the (necessarily
finite length) $S_P$-module
 $M_P$. We note that if $S_{[P]}=S[\mathcal{C}^{-1}]_0$, where
 $\mathcal C$ denotes the set of homogeneous
 elements in $S\smallsetminus P$, then
  $n_{\widetilde{M},P}$ also equals the length of the
  $S_{[P]}$-module $\gr_\Lambda(\widetilde{M})S[\mathcal{C}^{-1}]_0$.
 The {\it characteristic cycle} $\Ch M$ is defined to be
 \begin{equation}\label{charcycle-defn}
\Ch M \ = \  \sum_{\minp \widetilde{M}}
n_{\widetilde{M},P}\V(P) .
\end{equation}
Clearly this is a finite sum with finite coefficients. By the next
result it is also a well-defined notion. We remark that the more usual
associated cycle of $\gr_\Lambda M$ will not be used in this paper and so
there should be no confusion in our notation.

\begin{lem}\label{ch-is-wd} Keep the hypotheses (1--4) of
\eqref{filter-discuss} and let $M$ be a left $R_0$-module.
 Then the characteristic cycle $\Ch M$
is independent of the choice of the good filtration $\Lambda$
used in its definition.
\end{lem}

\begin{proof} Suppose that the module $M$ has two good filtrations
$\Lambda$ and $\Gamma$ and, for the moment, write
 $\Ch_\Lambda {M}$ for the characteristic cycle
 of ${M}$ induced from $\Lambda.$  As in
 the proof of Lemma~\ref{graded-lemma}(3) we may assume
 that $t=1$ in the conclusion of
 part (2) of that lemma, and hence that we have the two
 short exact sequences \eqref{trick1}. It follows immediately from
 those sequences that $(\gr_\Lambda\widetilde{M})_P$ and
 $(\gr_\Gamma\widetilde{M})_P$ have the same (finite) length
 for any $P\in \minp \widetilde{M}$. In other words,
 $\Ch_\Lambda {M}=\Ch_\Gamma {M}$.
    \end{proof}

\subsection{}\label{cycles-ses}
Assume that $R_0$ satisfies  hypotheses (1--4) of
\eqref{filter-discuss}.
 We will need to understand how characteristic cycles relate to short
exact  sequences of $R_0$-modules. In general one does not obtain additivity on
such sequences since  embedded primes cause problems---a standard example is
given by  the short exact
 sequence of $\C[x]$-modules $0\to (x)\to \C[x]\to \C\to 0$.
 The standard way around this is to only look at the components of the
 characteristic variety of maximal dimension
 and we do the same here.

  Thus, if $M$ is a finitely generated filtered
 $R_0$-module, write $\minp'\widetilde{M}$ for the subset of
 $\minp \widetilde{M}$ consisting of those prime ideals
 $P$ for which $\dim \V(P) = \dim \V(N_{\widetilde{M}})$ and
 let $\Chw M$  be the corresponding {\it restricted characteristic
 cycle}:
\begin{equation}\label{cyclew-defn}
\Chw M  = \sum_{\minp' \widetilde{M}}
n_{\widetilde{M},P}\V(P).
\end{equation}
By Lemma~\ref{ch-is-wd} this too is well-defined.

\begin{lem}\label{ses-cc} Assume that $R_0$ satisfies hypotheses (1--4) of
\eqref{filter-discuss} and let $0\to A\to B\xrightarrow{\beta}
C\to 0$ be a short exact sequence of finitely generated left
$R_0$-modules.
 Then precisely one of the following
 cases occurs:
\begin{enumerate}
\item[(a)]  $\dim \Chr( A )<\dim \Chr(
B)=\dim\Chr( C)$ and $\Chw  B = \Chw
 C$; \item[(b)]  $\dim \Chr( A)=\dim
\Chr( B)>\dim\Chr( C)$ and $\Chw
B = \Chw  A$; \item[(c)] $\dim \Chr( A)=\dim
\Chr( B)=\dim\Chr( C)$ and $\Chw
B = \Chw  A+\Chw  C$.
\end{enumerate}
\end{lem}

\begin{proof}
Choose a good filtration $\Gamma$ on $B$ and
 give $A$ and $C$ the induced
filtrations $ \Gamma^t A =  \Gamma^t B \cap A$ and $ \Gamma^t C =
\beta(\Gamma^t B).$ It is an easy and standard exercise  to check
that this induces a short exact sequence of $\gr R_0$-modules and
that the analogue of the lemma holds for
 the restricted characteristic cycles of those graded modules
 (see, for example, \cite[Formula~5.8]{Bj2}).
Unfortunately the induced tensor product filtration on
$\widetilde{A}$, etc,  need not be filtered exact and so one will
not necessarily obtain a short exact sequence of associated graded
modules. So the proof needs to be a little more involved.

In order to simplify the notation, we want to assume that $S$ is
generated in degree one, which  is easy to arrange.  By
\cite[Lemme~2.1.6(v) and Proposition~2.4.7(i)]{EGA}, some
 Veronese ring $S^{(n)}=\bigoplus S_{in}$ is generated
by $S^{(n)}_0=S_0$ and $S^{(n)}_1=S_n$ while
$\mathrm{Proj}(S)=\mathrm{Proj}(S^{(n)})$ for any such $n$.
 Thus if we pass to  the Veronese $\Z$-algebra
  $R^{(n)}=\bigoplus R_{in,jn}$, with associated graded ring
$\gr_FR^{(n)}=\widehat{S^{(n)}}$,
  then the associated  cycle of $R^{(n)}\otimes_{R_0}B$
 will equal that of
$\widetilde{B}$. Thus we may assume that $S$ is generated in degree one.

Consider the induced short exact sequence
 $0\xrightarrow{} \widetilde{A}\xrightarrow{}
\widetilde{B}\xrightarrow{{\widetilde{\beta}} }
 \widetilde{C}\xrightarrow{} 0$. Give $\widetilde{B}$ the
tensor product filtration $\Gamma$, which it is convenient to
write as $\widetilde{\Gamma}$.
Give $\widetilde{A}$ the induced
 filtration   $\widetilde{\Gamma}$; thus
$\widetilde{\Gamma}^mA(k) = A(k)\cap \widetilde{\Gamma}^mB(k)$
for all $m,k$. It is unnecessary to do this with $\widetilde{C}$
since, as in the proof of Lemma~\ref{graded-lemma}(1),
${\Gamma}^mC(k) =\widetilde{\beta}( \widetilde{\Gamma}^mB(k))$ for all
$m,k$. It is immediate that $0\to
\widetilde{\Gamma}^m\widetilde{A}\to
\widetilde{\Gamma}^m\widetilde{B}\to
 {\Gamma}^m\widetilde{C}\to 0$
is exact for each $m$ and hence that the associated complex of
graded modules
\begin{equation}\label{ses-cc1}
0\to \gr_{\widetilde{\Gamma}}\widetilde{A}\to
\gr_{\widetilde{\Gamma}}\widetilde{B}\to \gr_{
{\Gamma}}\widetilde{C}\to 0
\end{equation}
is exact.

 By Lemma~\ref{graded-lemma}(1),  $\gr_{\widetilde{\Gamma}}\widetilde{B}$
 is finitely generated and therefore, by
Hypothesis~\ref{filter-discuss}(4),
 so  is $\gr_{\widetilde{\Gamma}}\widetilde{A}$. Pick an integer $t$
such that
 $\gr_{\widetilde{\Gamma}}\widetilde{A}$ is generated
 by $\bigoplus_{j\leq t}\gr_{\widetilde{\Gamma}}A(j)$.
 Define a new filtration $\Lambda$ on
 $\widetilde{A}_{\geq t} =\bigoplus_{i\geq t}A(i)$
 by defining $\Lambda =\widetilde{\Gamma}$ on
 $A(t)$ and  the   tensor product filtration
 thereafter; thus $\Lambda^mA(i) = \sum_j F^jR_{it}\otimes
 \widetilde{\Gamma}^{m-j}A(t)$, for $i>t$.
 The choice of $t$ and the fact that $S$ is generated in degree one
 ensures that
$$\gr_{\widetilde{\Gamma}}A(m) =
 \sum_{j\leq t} S_{m-j}\gr_{\widetilde{\Gamma}}A(j)
= S_{m-t}\gr_{\Lambda}A(t),$$ for all $m\geq t$. Pulling this
back to   $\widetilde{A}_{\geq t}$ gives
$$\widetilde{\Gamma}^jA(m)
= \sum_{v\geq 0}   F^vR_{mt} \otimes \Lambda^{j-v}
(R_{t0}\otimes A)=\Lambda^jA(m),$$ for any $ m\geq t $ and $
j\geq0.$ In other words, $\Lambda=\widetilde{\Gamma}.$

Now consider characteristic cycles. Under the tensor product
filtration $\gr_\Gamma \widetilde{B}$ and $\gr_{{\Gamma}}
\widetilde{B}_{\geq t}$ only differ in the first $t$ terms and so
in $\mathrm{Proj}(S)$ they have the same characteristic variety
and restricted characteristic cycle. So we may work with
$\widetilde{B}_{\geq t}$. But now the choice of good  filtration
on  $A(t)$ is irrelevant, so we may choose the filtration
$\Lambda=\widetilde{\Gamma}$ on $A(t)$ and the original
filtrations $\widetilde{\Gamma}=\Gamma$ on $B(t)$ and $C(t)$. By
\eqref{ses-cc1} we therefore obtain a short exact sequence
\begin{equation}\label{ses-cc3}
0\to \gr_{\Lambda}\widetilde{A}_{\geq t} \to
\gr_{\Gamma}\widetilde{B}_{\geq t} \to \gr_\Gamma
\widetilde{C}_{\geq t}\to 0,
\end{equation}
for which the restricted characteristic cycles of the three terms
are $\Chw A$, respectively $\Chw B$ and
$\Chw C$.

It is now routine to see that the conclusion of the lemma is
satisfied. Indeed, if $P\in \minp'\widetilde{B}$, then localising
\eqref{ses-cc3} at $P$ gives a short exact sequence of finite
dimensional $S_P$-modules and so the dimensions add:
$n_{\widetilde{B},P}= n_{\widetilde{A},P}+n_{\widetilde{C},P}.$
\end{proof}

%%%%%%%%%%%%%%%%%%%%%%%%%%%%%%%%%%%%%%%%%%%%
\section{Rational Cherednik algebras and Hilbert schemes}
\label{rat-intro}

In this  section we collect some of the basic material we need
concerning rational Cherednik algebras and Hilbert schemes.

\subsection{The rational Cherednik algebra of type $A$}\label{chered-defn}
 Let $\WW=\mathfrak{S}_n$
\label{symmetric-defn} be the symmetric group on $n$ letters,
regarded as the Weyl group of type $A_{n-1}$ acting on its
$(n-1)$-dimensional reflection  representation $\h\subset \C^n$ by
permutations. We will  always identify
 $\h^*$ with its image in the coordinate ring
  $\C[\h]$ and we
 fix a basis $\{x_i\}$ of $\h^*$.
 Let $\mathcal{S}=\{s=(i,j)\ \text{with}\
i<j\}\subset \WW$ \label{involution-defn} denote the reflections,
with reflecting hyperplanes  $\alpha_s=0$. The  polynomial
$\delta = \prod_{s\in \mathcal{S}} \alpha_s \in \C[\h]$ \label{delta-defn}
will be important, in part  because
  $\hr=\h \setminus \{\delta=0\}$ is
the subset of $\h$ on which the action of $\WW$ is free.
Notice that $\delta$ transforms under $\WW$ by the
{\it sign representation $\sign$}.\label{sign-defn}

Given a variety $Z$ we will write $D(Z)$ for the  ring of differential
operators on  $Z$. Let $D(\hr)\ast\WW$  denote the skew group ring
of $D(\hr)$ by $\WW$; thus, by definition, $wf=w(f)w$, for all $f\in D(\hr)$ and
$w\in \WW$. For $c\in \C$ the {\it rational Cherednik algebra of type $A_{n-1}$}
is the  $\C$-subalgebra $H_c$\label{hc-defn} of
$D(\hr)\ast \WW$  generated by the
multiplication operators $\h^*\subseteq \C[\h]=\C[x_1,\dots,x_{n-1}]$,
 the group algebra $\C\WW$, and the {\it Dunkl operators}
\begin{equation} \label{dunkop} y_i =
\partial_{i} - \sum_{s\in S} c\alpha_s(y_i)\alpha_s^{-1}(1-s),
\qquad \text{where}\  \partial_{i}= \partial/\partial
x_i.\end{equation}

\subsection{} By \cite[Theorem~1.3]{EG} there is a
 Poincar\'{e}-Birkhoff-Witt isomorphism of $\C$-vector spaces
\begin{equation}
\label{PBW} \C[\h]\otimes_{\C} \C \WW \otimes_{\C} \C[\h^*]
\xrightarrow{\sim} \ H_c,
\end{equation}
under which $\C[\h^*]$ identifies with the subalgebra of $H_c$ generated
by  $\h=\bigoplus_{i=1}^{n-1} \C y_i$.

Equation~\ref{PBW} can be interpreted as follows.  For any variety $Z$,
one has a natural filtration on
$D(Z)$ by order of operators and this induces a filtration on
$D(\hr)\ast \WW$ and its subalgebras by defining elements of $\WW$
to have order zero. If  $R$ is a subalgebra (or subset) of
$D(\hr)\ast \WW$, let $\ord^nR$ denote the operators of order $\leq n$ in $R$.
 \label{order-filt-defn} The  associated graded ring
of $R$ will be written $\ogr(R)=\bigoplus \ogr^n(R)$, where
$\ogr^n(R)=\ord^n(R)/\ord^{n-1}(R)$ and the resulting graded
structure of $\ogr(R)$ will be called the {\it order} or {\it
$\ogr$ gradation}.\label{ogr-defn}

If we filter $H_c$ in this way, then \eqref{PBW}
implies that \label{filt-defn} $\ord^0 H_c = \C [\h]\ast W$, $\ord^1H_c
= \h + \ord^0H_c$ and
 $\ord^iH_c = (\ord^1H_c)^i$ for $i>1$. Moreover,  the  associated graded ring
 $\ogr H_c=\bigoplus \ogr^n H_c$ is  isomorphic to the skew group ring
 $\C[\h\oplus \h^*]\ast \WW$.

\subsection{The spherical subalgebra} Let $e\in \C \WW$ be the trivial
idempotent and $e_-\in \C\WW$\label{e-defn} be the sign idempotent; thus
$e = |\WW|^{-1} \sum_{w\in \WW} w$  and
 $e_- = |\WW|^{-1} \sum_{w\in\WW} \text{sign}(w)w.$
The {\it spherical subalgebra} of $H_c$ is the algebra
$\UU_c=eH_ce$\label{spherical-defn} while  the related algebra
$\UU^-_c=e_-H_ce_-$ is called the {\it anti-spherical subalgebra}.
By    \cite[Proposition~4.1]{BEGfd} these algebras are related
through the   identity
\begin{equation} \label{conj} \UU_c \ =\ \delta^{-1} \UU^-_{c+1}\delta
\ = \ e\delta^{-1} H_{c+1}\delta e.
\end{equation}

\subsection{} \label{main-hyp}
By \eqref{conj},  $eH_{c+1}\delta e$ is an $(U_{c+1}, U_c)$-bimodule but
 the results of this paper require that it is also a progenerator
on both sides. This fails for special values of $c\in (-1,0)$
(see \cite[Remarks~3.14]{GS}) and so throughout this paper we will need to
 make the following assumption:
\begin{hypothesis}  Set
  $\mathcal{C} = \{z: z=\frac{m}{d}\ \mathrm{where}\
  m, d \in \Z\text{ with }  2\leq d\leq n \text{ and } z\notin
  \Z\}.$\label{curlyC-defn}
  Assume that  $c\in \C$ is such that $c\notin \frac{1}{2} +
\mathbb{Z}$. If $c$ is a negative rational number assume further that
$c\not\in \mathcal{C}$.
\end{hypothesis}

\noindent {\bf Remarks.} (1)
If $n=2$ and $c\in \frac{1}{2}+\NN$, then all the results of this
paper do  also hold as stated. The reason is that \cite[Theorem~3.3]{GS} is   easy
to prove directly for Cherednik algebras of type $A_1$ (see
\cite[Remark~3.14(1)]{GS})  and this is the only place where the restriction
$c\not\in \frac{1}{2}+\Z$  is used.

(2)  Using the identities $H_c\cong H_{-c}$ and $U_c\cong U^-_{-c}$ from
 \cite[Section~2]{De}, analogues of the results in this paper do hold for
  $c\in \mathcal{C}$.
   The details are similar to the arguments of
 \cite[Corollary~3.13]{GS} and  are left to the interested reader.

\subsection{} \label{balg-defn} Assume that $c$ satisfies
Hypothesis~\ref{main-hyp}.
Then \cite[Corollary~3.13]{GS} implies that the  {\it shift functors}
 \label{shift-defn}
$$S_c: \UU_c \md \to \UU_{c+1}\md: \qquad
N\mapsto  eH_{c+1}\delta e \otimes_{\UU_c} N
$$
and
 $$\widetilde{S}_c: H_{c} \md \to H_{c+1}\md: \qquad
M\mapsto  H_{c+1}\delta e \otimes_{\UU_{c}} eM$$
are  equivalences of categories.
For $i> j\geq 0$ set
   \begin{equation}\label{Mij-defn}
    B_{jj}(c) = B_{jj} = \UU_{c+j} \qquad\text{and}\qquad
   B_{ij}(c) = B_{ij} =
   (eH_{c+i}\delta e)(eH_{c+i-1}\delta e)\cdots
  (eH_{c+j+1}\delta e),
   \end{equation}
   where the multiplication
is taken in $D(\hr)\ast W$. Since each $eH_{c+i}\delta
e$ is projective, multiplication  gives a natural isomorphism
\begin{equation}\label{tpdef}
eH_{c+i}\delta e\otimes_{\UU_{c+i-1}} eH_{c+i-1}\delta e
   \otimes_{\UU_{c+i-2}}\cdots
   \otimes_{\UU_{c+j+1}} eH_{c+j+1}\delta e \cong B_{ij}.
   \end{equation}
As a result we have a Morita $\Z$-algebra\label{B-ring-defn}
 $$B = B(c) = \bigoplus_{i\geq
j\geq 0} B_{ij}$$ associated to the data $\{U_{c+i}, eH_{c+i}\delta e; i\in
\NN\}$. By setting $E = D(\hr)\ast \WW$ with
the order filtration, we see that the requirements (1--3) of
\eqref{filter-discuss} are fulfilled.

\subsection{The Hilbert scheme} We next want to relate $U_c$ and $B$ to Hilbert
schemes.
Let $\hin$\label{hin-defn} be {\it the Hilbert scheme of $n$ points on
the plane,} which we realise as the set of ideals of colength $n$
in the polynomial ring $\C[\C^2]$. If we identify the variety
$S^n\C^2$ of \textit{unordered} $n$-tuples of points in $\C^2$
with the categorical quotient $ \C^{2n}/\WW$, then the
map\label{tau-defn}
  $\widehat{\tau} : \hin \longrightarrow S^n\C^2$
which sends an ideal to its support (counted with multiplicity) is
a resolution of singularities
\cite[Theorem~1.15]{Nak}.

Consider the $\WW$-equivariant map $\h\oplus \h^* \hookrightarrow \C^{2n}$.
 To fix notation, let $\h$ be
the hypersurface $\mathbf{z}=0$ in the first copy of $\C^{n}$ and
similarly let $\h^*$ be the
hypersurface $\mathbf{z}^*=0$ in the  second copy of $\C^{n}$; thus
$\C[\C^{2n}]=\C[\h\oplus\h^*][\mathbf{z},\mathbf{z}^*]$.
Define
$$\hi=\widehat{\tau}^{\,-1} (\h\oplus\h^*/\WW). \label{hi-defn}$$
In \cite{hai3, hai1}, Haiman  proves a number of fundamental results
 about $\hin$. The analogous results  for $\hi$ will be needed in this paper
but, by \cite[Lemma~4.9 and Corollary~4.10]{GS},
they follow routinely from Haiman's work.
 In particular, $\widehat{\tau}$ restricts to a crepant
 resolution of singularities $\tau: \hi\to \h\oplus\h^*/\WW.$

\subsection{}
Write  $A^0=\cxy^\WW$, $J^0=\cxy$ and set
 $A^1
= \cxy^{\sign}$ and $ J^1 = \cxy A^1.$
\label{A-1-defn} For $k>1$, define
 $A^k=(A^1)^k$ and $J^k = (J^1)^k$ for  the
respective $k^{\text{th}}$ powers
  using   multiplication in $\cxy$. Finally, write
  \begin{equation}\label{A-alg-defn}
 A=\bigoplus_{k\geq 0}A^k\delta^k \cong A^0[A^1\delta t]
 \qquad\text{and}\qquad S = \bigoplus_{k\geq 0}J^k\delta^k
 \cong \cxy[J^1\delta t]
 \end{equation}
 for the corresponding Rees rings at the ideals $A^1\delta$, respectively
 $J^1\delta$.
 By \cite[Corollary~4.10]{GS}, $\hi = \prj (A).$ Moreover
$X_n=\prj(S)$ is the reduced fibre product
 \begin{equation}\label{Xn-defn}
 \begin{CD} X_n @> >> \h\oplus \h^*
\\ @V \rho VV @VVV \\ \hi @> \tau >> \h\oplus \h^*/\WW. \end{CD}
\end{equation}
 and the map $\rho$ is flat of degree $n!$.

Associated to these objects we have three important vector
bundles: \label{tauto-defn}the {\it rank $n$ tautological bundle $\BB$}
on $\hi$ whose fibre above $I\in \hi$ is
 the $n$-dimensional vector space $\C [\C^2]/I$; the {\it Procesi bundle} $\PP = \rho_{\ast}
\OO_{X_n}$\label{PP-defn}
 of rank $n!$; and the line bundle $\LL= \wedge^n \BB$
 \label{LLL-defn}. By \cite[Proposition~2.12]{haidis} $\LL \cong
 \OO_{\hi}(1)$, the ample line bundle associated to the isomorphism $\hi = \prj
 (A)$.

\subsection{}\label{good-P-subsec}
The given definitions of $A$ and $S$ agree with
those used in Section~6 of \cite{GS}, but not with those used in
the earlier sections of that paper where the (isomorphic) rings
$A^0[A^1t]$ and $\cxy [J^1t]$ are used. The reason for caring about the  distinction
is that  $\WW$ acts on $S$ via its natural
action on $\cxy$, and on $\PP = \rho_{\ast} \OO_{X_n}$ by the
permutation action on $X_n$.  By
\cite[(4.6) and Corollary~4.10]{GS}
 the induced  actions  agree under our chosen action:
\begin{equation}\label{good-P}
\mathrm{H}^0(\hi, \PP\otimes \LL^d) \cong J^d\delta^d
\qquad\mathrm{as}\ \WW\mathrm{-representations.}
\end{equation}
In contrast,   $\WW$ acts on
$\delta$ with the $\sign$ representation and so acts on $J^d$ by $\sign^{\otimes
d}$ \cite[loc.cit.]{GS}.

\subsection{}
\label{OMT} Let $\widehat{A}= \bigoplus_{i\geq j \geq 0}
A^{i-j}\delta^{i-j}$ \label{Aij-defn} be the $\Z$-algebra
associated to $A$ as in \eqref{zalgex1}. By
\cite[Theorem~6.4(2)]{GS}
\begin{equation}\label{OMT1}\ogr B =e \widehat{A} e \cong \widehat{A}
\end{equation}
 whenever $c$ satisfies
Hypothesis~\ref{main-hyp}. Hence
requirement (4) of \eqref{filter-discuss} is also satisfied. Combined with
\eqref{shift-defn} this implies that:
\begin{itemize}\item[]{\it If $c\in \C$ satisfies
 Hypothesis~\ref{main-hyp}, then Theorem~\ref{intro-1.2} holds
 for $U_c$ and, moreover,
 we can apply all the results from Section~\ref{zalg}.}
\end{itemize}

\subsection{Partitions}\label{dominance} We will use the same conventions for
partitions and $\WW$-representations as was used in \cite{GS}.
Thus $\irr{\WW}$\label{irred-defn}
denotes the set of irreducible representations of $\WW$, up to isomorphism.
These  irreducible representations   will be parametrised
by partitions $\mu = (\mu_1\geq \mu_2 \geq \cdots \geq \mu_\ell > 0)$
of $n$,
where the understanding is that $\mu_i =0$ for
$i>\ell$. The \textit{Ferrers diagram} of $\mu$ is the set of lattice
points\label{d-mu-defn}
 $$d(\mu) = \{ (i,j)\in \NN\times \NN : j <
\mu_{i+1}\}.$$ Following the French style, the diagram is drawn
with the $i$-axis vertical and the $j$-axis horizontal, so the
parts of $\mu$ are the lengths of the rows,
 and $(0,0)$ is the lower left corner. The \textit{arm} $a(x)$ and the
\textit{leg} $l(x)$ of a point $x\in d(\mu)$ denote the number of
points strictly to the right of $x$ and above $x$, respectively.
See \cite[(2.6.1)]{GS} for a  typical example. The {\it transpose
 partition $\mu^t$} is obtained from $\mu$ by
reflecting the
 Ferrers diagram about the line $y=x$, in other words exchanging the rows
 and columns of $\mu$.

The  {\it partition statistic} of a partition $\mu$ is
$ n(\mu)= \sum_{i} \mu_i (i-1).$
We will always use the {\it dominance
ordering}\label{dominance-defn} of partitions as in, for example
\cite[p.7]{MacD};
 thus if $\lambda$ and $\mu$
  are partitions of $n$ then $\lambda\geq \mu$ if and only if
$\sum_{i=1}^k \lambda_i \geq \sum_{i=1}^k\mu_i$
 for all $k\geq 1$.  Thus, as in \cite[Example~1, p.116]{MacD},
the  {\it trivial representation $\triv$} \label{triv-defn}
 is labelled by $(n)$ while the {sign representation} $\sign$
 is parametrised by $(1^n)$ and so $\triv>\sign$.

\subsection{}
For a pair of partitions $\mu,\lambda$ let $K_{\mu\lambda}(t,s)$
 \label{kostka-defn} be {\it the Kostka-Macdonald coefficients}
defined in \cite[VI, (8.11)]{MacD}. Their
specialisations yield the {\it classical Kostka numbers} $K_{\mu\lambda}
= K_{\mu\lambda}(0,1)$. By definition, these numbers give the
transition matrix between the bases of the ring of symmetric
functions given by the Schur functions, $s_{\mu}$, and by the
monomial symmetric functions, $m_{\mu}$; thus
 \label{kos-tran}
$ s_{\mu } = \sum_{\lambda \vdash n} K_{\mu\lambda} m_{\lambda}. $
By \cite[I, (6.5)]{MacD}, $K_{\mu \mu} = 1$ and $K_{\mu\lambda} =
0$
 unless $\mu \geq \lambda$.

\subsection{The punctual Hilbert scheme}\label{punctual-sect}
Let ${\bf 0}\in S^n\C^2$ be the zero orbit. The \textit{punctual
Hilbert scheme}, $\hio$, is defined to be $\widehat{\tau}^{\,-1}({\bf 0})$, the
 fibre of the resolution $\widehat{\tau}: \hin \rightarrow S^n\CC^2$
 above ${\bf 0}$. By \cite[Proposition~2.10]{haidis} $\hio$ is
 reduced. Since ${\bf 0}\in \h\oplus \h^*/\WW$ we also
have $\hio = \tau^{-1}({\bf 0}).$
Let $\mathfrak{m} = \cxy^\WW_+ \lhd \cxy^\WW$
 be the maximal ideal corresponding to
${\bf 0}$. Since $\hi = \prj A$, we have
$ \hio =\prj (A/A\mathfrak{m})$.

\begin{lem} \label{punctcoh}
  Let $\PP$ and $\LL$ be the bundles defined in \eqref{LLL-defn}.
Then,  for all $k\geq 0$, there is a vector space isomorphism
$\mathrm{H}^0(\hio,\, \LL^k)\cong A^k\delta^k/A^k\delta^k\mathfrak{m}$
and a $\WW$-equivariant isomorphism of vector spaces
\begin{equation}\label{punctcoh-eq}
 \mathrm{H}^0(\hio,\, \PP\otimes \LL^k) \cong
 {J^k\delta^k}/{J^k\delta^k\mathfrak{m}}.
\end{equation}
\end{lem}

\begin{proof} Once we set up the notation, this is an easy consequence of
\cite[Theorem~2.2]{hai1}.  We first prove \eqref{punctcoh-eq}.
Let $\PP_1$, $\BB_1$ and $\LL_1=\bigwedge^n \BB_1$
denote the  vector bundles on $\hin$ analogous to
$\PP$, $\BB$  and $\LL$ on $\hi$, as in \cite[(4.5)]{GS}. By
construction, the restriction of $\PP_1$ (respectively $\LL_1$) to
$\hio$ equals the restriction of $\PP$ (respectively $\LL$). Set
$\ell =kn$ and define $R(n,\ell)=\mathrm{H}^0(\hin, \PP_1\otimes
\BB_1^\ell),$ where $\BB_1^\ell=\BB_1^{\otimes \ell}$. Finally,
let $\WW^k \subset  \mathfrak{S}_\ell$ act on $\BB_1^\ell$ by
permutations.

Let $\mathfrak{m}_1 \lhd \C[\C^{2n}]$
be the maximal ideal corresponding to
${\bf 0}$.
It follows from \cite[Theorem~2.2]{hai1}  that
$R(n,\ell)/R(n,\ell)\mathfrak{m_1}
=\mathrm{H}^0(\hio,\, \PP_1\otimes \BB_1^\ell)$ while, by definition,
 if $\ep_k$ denotes the sign
representation of $W^k$ then $\LL_1^k = (\BB^\ell_1)^{\ep_k}$.
As the action of $\WW^k$ is trivial on  $\PP_1$, we therefore obtain
 $$\mathrm{H}^0(\hio,\,\PP\otimes \LL^k)
=\mathrm{H}^0(\hio,\,\PP_1\otimes \LL_1^k)
 = \mathrm{H}^0(\hio,\, (\PP_1\otimes \BB_1^\ell)^{\ep_k})
= \mathrm{H}^0(\hio,\, \PP_1\otimes \BB_1^\ell)^{\ep_k}
=R(n,\ell)^{\ep_k}/R(n,\ell)^{\ep_k}\mathfrak{m_1}.$$
By \cite[Lemma~4.9(1) and (4.6.2)]{GS}
$R(n,\ell)^{\ep_k}/R(n,\ell)^{\ep_k}\mathfrak{m_1}\cong
 J^k\delta^k/J^k\delta^k\mathfrak{m} $
 $\WW$-equivariantly. Combined with the last
  display, this  proves \eqref{punctcoh-eq}.

By \cite[Lemmas~4.4(1) and 4.9(1)]{GS},
 $A^k \delta^k$ equals the $\triv$-isotypic component of
  $J^k\delta^k$ , while  $\OO_{\hi} = \PP^\WW$.
Substituting these observations into  \eqref{punctcoh-eq} gives
$$\mathrm{H}^0(\hio,\, \LL^k)
\cong \mathrm{H}^0(\hio,\, (\PP\otimes \LL^k))^\WW\cong
{A^k}\delta^k/{A^k\delta^k\mathfrak{m}},$$ as required.
\end{proof}

\subsection{Torus action} \label{toract1}
In Section~\ref{var-sect} we will need two refinements on
 Lemma~\ref{punctcoh}: first that the isomorphisms are bigraded
 under the appropriate torus action and second that
 we can determine the dimensions of $\mathrm{H}^0(\hio,\LL^k)$
 and $\mathrm{H}^0(\hio,\PP\otimes \LL^k)$.  Both results are
 implicit in Haiman's papers, but take  a while to explain.

 We begin with the torus action.
The torus $\TT^2=(\C^*)^2$ acts linearly on $\C^2$
\label{bigrad-defn0}  as the group of diagonal matrices of the form
$\tau_{s,t}=\mathrm{diag}\{s^{-1},t^{-1}\}$; thus, as in
\cite[(12)]{hai1},  $\TT^2$ acts on $\C[\C^2] = \C[x,y]$ by
$\tau_{s,t} x = sx, \ \tau_{s,t}y = ty.$ This induces a
$\TT^2$-action on $S^n\C^2$ and $\hin$; furthermore, since
$\h\oplus \h^*/\WW \subset S^n\C^2$ is $\TT^2$-stable, there is an
induced $\TT^2$-action on $\hi$. This action further  restricts to
the punctual Hilbert scheme $\hio$ since the zero orbit in
$\h\oplus \h^*/\WW$ is invariant under~$\TT^2$.

As in   \cite[Section~2, p.377]{hai1},   there is an induced
$\TT^2$-equivariant structure on the
tautological bundle $\BB$ and on the Procesi bundle $\PP$
 arising  from the action of $\TT^2$ on $\C[\C^2]$, respectively
 $\C[\h\oplus \h^*]$. Consequently,
$\LL = \wedge^n \BB$ is also a $\TT^2$-equivariant sheaf.
Analogously, the vector bundles $\BB_1, \PP_1, \LL_1$ over $\hin$
introduced in the proof of Lemma~\ref{punctcoh} have natural
$\TT^2$-equivariant structures.

Of course, one also obtains induced actions
of $\TT^2$ on the sections of each of these bundles. This
action can be equivalently described by a $\Z^2$-grading:
 an element $f$ is homogeneous of weight $(i,j)$ if
 $\tau_{s,t} f = s^it^j f$. If $M = \bigoplus_{i,j} M^{ij}$ is a bigraded
decomposition of a module $M$ arising from such a $\TT^2$-action,
then the {\it Poincar\'{e} series} of $M$ is the Laurent series
$p(M,s,t) = \sum_{i,j} s^it^j \dim M^{ij}.$  \label{bi-poincare}
The reader should note that the variables $s$ and $t$ appear in the opposite
order in \cite[(46)]{hai1}.

\subsection{} The action of $\TT^2$ on $\C[\h\oplus \h^*]$ also  induces
a bigrading on $A^k$ and $J^k$ and the element $\delta$ is
bihomogeneous for this action with  weight $(N,0)$, where
$N=n(n-1)/2$. We have the following refinement of
 Lemma~\ref{punctcoh}.

\begin{cor} \label{punctcohref}
The isomorphisms of Lemma~\ref{punctcoh} restrict to
identifications of bigraded components
$$
\mathrm{H}^0(\hio , \LL^k)^{ij} \cong \left(
\frac{A^k\delta^k}{A^k\delta^k\mathfrak{m}}\right)^{i+Nk,j}
\qquad\mathrm{and}\qquad
 \mathrm{H}^0(\hio , \PP\otimes \LL^k)^{ij} \cong \left(
\frac{J^k\delta^k}{J^k\delta^k\mathfrak{m}}\right)^{i+Nk,j}.
$$
\end{cor}

\begin{proof}
Set $\ell=nk$, for some $k$, and let $\JJJ^k$ denote the $\C[\C^{2n}]$
analogue of $J^k$, as in \cite[(4.3)]{GS} . By \cite[(4.6.1)]{GS},
there is an isomorphism
\begin{equation}\label{pun-eq1}
\mathbb{J}^k \cong R(n, \ell)^{\epsilon_k} \cong
 \mathrm{H}^0( \hin, \PP_1\otimes \LL_1^k).
 \end{equation}
We claim that \eqref{pun-eq1}  is $\TT^2$-equivariant. For the second
isomorphism this is just \cite[(69)]{hai1}.  Note, here, that the
extra variables $a_i$ and $b_j$ used in the definition of
 $R(n,\ell)$ in  \cite[(17)]{hai1} have bidegree $(1,0)$,
 respectively $(0,1)$.   The first isomorphism in
 \eqref{pun-eq1}  is constructed explicitly in the third
 paragraph of  the proof of
 \cite[Proposition~4.11.1]{hai3}
and is clearly $\TT^2$-equivariant, as claimed.

As   in the proof of  Lemma~\ref{punctcoh},
${J}^k/{J}^k\mathfrak{m}  = \JJJ^k/\JJJ^k\mathfrak{m}_1
\cong \mathrm{H}^0(\hio ,\PP\otimes \LL^k)$ and, by
\eqref{pun-eq1}, this isomorphism is $\TT^2$-equivariant.
Since $\delta$ is bihomogeneous of degree $(N,0)$,  this
proves the second isomorphism of the
lemma. As in the proof of Lemma~\ref{punctcoh}, the first
equation follows by taking $\WW$-fixed points.
\end{proof}

\subsection{} \label{fddim}  For the next result we need some notation.
By  \cite[Proposition~3.1]{hai1}
 the $\TT^2$-fixed points on $\hi$ are labelled by partitions of
$n$:   to any partition $\eta$ we associate the monomial  ideal
 \begin{equation}\label{I-eta-defn}
 I_{\eta} = \C \cdot \{ x^ry^s : (r,s)\notin
d(\eta)\} \unlhd \C [x,y],
\end{equation}
 regarded as a point of $\hi$.
 We also  recall
that the {\it Catalan number} $C_n^{(k)}$
 is defined to be $$C_n^{(k)}\ = \ \frac{1}{kn+1}\binom{(k+1)n}{n}.$$

\begin{lem}
 For all   $k\geq 1$, we have \
 $\dim \mathrm{H}^0(\hio ,\, \PP \otimes \LL^{k-1}) = (kn+1)^{n-1}$
and  $\dim \mathrm{H}^0(\hio ,\, \LL^k) =C_n^{(k)}.$
\end{lem}

\begin{proof}
 As in the proof of Lemma~\ref{punctcoh}, we may replace $\PP$ and
$\LL$ by $\PP_1$ and $\LL_1$ and so, by \cite[Theorem~2.2]{hai1}, the hypothesis of
 \cite[Theorem~2]{haidis} is now valid. Thus, by
  \cite[(1.7) and Theorem~2]{haidis}, and in the notation of that paper,
$\dim \mathrm{H}^0(\hio ,\, \LL_1^k)= C_n^{(k)}(1,1)=C_n^{(k)}$.
This  proves the second formula.

  For a partition $\mu$ of $n$, set
$\Omega(\mu)= \prod_{x\in
d(\mu)}(1-s^{1+l(x)}t^{-a(x)})(1-s^{-l(x)}t^{1+a(x)})$, and write
$$B_{\mu} =  \sum_{(x,y)\in d(\mu)} s^{x} t^{y}, \qquad
\Pi_{\mu} = \prod_{ {(x,y)\in d(\mu)}\atop{(x,y)\neq (0,0)}}
(1-s^x t^y)\qquad \mathrm{and}\qquad P_{\mu} = \sum_{\lambda}
s^{n(\mu)}K_{\lambda \mu}(t,s^{-1}) \dim \lambda.$$

We claim that $\mathrm{H}^0(\hio, \PP_1 \otimes\LL_1^{k-1})$ has
the  Poincar\'{e} series
\begin{equation} \label{bigradedactchar}
 p( \mathrm{H}^0(\hio, \PP_1 \otimes
\LL_1^{k-1}), s,t) = \sum_{\mu}
P_{\mu}s^{(k-1)n(\mu)}t^{(k-1)n(\mu^t)}(1-s)(1-t)
\Pi_{\mu}B_{\mu}\Omega(\mu)^{-1}.
\end{equation}
 We prove the claim by following the derivation of
\cite[Theorem~3.2]{hai1}.  A locally free, $\TT^2$-equivariant
resolution $V_{\bullet} \to \mathcal{O}_{\hio}$ of
$\OO_{\hin}$-modules is given by \cite[Proposition~2.2]{hai1}.
 Since $\LL_1^{k-1}=(\bigwedge^n
\BB_1)^{\otimes(k-1)}$
 is a direct summand   of $\BB_1^{\otimes n(k-1)}$,
 \cite[Theorem~2.2]{hai1} shows that
 $  \mathrm{H}^i(\hin, \PP_1 \otimes
\LL_1^{k-1}\otimes \OO_{\hio})=
\mathrm{H}^i(\hio, \PP_1\otimes \LL_1^{k-1})
= 0$ for all $i >0$.
 Combining these observations means that
\begin{eqnarray}
\label{ABL1} p\left( \mathrm{H}^0(\hio, \PP_1 \otimes
\LL_1^{k-1}), s,t \right) &=& p\left( \mathrm{H}^0(\hin, \PP_1 \otimes
\LL_1^{k-1}\otimes \OO_{\hio}), s,t \right)\notag
 \\ &=&
\sum_{j\geq 0} (-1)^j p\left( \mathrm{H}^0(\hin, \PP_1 \otimes
\LL_1^{k-1}\otimes V_j), s,t \right). \notag
\end{eqnarray}
By \cite[(35)]{hai1}, the sheaves $ \PP_1 \otimes\LL_1^{k-1}\otimes V_j$
are acyclic for the global section functor.  Thus,
combining the last displayed equation  with
 the Atiyah-Bott-Lefschetz formula \cite[Proposition~3.2]{hai1}
 gives
\begin{equation}\label{ABL2}
\begin{array}{rl}\displaystyle
 p\left( \mathrm{H}^0(\hio, \PP_1 \otimes
\LL_1^{k-1}), s,t\right) \ =&\displaystyle
 \sum_{\mu} \sum_{j\geq 0} (-1)^j
 p\left( \PP_1(I_{\mu}) \otimes \LL_1^{k-1}(I_{\mu})
 \otimes V_j(I_{\mu}), s,t \right)\Omega(\mu)^{-1} \\
 \noalign{\vskip 7pt}
 =& \displaystyle
  \sum_{\mu} p\left(\LL_1(I_{\mu}),s,t \right)^{k-1}\,\Omega(\mu)^{-1}  \,
 \sum_{j\geq 0} (-1)^j  p\left( \PP_1(I_{\mu})
 \otimes V_j(I_{\mu}), s,t \right).
\end{array}
\end{equation}
As was shown in the proof of \cite[Proposition~4.8]{GS}, $
p(\LL_1(I_{\mu}), s,t) = s^{n(\mu)}t^{n(\mu^t)}, $
 while the final sum in \eqref{ABL2} is
 computed in  \cite[(89)]{hai1} (the individual terms in \cite[(89)]{hai1} are
 defined in Equations (71), (73) and (88) of \cite{hai1};
 note, however, that the convention \cite[(46)]{hai1} for Poincar\'e series
 differs from ours). These computations combine to give
\eqref{bigradedactchar}, as claimed.

By \cite[(1.9)]{haidis} we have $ \sum_{x\in d(\mu)}l(x) = n(\mu)
$ and, similarly, $\sum_{x\in d(\mu)} a(x) = n(\mu^t)$.
 Applying these formul\ae\ to \eqref{bigradedactchar} gives
\begin{equation}
\label{haifd} p( \mathrm{H}^0(\hio, \PP_1 \otimes \LL_1^{k-1}), s,t)
\ = \   \sum_{\mu} \frac{P_{\mu}(1-s)(1-t)s^{kn(\mu)}t^{kn(\mu^t)}
 \Pi_{\mu}B_{\mu}}{\strut\prod_{x\in d(\mu)}
  (t^{a(x)}-s^{1+l(x)})(s^{l(x)}-t^{1+a(x)})}.
\end{equation}
The right hand side of this equation is the same as
\cite[(3.27)]{garhai2},
provided one  specialises the formal
variable $X=(x_1,x_2,\ldots )$ from \cite[(3.27)]{garhai2} to
$x_1=\cdots =x_n = 1$ and $x_i=0$ for $i> n$. Now
\cite[(3.27)]{garhai2}  can be computed from
\cite[Theorem~4.1]{garhai2}. Let $e_n$ denote the
$n^{\mathrm{th}}$ elementary symmetric function and note that, if
$q=1$ and $X$ is as above, then  $e_n[X(1+q+\ldots + q^{kn})]=e_n
[X(kn+1)] =(kn+1)^{n}$ in the notation of
\cite[Introduction]{garhai2}. Therefore, specialising
\cite[Theorem~4.1]{garhai2} to $q=1$ (which means taking $s=t=1$
in our notation)
 gives
$$\dim \mathrm{H}^0(\hio, \PP_1 \otimes \LL_1^{k-1}) \ =\
\frac{1}{kn+1}e_n[X(1+q+\ldots + q^{kn})]  =(kn+1)^{n-1},$$
as required.
\end{proof}

\subsection{}
\label{bichar2} We note for use in \eqref{bigrfd} that the
derivation of \eqref{bigradedactchar} also yields
\begin{equation}
\label{bigradedactchar2} p( \mathrm{H}^0(\hio, \LL^{k}), s,t) =
 \sum_{\mu} s^{kn(\mu)}t^{kn(\mu^t)}(1-s)(1-t)\Pi_{\mu}B_{\mu}
 \Omega(\mu)^{-1}.
\end{equation}
This is the two variable Catalan-like number $C_n^{(k)}(s,t)$ of
Garsia  and Haiman, as defined in  \cite[(1.10)]{haidis}.

\section{Representations of $H_c$ and coherent sheaves on $\hi$}
\label{cohsh}

\subsection{} Throughout the section we assume that $c\in\C$ satisfies
Hypothesis~\ref{main-hyp} although, by the observations  from
Remark~\ref{main-hyp}(2),
analogues of the results proved here do hold for more general
values of $c$.

The results of the previous sections allow us to associate a
coherent sheaf $\Phi_{\Lambda}(M)$ on $\hi$ to any $\UU_c$- or
$H_c$-module  $M$ with a good filtration $\Lambda$ and in this
section we investigate the general consequences of this
construction. These results will begin to explain one of the
central theses of this paper:  $\Phi_{\Lambda} (M)$ carries much
deeper information about the structure of $M$ than does, for
example, the associated graded $A^0$-module $\gr_\Lambda M$. One
reason for this is that the fundamental Morita equivalence between
$U_c$ and $U_{c+1}$, and with it the more subtle information about
$U_c$-modules, is lost in passing to the associated graded rings
but survives in our approach---in the notation of \eqref{A-1-defn}
it becomes the shift functor $(\LL\otimes -)$ for coherent
sheaves (see Lemma~\ref{tensorline}).  In
Proposition~\ref{compBKR}, we show that there are intriguing
connections between the map $\Phi$ and the Bridgeland-King-Reid
equivalences.

\subsection{Associated sheaves} \label{cohehe}
Given a filtered noetherian ring $R$, write
$R\fMd$  for the category of $R$-modules equipped with increasing
exhaustive filtrations and whose morphisms respect those filtrations. Let
$R\fmd$ be the full subcategory of $R\fMd$ consisting of modules
with good filtrations. By  \eqref{zalgex1} and  \eqref{OMT1},
graded modules over $\ogr B\cong \widehat{A}$ are the same as those over $A$.
In the notation of \eqref{graded-section} we therefore
have a functor\label{Phi-defn}
$$\Phi: \UU_c\fMd\ \longrightarrow\ \Qcoh \hi
\qquad (M,\Lambda) \mapsto \Phi_{\Lambda}(M) =
\pi(\gr_{\Lambda} \widetilde{M}), $$
 where $\pi$ is the quotient map from $A\lGr$ to $A\lQgr\simeq \Qcoh \hi$.
 Lemma~\ref{graded-lemma}(1) implies that  $\Phi$ restricts to a functor
$\Phi : \UU_c\fmd \rightarrow \coh \hi.$

We will often abuse notation by writing just $\Phi(M)$
when the filtration on $M$ is clearly understood, but we will write
$\Phi_{\Lambda}(M)=\Phi^{c}_{\Lambda}(M)$ if we need to specify the initial
algebra $U_c$.

Similarly, we have a functor $\widehat{\Phi}: H_c\fMd \to
 \Qcoh\hi$.   More precisely, if
$(N,\Gamma)\in  H_c\fMd$  then applying the functor $E_c=(eH_c\otimes -)$ to
$(N,\Gamma)$ gives an object
$(eN , \Gamma_E)\in \UU_c\fMd$
by setting $\Gamma_E^k(eN) = e \Gamma^k(N)$. Thus we can extend $\Phi$ to
\label{Phi-hat-defn}
$$\widehat{\Phi}: H_c\fMd \longrightarrow \Qcoh \hi  \qquad
(N,\Gamma) \mapsto \Phi_{\Lambda}(M) \quad\mathrm{for}\
(M, \Lambda) = (eN, \Gamma_E) .$$
If $(N,\Gamma)\in H_c\fmd$, then
 $\gr_\Lambda eN = \bigoplus e\Gamma^n N/e\Gamma^{n-1}N$ is a homomorphic
image of the finitely generated $\cxy^\WW$-module $\gr_\Gamma N$ and so
$\Lambda$ is also a good filtration.
Thus, as before, it follows from Lemma~\ref{graded-lemma} that $\widehat{\Phi}$
 restricts to a functor
$\widehat{\Phi}: H_c\fmd \rightarrow \coh \hi.$

\subsection{Warning}\label{warning-subsec}
 The careful reader has noticed that we have two
filtrations in play at the moment. In \eqref{OMT1} we used the order
filtration on the bimodules $B_{ij}$, whilst in \eqref{indfil} and
\eqref{cohehe} we have
insisted on the tensor product filtration for modules.
For   general modules these two
filtrations can differ (an example is given in \cite[(7.4)]{GS})  and so
one can ask whether $\Phi (U_c)$ actually
equals $\OO_{\hi}$. It does, by \cite[Lemma~7.2]{GS}.

This potential ambiguity   occurs several times in this paper---one
instance occurs in the next lemma---but in each case it is resolved by
\cite[Lemma~7.2]{GS}.

\subsection{Shift Functors} \label{tensorline}
 Recall from \eqref{A-1-defn}  that  $\LL \cong \OO_{\hi}(1)$,
 the Serre twisting sheaf corresponding to the Rees ring $A$. The
role of the  shift functor $(\LL\otimes -)$ in $\Qcoh \hi$
 is played by the shift functor
$S_c$ in the noncommutative world:

\begin{lem} Assume that $c\in \C$ satisfies Hypothesis~\ref{main-hyp} and
let $(M, \Lambda) \in \UU_c\fMd$. Define
$(M', \Lambda') \in \UU_{c+1}\fMd$ by setting
 $M'=S_{c}(M)=(eH_{c+1} \delta e)\otimes_{U_c}M$
 together with the tensor product filtration $\Lambda'$ from by
\eqref{indfil}. Then
$\Phi^{c+1}_{\Lambda '}(M') \cong \Phi^c_{\Lambda}(M)\otimes \mathcal{L}$.
\end{lem}

\begin{proof} Recall from \eqref{Mij-defn} the definition of
$B_{ij}(c)$. By construction  $M(i) = B_{i0}(c)\otimes_{\UU_c} M$
and
\begin{equation}\label{proofample}
 M'(i) = B_{i0}(c+1)\otimes_{\UU_{c+1}} M'
= B_{i+1, 0}(c)\otimes_{\UU_c} M = M(i+1).\end{equation}
As in \eqref{warning-subsec} we should be careful to compare the two filtrations
on $M'(i)$ given by this equality. By definition and \cite[Lemma~7.2]{GS},
 $\Lambda'$ is given by
\begin{equation}\label{proofample2}
\begin{array}{rl}
(\Lambda')^k M'(i) \ = &
\sum_{u+v+w= k} \ord^u B_{i0}(c+1)\otimes \ord^vB_{10}(c)\otimes \Lambda^wM
\\ \noalign{\vskip 4pt}
& \qquad =\ \sum_{x+w=k} \ord^x B_{i+1,0}(c)\otimes \Lambda^wM
\ = \ \Lambda^k M(i+1).
\end{array}
\end{equation}
 Thus \eqref{proofample} is also an equality of filtered modules.

The sheaf $\Phi_{\Lambda}(M)$ on $\hi$ corresponds
to the graded $A$-module $\gr_{\Lambda} \widetilde{M} =
 \bigoplus_{i\geq 0} \gr_{\Lambda} M(i)$. The action of $A^j$ taking
 $\gr_{\Lambda} M(i+1)$ to
$\gr_{\Lambda} M(i+j+1)$ is induced from the mapping $$B_{i+j+1,
i+1}(c)\otimes_{\UU_{c+i+1}} M(i+1) \longrightarrow M(i+j+1).$$
Similarly the action of $A^j$ on $\gr_{\Lambda'} M'(i)$ is induced from
the map
$B_{i+j, i}(c+1)\otimes_{\UU_{c+i+1}} M'(i)\rightarrow M'(i+j).$
Since $B_{i+j+1,i+1}(c)=B_{i+j, i}(c+1)$, it follows from
 \eqref{proofample}  that  these actions are equal for all $j$ and all
$i\geq 0$. Thus,   \eqref{proofample2} implies that
  $(\gr_{\Lambda}\widetilde{M})[1]=\bigoplus_{i\geq 0} {\gr_{\Lambda} M(i+1)}$
 and
$\gr_{\Lambda'}\widetilde{M'}=\bigoplus_{i\geq 0} \gr_{\Lambda'}M'(i)$
 are equal  as  graded $A$-modules. Since  $\LL \cong \OO_{\hi}(1)$,
this is equivalent to the conclusion of the lemma.
 \end{proof}

\subsection{} \label{cohpro}
Many of the results of this paper can be interpreted as saying that important
$H_c$-modules correspond via $\widehat{\Phi}$ to important sheaves on $\hi$. We
give one example of this philosophy here---for the module $H_c$ itself.
Comparing this result with
 \eqref{Xn-defn} shows that one can think of $H_c$ as a noncommutative analogue
 of the isospectral Hilbert scheme $X_n$.

\begin{thm} Assume that $c\in \C$ satisfies Hypothesis~\ref{main-hyp}. Then
  $\ogr (\widetilde{eH_c}) =\bigoplus_{k\geq 0} eJ^k\delta^k$ and
   $\widehat{\Phi}_{\ord} ( H_c) \cong \PP$, the Procesi bundle on
$\hi$.
 \end{thm}

\begin{proof}  We emphasise that if we begin with $(N,\Lambda) = (H_c, \ord)$,
then $\widehat{\Phi}(H_c)$ is defined by the tensor product
filtrations on the summands $B_{k0}\otimes_{U_c} eH_c$ of
$\widetilde{eH_c} =\bigoplus B_{k0}\otimes_{U_c} eH_c \cong
\bigoplus N(k)$ rather than the $\ord$ filtration. However, by
\cite[Lemma~7.2]{GS}, these filtrations are  equal. So the
identity  $\ogr (\widetilde{eH_c}) = \bigoplus eJ^k\delta^k$ of
graded vector spaces is immediate from \cite[Proposition~6.5]{GS}.
That this is an identity of $A$-modules then follows from
\cite[Lemma~6.7]{GS}. Since $X_n=\prj{S}$ for $S=\bigoplus
J^k\delta^k$, in the notation of \eqref{A-alg-defn}, this  implies
that $\widehat{\Phi}(H_c) \cong  \pi(S)=\PP$, as required.
\end{proof}

\subsection{}\label{tensor-ann-sec} Recall from \eqref{charvar-defn} the
definitions of the  the characteristic  variety $\Chr M \subseteq \hi$
  and the associated variety
   $\Chro M  \subseteq \h\oplus\h^*/\WW$ for $(M, \Lambda)\in U_ c\fmd$.
As has been mentioned in the introduction,  $\Chr M$ carries  more
subtle information $M$  than $\Chro M$. Nevertheless these two
varieties are connected and we will show in
Proposition~\ref{compare-char} that $\Chro M = \tau(\Chr M)$,
under the  resolution of singularities $\tau$. In a weak sense
this shows that the natural  diagram
$$\begin{CD} B\qgr @>\sim >> U_c\md \\  @V \gr VV   @VV\gr V \\
 \coh \hi @>\tau_{\ast} >> A^0\md \end{CD}$$
 commutes.   One can show by example that the stronger commutativity result,
  $\tau_{\ast}(\gr_\Lambda(\widetilde{M}))=\gr_\Lambda(M)$, does not always hold.

We begin with an abstract lemma.

\begin{lem}\label{tensor-ann} Let $R=\bigcup F^iR$ and $S=\bigcup F^iS$ be
subrings of a filtered ring $D=\bigcup_{i\geq 0}F^iD$ and $T$ an
$(S,R)$-sub-bimodule of $D$. Assume that $\gr_FD$ is commutative, that
 $\gr_FR=\gr_FS$ under the induced filtrations and  that the induced
filtration $F$ on $T$ is good on both sides.
Let $(N, F)\in R\fmd$ and give $T\otimes_RN$ the tensor product filtration $F$.
Then  $I=\mathrm{ann}_{\gr_FR}(\gr_FN)\subseteq
\mathrm{ann}_{\gr_FS}(\gr_F(T\otimes N))$.  \end{lem}

\begin{proof}  We drop the subscript $F$ from  all graded modules.
The commutativity of $\gr D$ implies the following: take
$\overline{q}\in \gr^jT$ and $\overline{r}\in \gr^iR=\gr^iS$, and
lift $\bar{q}$ to
 $q\in F^jT$ and lift $\bar{r}$ to elements $r_1\in
F^iS$ and $r_2\in F^iR$. Then, $r_1q=qr_2$ modulo $F^{i+j-1}T$.

Recall from  \eqref{princ-symbol} that the principal symbol of $n\in N$ is
written $\sigma(n)$.
Now let $\overline{a}\in I_k=I\cap\gr^kR$ and pick a typical generator,
say $\sigma(q\otimes n)$ of
$\gr^t(T\otimes N)$, where  $q\in F^iT$ and $n\in
F^{t-i}N$. Lift $\overline{a}$ to elements $a\in F^kS$ and   $b\in F^kR$. Since
$\sigma(b)\sigma(n)=0$, clearly $bn\in F^{k+(t-i-1)}(N)$. By the last
paragraph, $aq-qb=c\in F^{k+i-1}T$ and so
$$a(q\otimes n) = qb\otimes n + c\otimes n
=q\otimes bn + c\otimes n \in F^{k+t-1}(T\otimes N).$$
 In other words, $\bar{a}\sigma(q\otimes n)=0$ and hence $\bar{a} \gr(T\otimes
 N)=0$.
\end{proof}

\subsection{Corollary}\label{tensor-ann2}  {\it
Let $R_\Z$ be a Morita $\Z$-algebra satisfying Hypotheses (1--4)
of \eqref{filter-discuss} and suppose that  $M$ is a finitely
generated left $R_0$-module. Then $\Chro  M = \Chro
(R_{k0}\otimes_{R_0} M)$ for all $k\geq 0$.}

 \begin{proof} Set $Q=R_{k0}$.
 Pick a good filtration $\Lambda$ of $M$ and give $Q\otimes M$
 the tensor product filtration $\Gamma$; by \cite[Lemma~6.7(2)]{GS} this is a good
 filtration. Since $\Chro M$  is independent of the choice of good
 filtration \cite[Proposition~5.1]{Bj2},
  it  is  determined
 by $\mathrm{ann}_{\gr (R_0)}(\gr_\Lambda M).$ The hypotheses of
   Lemma~\ref{tensor-ann} are satisfied by $R=R_0$, $S=R_{k}$, $T=Q$
   and $N=M$ and so that lemma
 implies that $\Chro   M \supseteq \Chro  (Q\otimes M)$.

The lemma may also be  applied  to the rings $R=R_k$, $S=R_0$ with
the modules $T=P=Q^*$
and $N=Q\otimes M$, with its good tensor
  product filtration $\Gamma$.  By  hypothesis, $Q$ is a progenerator
   and so $P\otimes_{R_{k}}N = (P\otimes_{R_{k}}
  Q)\otimes_{R_{0}} M \cong  M.$   Thus the lemma  now implies that $\Chro
   (Q\otimes M) = \Chro   N  \supseteq  \Chro (P\otimes N)= \Chro   M$
  and hence that $\Chro  M =\Chro  (Q\otimes M)$.
  \end{proof}

 \subsection{}  We can now relate $\Chro M $ to
 $ \Chr M $ for a finitely generated $U_c$-module $M$,
 where $c$ satisfies Hypothesis~\ref{main-hyp}.
We will also need this result in other contexts and so we
generalise it slightly. Let $R_\Z$ be a Morita $\Z$-algebra that
satisfies Hypotheses (1--4) of \eqref{filter-discuss}. The graded
ring $S=\bigoplus S_i$
 of \eqref{filter-discuss}(4) satisfies $S_0=\gr R_0$ and so we have
 a natural surjection $\tau':\prj S\to \spc \gr R_0$ analogous to $\tau$.

  \begin{prop}\label{compare-char}
Let  $R_\Z$  be a Morita $\Z$-algebra that satisfies Hypotheses
(1--4) of \eqref{filter-discuss}  and let $M$ be a finitely
generated left $R_0$-module. Then $\tau'(\Chr M)=\Chro M$.
 \end{prop}

\noindent  {\bf Remark. }
  As was remarked in \cite[(1.3)]{GS}, Theorem~\ref{intro-1.2}
 can be regarded as a kind of analogue for Cherednik algebras of the
 Beilinson-Bernstein equivalence of categories for enveloping algebras.
Under this analogy, the proposition for $U_c$   corresponds to
\cite[Theorem~1.9(c)]{BB}.

  \begin{proof} Identify $S_0=\gr R_0$ and set
   $I=\mathrm{ann}_{S_0}(\gr M)$
 and $J =\mathrm{ann}_{S}(\gr \widetilde{M})$.
  If $P=\bigoplus_{i\geq 0}P_i$ is a graded prime ideal of $S$  that does not
  contain the irrelevant ideal $S_+=\bigoplus_{i>0} S_i$, then the
  corresponding varieties satisfy
  $\tau'(\mathcal{V}(P))=\mathcal{V}(P_0)$.
  For an arbitrary  graded ideal of $S$ we have to worry about occurrences of
  $S_+$. Thus, let $K=\bigoplus_{i\geq 0}K_i$
  denote the largest ideal $L$ of $S$ such that
  $L(S_+)^m\subseteq J$ for some $m\geq 0$. By definition
  $\mathcal{V}(J)=\mathcal{V}(K)$ but now the prime ideals $P$ minimal over
  $K$ do not contain $S_+$. Thus $\tau'(\mathcal{V}(J))=\tau'(\mathcal{V}(K))=
  \mathcal{V}(K_0)$.

  A standard exercise shows that
  $K=\mathrm{ann}_S\big((\gr \widetilde{M})_{\geq t}\big)$, for $t \gg 0$.
   For some such $t$ pick $s$ such that $(\gr \widetilde{M})_{\geq t}$
  is generated by  $ \bigoplus_{j=t}^s \gr {M}(j)$ as an $S$-module.
  Clearly $K_0 = \bigcap_{j=t}^s \mathrm{ann}_{S_0}\gr M(j)$.
  By Corollary~\ref{tensor-ann2} $\sqrt{I} = \sqrt{\,\mathrm{ann}_{S_0}(\gr M(j))}$
 for each such $j$ and so $\sqrt{I}=\sqrt{K_0}$.
  By the conclusion of the last paragraph, this implies that
    $\tau'(\Chr M) = \Chro M.$
 \end{proof}

\subsection{}\label{pure-question}
 Let $M\in U_c\md$ be pure in the sense that
each nonzero submodule has the same Gelfand-Kirillov dimension as $M$. By
Watanabe's Theorem \cite{Wa}, $\cxy^\WW$ is Gorenstein and so
Gabber's theorem says that the associated variety $\Chro M$ is
equidimensional (see, for example, \cite[(5.2)]{Bj2} and
\cite[Theorem~7.1]{bj}). It would be interesting to know if the
analogous result holds for $\Chr$:

\begin{question} If  $M\in U_c\md$ is pure, is $\Chr M$
equidimensional?
\end{question}

\subsection{The BKR equivalence}\label{bkr}
In the final result of this section, Proposition~\ref{compBKR}, we show that
there is  an interesting connection between $\widehat{\Phi}$ and the
Bridgeland-King-Reid equivalence of categories.

We start by giving a version of that equivalence of categories that is
appropriate
for $\h\oplus\h^*/\WW$.
Recall  the following commutative diagram from \eqref{Xn-defn}.
 $$\begin{CD} X_n @>  f >> \h\oplus \h^*
\\ @V \rho VV @VVV \\ \hi @> \tau >> \h\oplus \h^*/\WW \end{CD}
$$
Let $D(\hi)$ denote the derived category of complexes of
quasicoherent sheaves on $\hi$ with bounded  coherent
cohomology and let $D^\WW(\h\oplus \h^*)$ be the derived category
of complexes of $\WW$-equivariant quasicoherent sheaves on
$\h\oplus \h^*$, again with bounded coherent cohomology.
Equivalently, $D^\WW(\h\oplus \h^*)$ can be identified with the
derived category of bounded complexes of finitely generated
$\C[\h\oplus \h^*]\ast\WW$-modules.

\begin{thm} {\rm (\cite[Theorem~1.1]{BKR} and
\cite[Corollary~5.1]{hai1})} The
functor $\Psi : D^\WW(\h\oplus\h^*)
\to D(\hi)$ given by $\mathcal{F} \mapsto (\rho_{\ast}
\circ Lf^{\ast} (\mathcal{F}))^{\WW}$\label{Psi-defn}
is an equivalence of triangulated
categories.
\end{thm}

\begin{proof}
By \cite[Corollary~4.10]{GS},  $\tau : \hi\to \h\oplus \h^*/\WW$
is a crepant resolution of singularities and $\hi$ is irreducible.
Since $\WW$ preserves the symplectic form on
$\h\oplus\h^*$, the result will therefore  follow from
\cite[Corollary~1.3]{BKR} once we show that $\hi\cong
\h\oplus\h^*\twoslash\WW$, the $\WW$-Hilbert scheme on $\h\oplus
\h^*$. By definition the  points  of $\h\oplus\h^*\twoslash\WW$
are (isomorphism classes of) $\WW$-clusters on $\h\oplus \h^* $;
that is, cyclic $\cxy\ast \WW$-modules which both carry the
regular representation as $\WW$-modules and are generated by the
trivial representation.

Let $\C^{2n}\twoslash\WW$ denote the $\WW$-Hilbert scheme on
$\C^{2n}$. Any central element of $\C[\C^{2n}]\ast \WW$ acts by
scalar multiplication on a $\WW$-cluster in $\C^{2n}$, so there is
a morphism $p: \C^{2n}\twoslash\WW \longrightarrow S^n\C^2$. Since
$\h\oplus \h^*\twoslash\WW$ consists of the clusters which are
annihilated by the elements ${\bf z}$ and ${\bf z}^*$ of
\eqref{hi-defn}, it follows that $p^{-1}(\h\oplus \h^*/\WW)
=\h\oplus \h^*\twoslash\WW$. By \cite[Theorem~5.1]{hai3} there is
an isomorphism of $S^n\C^2$-varieties $\iota: \C^{2n}\twoslash\WW
\longrightarrow \hin.$ Since $\hi = \tau^{-1}(\h\oplus \h^*/\WW)$
this induces the required isomorphism $
\h\oplus\h^*\twoslash\WW\cong  \hi.$

Finally, we compute the functor giving the equivalence of
categories, since the functor we want is not quite  that of
\cite{BKR} or \cite{hai1}. As is explained by \cite[Corollary~5.1
and Proposition~5.3]{hai1}, adjusted to our setup over
$\h\oplus\h^*$, the equivalence $D(\hi) \to D^\WW (\h\oplus\h^*)$
is given by $\alpha=R\Gamma_{X_n}(\PP\otimes -)$, with inverse
$\OO(1)\otimes (\rho_* \circ  Lf^*)^{\sign}.$ Consequently, if we
replace $\alpha$ by the functor $R\Gamma_{X_n}(\PP(-1)\otimes
-)\otimes \sign$, which is obviously still an equivalence, then
its inverse will indeed be the functor $ (\rho_{\ast} \circ
Lf^{\ast} )^{\WW}$ that we require.
\end{proof}

\subsection{} \label{compBKR}
Theorem~\ref{bkr} gives   a second  way in which to
 construct coherent sheaves on
$\hi$ from objects of $H_c\fmd$. Let $(N, \Gamma)\in H_c\fmd$. Then
 $\gr_{\Gamma} N$ is a finitely generated
$\C[\h\oplus \h^*]\ast \WW$-module, which we consider  as
an object of $D^\WW(\h\oplus \h^*)$ concentrated in degree zero.
Applying the functor of Theorem~\ref{bkr} to this module
gives a {\it complex} $\Psi(N)=\Psi_{\Gamma}(N)$ of
coherent sheaves on $\hi$. Finally,
taking the zeroth homology of $\Psi$ gives the functor
$$\Psi_\Gamma^0 : H_c\fmd \longrightarrow \coh \hi, \qquad (N, \Gamma) \mapsto
(\rho_*f^*(\gr_{\Gamma}N))^{\WW}.$$ The following proposition shows this
functor is closely related to the functor $\widehat{\Phi}$ obtained from
our $\Z$-algebra construction.

\begin{prop} For every $(N, \Gamma)\in H_c\fmd$ there is a surjective mapping
$\Psi^0_{\Gamma} (N) \twoheadrightarrow \widehat{\Phi}_{\Gamma}(N)$
of sheaves on $\hi$.
\end{prop}

\begin{proof}
To prove the proposition we will first study free
$H_c$-modules, and then pass to the general case via free $H_c$-resolutions.
The result for free modules follows from the following result.

\begin{sublemma}   Let $R_x: H_c\to H_c$ be the mapping given
by right multiplication by an element $x \in H_c$. Then $\Psi^0_{\ord}(H_c) =
\PP=\widehat{\Phi}_{\ord}(H_c)$ and the induced map
$\Psi_{\ord}^0( R_x) : \PP\to \PP$
coincides with $\widehat{\Psi}(R_x)$.
\end{sublemma}

\noindent
{\it Proof of the sublemma.} This just amounts to unravelling the definitions
and we begin with $\widehat{\Psi}$.

Recall from \eqref{warning-subsec} that the $\ord$ and tensor
product filtrations coincide on $\widetilde{eH}_c$. It is
therefore clear that the map $\ogr R_x : \ogr\, \widetilde{eH}_c
\to \ogr\, \widetilde{eH}_c$ induced from $R_x$ is just given by
right multiplication by the principal symbol $\sigma(x)\in
\cxy\ast\WW$ of $x$. Write $\sigma(x) = \sum p_\gamma \gamma$, for
some $p_\gamma\in \cxy$ and $\gamma\in \WW$
and recall from \eqref{chered-defn} that
multiplication in $\cxy  \ast \WW$ is defined by $wf=w(f)w$ for
$w\in \WW$ and $f\in \cxy$. By
Theorem~\ref{cohpro}, $\ogr \widetilde{eH}_c  \cong S=
\bigoplus_{k\geq 0} J^k\delta^k$, where the isomorphism is induced
on each summand from the natural isomorphism $\Theta^{-1} :
e(\cxy\ast\WW)\xrightarrow{\sim} \cxy$ given by $
ep\gamma=e\gamma\gamma^{-1}(p)\mapsto \gamma^{-1}(p)$, for   $p
\in \cxy$ and $\gamma\in \WW$. Therefore, at the level of graded
$A$-modules,
 the map $\widehat{\Phi}(R_x)$ is  the map
$S\to S$
 defined  by
\begin{equation} \label{noncommfirst}
j \ \mapsto\ ej \ \mapsto\ \sum_{\gamma}
ejp_{\gamma}\gamma \ \mapsto\ \sum_{\gamma} \gamma^{-1}(jp_{\gamma}),
 \qquad\mathrm{for}\ j\in J^k\delta^k.\end{equation}

Now consider $\Psi^0$. The map  $f: X_n \longrightarrow \h\oplus \h^*$
is $\WW$-equivariant, so
the pullback $f^*$ sends $\WW$-equivariant coherent sheaves on
$\h\oplus \h^*$ to $\WW$-equivariant coherent sheaves on $X_n$.
 We identify the category of  $\WW$-equivariant coherent
sheaves on  $\h\oplus\h^*$ with $\cxy\ast\WW\md$  and, by
\eqref{A-1-defn}, write  $X_n = \prj S$.  In
this language, if $M$ is an object of $\cxy\ast \WW\md$, then $f^{\ast}(M) =
\ddd\otimes_{\cxy} M$ with $\WW$ acting diagonally on the right hand side.
Similarly,  the functor $\rho_*: \coh X_n\to \coh \hi $
corresponds to the restriction functor from
graded $S$-modules to graded $A$-modules.

The map $\ogr R_x:  \ogr H_c=\cxy \ast \WW \to
\cxy\ast \WW$ is clearly given by right multiplication by $\sigma(x)$
and so, by the last paragraph,   $\Psi^0(R_x)$ is given by
  the homomorphism of graded $A$-modules
\begin{equation}
\label{imagefunctor} R_{\sigma(x)} : (\ddd\otimes_{\cxy} (\cxy \ast
\WW))^{\WW} \longrightarrow (\ddd \otimes_{\cxy} (\cxy \ast
\WW))^{\WW}\end{equation} defined by right multiplication
by $\sigma(x)$.
We need to compare this to \eqref{noncommfirst}. There is a natural identification
$$\chi:\ddd\ \buildrel{\sim}\over{\longrightarrow}\
(\ddd\otimes_{\C} \C\WW)^{\WW}
\  \buildrel{\sim}\over{\longrightarrow} \
(\ddd\otimes_{\cxy} (\cxy\ast \WW))^{\WW}$$
given by
$$ b\  \mapsto\ \frac{1}{n!} \sum_{g\in\WW} g(b)\otimes_{\C} g
\ \mapsto\  \frac{1}{n!} \sum_{g\in\WW} g(b)\otimes_{\cxy} g
\qquad \mathrm{for}\ b\in S.$$
Write $\sigma(x)=\sum p_\gamma \gamma$, as before. Composing the map $\chi$ with
\eqref{imagefunctor} yields the mapping $\ddd\to \ddd$
given by \begin{eqnarray*} b \ \mapsto\  \frac{1}{n!} \sum_g g(b)\otimes g &
 \buildrel{\strut\cdot \sigma(x)}\over{\longmapsto} &
 \frac{1}{n!} \sum_{g,\gamma}  g(b)  \otimes gp_\gamma\gamma \\&
= & \frac{1}{n!} \sum_{g,\gamma} g\gamma( \gamma^{-1}(bp_{\gamma})) \otimes
g\gamma\ \mapsto \ \sum_{\gamma} \gamma^{-1}(bp_{\gamma}).\end{eqnarray*} Thus
$\Psi^0 (\ogr R_x)$ coincides with $\widehat{\Phi}(R_x)$ as graded $A$-module
maps
from $S$ to $S$. As in the proof of Theorem~\ref{cohpro}, $\pi(S)=\PP\in \coh
\hi$ and so the sublemma follows.
\hfill\qed

We return to the proof of the proposition and let
  $(N,\Gamma)\in H_c\fmd$. By \cite[Section~2.3.11]{bj}, we can
construct an exact sequence in $H_c\fmd$
 \begin{equation} \label{freeres} \begin{CD} \cdots
@> d>>  F_i @> d>> F_{i-1} @> d>> \cdots
@> d>>  F_0 @> \eta>> N \longrightarrow 0\end{CD} \end{equation}
where the $F_i$ are filtered free $H_c$-modules   whose
 filtrations
shift the $\ord$ filtration and
the induced complex
\begin{equation} \label{2ndcomplex} \begin{CD} \cdots
@> \gr d>> \gr F_i @> \gr d>> \gr F_{i-1} @> \gr d>> \cdots
@> \gr d>>  \gr F_0 @> \gr \eta>> \gr_{\Gamma} N \longrightarrow 0\end{CD}
\end{equation}
is a free resolution of $\gr_{\Gamma}N$. In particular the maps $d$
are given by right multiplication by  matrices with entries from $H_c$
and so the sublemma can be applied to them.
By Theorem~\ref{cohpro}, if we  apply $\widehat{\Phi}$ to \eqref{freeres}
we obtain the
complex of coherent sheaves on $\hi$
\begin{equation}
\label{vbres} \begin{CD} \cdots @> \gr \tilde{d} >>
\mathcal{P}^{\oplus n_i} @> \gr \tilde{d} >>
\mathcal{P}^{\oplus n_{i-1}} @> \gr \tilde{d} >> \cdots @>
\gr \tilde{d} >> \mathcal{P}^{\oplus n_0}
@> \gr \tilde{\eta} >> \widehat{\Phi}(N) @>>>  0,\end{CD} \end{equation}
 where $n_i$ is the rank of $F_i$ over $H_c$ and $\gr \tilde{d}$ is
 obtained from $d$ by repeated use of \eqref{noncommfirst}.
 As in the proof of Lemma~\ref{graded-lemma}(1),
 the surjectivity of $\eta$ implies that
 final map $\gr \tilde{\eta}$
 in  \eqref{vbres}  is also surjective.
Unfortunately, as happened in the proof of Lemma~\ref{ses-cc}, one cannot
 expect that $\gr \tilde{d } $ is right exact and so
 we can only conclude that  the zeroth cohomology of the deleted complex
 \eqref{vbres} surjects onto $\widehat{\Phi}(N)$.

  On the other hand,
 $(\rho_{\ast}f^{\ast}(-))^\WW$ is a right exact functor
 and so  $\Psi^0(N)$
 is the zeroth cohomology of the functor $(\rho_{\ast}f^{\ast}(-))^\WW$ applied
 to the deleted complex  \eqref{2ndcomplex}.
 By the sublemma, this new
 complex is just the deleted complex \eqref{vbres}.
The proposition follows. \end{proof}

\subsection{} \label{cohvan} It is natural to ask whether
Proposition~\ref{compBKR} is the shadow of a much stronger statement.
\begin{question1} Let $(N,\Gamma)\in H_c\fmd$.
Is there a quasi-isomorphism of coherent sheaves between
$\widehat{\Phi}_{\Gamma} (N)$ and $ \Psi_{\Gamma}(N)$?
 In other words  is
$\Psi^0_{\Gamma}(N)\cong \widehat{\Phi}_{\Gamma} (N)$ and do the
homology sheaves $\mathrm{H}^i( \Psi_{\Gamma}(N))$ vanish for
$i>0$?\end{question1}

It is instructive to  rephrase the  final part of this question.
As we observed in the proof of Theorem~\ref{bkr}, the inverse
 to $\Psi$ is the functor $R\Gamma(\PP(-1) \otimes -) \otimes \sign$.
 By \cite[(45)]{hai3} $\PP (-1) \cong \PP^*$, the dual bundle of $\PP$. Thus the
following question is a special case of Question~1.

\begin{question2}
If $(N, \Gamma) \in H_c\md$, does
$\mathrm{H}^i \big(\hi,\, \PP^*\otimes \widehat{\Phi}_{\Gamma}(N)\big) =0$ for $i > 0$?
\end{question2}

This question is rather similar to  \cite[Conjecture~3.2]{hai1}
which asks the following: Let $\BB_1$ denote the tautological
rank $n$ bundle and $\PP_1$ the Procesi bundle
on $\hin$, as in the proof of Lemma~\ref{punctcoh}.
  Then:

 \begin{conjecture} \cite{hai1}  Does
 $\mathrm{H}^i(\hin,\, \PP_1^*\otimes \BB_1^{\otimes \ell})=0$ for
all $i>0$ and all $\ell \geq 0$.\end{conjecture}

\section{Finite dimensional $U_c$-modules and their characteristic
varieties}\label{var-sect}

\subsection{}  \label{howaboutfd}
We assume throughout the section that $c\in\C$ satisfies
Hypothesis~\ref{main-hyp}. Some of the most interesting $H_{c}$-modules are the
finite dimensional ones, which exist if (and only if) $c=r/n$, for some $r\in
\NN$ with $(r,n)=1$ \cite[Theorem~1.2]{BEGfd}. For such $c$, the irreducible
finite dimensional modules are just the  modules $L_c(\triv)$, as defined below
 in \eqref{subsec-3.7}. In particular the module $L_{1/n}(\triv)$ is
one-dimensional (see  \cite[Proposition~2.1]{BEGfd})
and each $L_{1/n+k}(\triv)$ is
induced from   $L_{1/n}(\triv)$ via the shift functors $\widetilde{S}$.

Our $\Z$-algebra techniques are particularly well suited to studying these
modules and we give the details in this section. We first prove in
Proposition~\ref{1d-prop} that $\widehat{\Phi}(L_{1/n}(\triv))$ is simply
$\OO_{\hio}$, where  $\hio$ is the punctual Hilbert scheme.  This statement
contains within it information about {\it each} $L_{1/n+k}(\triv)$ and  this
allows us to answer conjectures of Berest, Etingof and Ginzburg from
\cite[Section~7]{BEGfd}. Essentially, these questions ask for a bigraded
decomposition of the associated graded modules of $ L_{1/n+k}(\triv)$ and
$eL_{1/n+k}(\triv)$. However, the formal statement uses slightly different
filtrations from the ones use so far and so we refer the reader to
Theorem~\ref{beg-conj} for the precise statement.

\subsection{Category $\mathcal{O}_c$}\label{subsec-3.7}
 As in \cite{GGOR} or  \cite[Definition~2.4]{BEGqi},
we define $\mathcal{O}_c$\label{cat-O-defn}
 to be the abelian category of finitely-generated
$H_c$-modules $M$ for which  the action  of
$\mathbb{C}[\h^*]$ on $M$ is locally nilpotent. By
\cite[Theorem~3]{guay} $\OO_c$ is a highest weight category. Given $\mu \in
\irr{\WW}$, we define   the
\textit{standard module $\Delta_c(\mu)\in \mathcal{O}_c$},\label{standard-defn}
 to be the induced module
$$\Delta_c(\mu) = H_c\otimes_{\mathbb{C}[\h^*]\ast \WW} \mu,$$
where $\mathbb{C}[\h^*]\ast \WW$ acts on $\mu$ by   $pw\cdot m =
p(0) (w\cdot m)$, for $p\in \mathbb{C}[\h^*]$, $w\in \WW$ and
$m\in \mu$. It follows from the
PBW Theorem~\ref{PBW} that   $\Delta_c(\mu)$ is
a free left $\mathbb{C}[\h]$-module of rank $\dim(\mu)$.
It is shown in \cite[Section~2]{BEGqi} that   each
$\Delta_c(\mu)$ has a unique simple quotient
$L_c(\mu)$,\label{L-defn}
 that the set $\{ L_c(\mu) : \mu \in
\irr{W}\}$ provides a complete list of non-isomorphic simple
objects in $\mathcal{O}_c$, and that every object in
$\mathcal{O}_c$ has finite length.

\subsection{}\label{1d-trivia} The following observation will be used several
times.

\begin{lem} Let $U'$ and $U$ be subrings of a filtered $\C$-algebra
$D=\bigcup_{i\geq 0} F^iD$ and $T$ a $(U',U)$-submodule of $D$. Give $U'$, $U$
and $T$ the induced filtration and suppose that $U$ has a $1$-dimensional
module $\C=U/ \mathfrak{n}$ with the \emph{trivial  filtration}
 $F^r(\C)=F^0(\C)=\C$ for all $r\geq 0$.
Then, under the tensor product filtration on $T\otimes_U\C$, we have
$F^i(T\otimes \C)=(F^iT)\otimes \C$ and a  surjective $\gr_FU'$-module
map $\gr_FT/\gr_FT\gr_F \mathfrak{n} \twoheadrightarrow
\gr_F(T\otimes \C)$. \end{lem}

\begin{proof} Clearly $F^i(T\otimes \C)=\sum_r F^{i-r}T\otimes F^r(\C)
=F^i(T)\otimes \C = (F^iT+T \mathfrak{n})/T \mathfrak{n}$.
 This induces a  surjective map $F^iT/F^{i-1}T \twoheadrightarrow
  \left(F^iT+T \mathfrak{n}\right)\left/\right.
\left(F^{i-1}T+T \mathfrak{n}\right)
=  F^i(T\otimes \C)/F^{i-1}(T\otimes \C)$, and hence a surjection
$\gr_FT \twoheadrightarrow \gr_F(T\otimes \C)$.
The kernel of this map contains
$\gr_F(T\mathfrak{n}) \supseteq \gr_F(T)\gr_F(\mathfrak{n}).$
\end{proof}

\subsection{} \label{1d} Let $c=1/n$.
By \cite[Proposition~2.1]{BEGfd} or \cite[Lemma~4.6]{gordc} the
$H_c$-module $ L_c(\triv)$ is one dimensional.
Thus, if we give  $L_c(\triv)$ the trivial filtration
  $\Gamma^i (L_c(\triv)) = L_c(\triv)$ for
 $i\geq 0$, then $\gr_{\Gamma} L_c(\triv) = \C$ is the trivial
 $\cxy\ast \WW$-module concentrated in degree 0.
The first main result of this section identifies the $\hi$-module
$\widehat{\Phi} (L_c(\triv)).$

\begin{prop}\label{1d-prop}
 Let $c=1/n$ and give
 $L_c(\triv)\cong \C$ the trivial filtration $\Gamma$. Then
 $\widehat{\Phi}_{\Gamma}(L_c(\triv))\cong \OO_{\hio}$.
 \end{prop}

\begin{proof}
Recall that $A=\bigoplus_{k\geq 0} A^k\delta^k.$
Identify $eL_c(\triv)$ with $N=
 U_c/\widehat{\mathfrak{m}}$ for the appropriate ideal $\widehat{\mathfrak{m}}$.
Filter $N$ by  $F^0(N)=N$,
 and write $N(k)=B_{k0}\otimes_{U_c}N$ as in \eqref{graded-section}.
By \cite[Remark, p.306]{BEGqi} $\h$ and $\h^*$
annihilate $L_c(\triv)$, so
$\mathfrak{m}=\ogr(\widehat{\mathfrak{m}})
 = \cxy^{\WW}_+$, the maximal ideal of $A^0=\cxy^\WW$
corresponding to the zero orbit.
 For any $k\geq 0$, \eqref{OMT1} and
Lemma~\ref{1d-trivia} combine to give  a  surjection
$$\eta_k: \frac{A^k\delta^k }{A^k\delta^k{\mathfrak{m}}}
\ \cong \ \frac{\ogr B_{k0}}{\ogr B_{k0}\ogr\widehat{\mathfrak{m}}}\
 \twoheadrightarrow \ \gr N(k)$$ and, hence,  a surjection
$\eta: A/A\mathfrak{m}\twoheadrightarrow \gr(\widetilde{N})$.
On the other hand,  Lemma~\ref{punctcoh} shows that $
 \mathrm{H}^0 (\hio, \LL^k) \cong  A^k\delta^k/A^k\delta^k\mathfrak{m}$
and so  $\mathcal{O}_{\hio}$ is the   $\hi$-module
 corresponding to the graded $A$-module $A/A\mathfrak{m}$.

 In order  to complete the proof of the proposition it therefore suffices to
 show  that $\eta$ is an isomorphism, which we do by computing dimensions.
  By comparing  Lemma~\ref{fddim} and \cite[Theorem~1.14]{BEGfd} we obtain
$$ \dim A^k\delta^k/A^k\delta^k \mathfrak{m} = \frac{1}{kn+1}\binom{(k+1)n}{n} =
 \dim eL_{c+k}(\triv) = \dim N(k).$$
 Since the   maps  $\eta_k$ are   surjections they are therefore
isomorphisms. Hence, so is $\eta$.
 \end{proof}

\subsection{}\label{findim-question}
As remarked in \eqref{howaboutfd}, the simple module $L_c(\triv)$
is finite dimensional whenever $c= r/n$ with $r\in \NN$ and
$(r,n)=1$. It would be interesting to know whether an analogue  of
Proposition~\ref{1d-prop} holds for  these modules:

\begin{problem} If $c = r/n$ with $(r,n)=1$,
find $\widehat{\Phi}(L_c(\triv))$ under some natural
filtration on $L_c(\triv)$.
\end{problem}

A shadow of the answer will be  given in \eqref{simp-cc-fd}
which computes the restricted character of $L_c(\triv)$.

\subsection{} By \cite[Proposition~3.16]{GS},
$B_{k0}\otimes_{U_c}eL_c(\triv)\cong eL_{c+k}(\triv)$ for $k>0$ and so
Proposition~\ref{1d-prop} also contains information about the modules
$eL_{1/n+k}(\triv)$. For the rest of this section we will study the
consequences of this observation. In particular we will   determine
the bigraded structure of the associated graded modules of
these modules, thereby answering the
conjectures from \cite[Section~7]{BEGfd}.

\subsection{Torus action revisited} \label{toract2}
 As we will discuss, this bigraded structure  is the standard one
for rings  of differential operators, but  it is not quite
 the bigrading used in \eqref{toract1} and \cite{hai1}.
Instead, we need to consider a second \label{bigrad-defn}action of
$\TT^2$ on $\C[\C^2] = \C[x,y]$, this time given by $\tau_{s,t}
\cdot x = st x, \ \tau_{s,t} \cdot y = st^{-1}y$. This induces new
$\TT^2$-structures on $\hi$, $\hio$, $\PP$, $\LL$, etc.

It is routine to check that,  if an element $f$ is bihomogeneous
of weight $(i,j)$ for the action in \eqref{toract1}, then under
this new action $\tau_{s,t} \cdot f = s^{i+j} t^{i-j}f$
 and so $f$ is bihomogeneous of weight  $(i+j,i-j)$.
Thus,  under the present torus action,  Lemma~\ref{punctcohref}
provides the  following identification of bigraded components:
\begin{equation}\label{punctcohref99}
\mathrm{H}^0(\hio , \LL^k)^{ij} \cong \left(
\frac{A^k\delta^k}{A^k\delta^k\mathfrak{m}}\right)^{i+Nk,j+Nk},
\qquad \mathrm{H}^0(\hio , \PP\otimes \LL^k)^{ij} \cong \left(
\frac{J^k\delta^k}{J^k\delta^k\mathfrak{m}}\right)^{i+Nk,j+Nk}.
\end{equation}

\subsection{Euler operator}\label{gradingsec}  As in \cite[(2.4)]{GS}
the rings of differential operators $D(\h)$ and $D(\hr)$ have
a graded structure, given by the adjoint action $[\EE,-]$ of the
 {\it Euler operator } $\EE=\sum x_i\partial_{i}\in D(\h)$.
 \label{Euler-defn} We will call this the {\it Euler grading}
 and write $\Edeg$\label{Euler-deg-defn}
 for the corresponding degree function; thus
 $\Edeg x_i=1$ and $\Edeg \partial_{i}=-1$.
 Since $\EE\in D(\h)^\WW$, $\EE$ commutes
 with $\WW$ in $D(\hr)\ast \WW$
 and so this grading extends to that ring with $\Edeg \WW=0$.
 By inspection, \eqref{dunkop} implies that
 the elements $y_i$ defined there
  also have degree $-1$ and so each $H_c$ is also graded under
 $[\EE,-]$ and we continue to call this the Euler grading.

It is well-known and easy to check that the $\EE$-grading is
compatible with the order filtration on $D(\hr)\ast \WW$, in the
sense that $[\EE, \ord^n D(\hr)\ast \WW]\subseteq \ord^n D(\hr)\ast
\WW$ for all $n\geq 0$. Thus,  we obtain an induced grading, again called the
$\EE$-grading, on the associated graded object $\ogr D(\hr)\ast
\WW\cong \C[\hr\oplus \h^*]\ast W$ and its subrings. Clearly this
is again given by $\Edeg \h^*=1$ (which we define to mean that
$\Edeg(x)=1$ for every   $0\not=x\in \h^*$) while
 $\Edeg \h=-1$ and $\Edeg \WW=0$.
In other words, the $\EE$-grading coincides with the grading
induced by the second component of $\TT^2$ in \eqref{bigrad-defn}; which
 is the point of  defining that action.

\subsection{} \label{totfil} The first conjecture from \cite[Section~7]{BEGfd}
asks whether the isomorphism from  Proposition~\ref{1d-prop}
 can be extended to an appropriate
$\TT^2$-equivariant isomorphism. However the $\TT^2$-action in \cite{BEGfd} is
not the one arising from the $\ord$ filtration and $\EE$-grading;
we must  replace
the  $\ord$ filtration by the filtration $\tot$  defined by the total
 degree of differential operators. Fortunately, as we next show, it is easy to pass
 between the two filtrations.

We begin with the formal definitions.
The filtration $\tot$\label{total-filt-defn}
  is defined on $D(\hr)\ast \WW$  and its subspaces by giving
$\WW$ total degree zero and $\h\oplus \h^*\smallsetminus \{0\}$
   total degree one.  As usual, for any subset $M$ of
$D(\hr)\ast \WW$, we write
$\tot^r M$ for the elements $m\in M$ of total degree $\tdeg(m)\leq r$
and   denote the associate graded object  by
$\tgr M = \bigoplus_{r\geq 0} \tgr^r M$,
where $\tgr^rM =\tot^rM/\tot^{r-1}M$.
It follows from \eqref{gradingsec} that commutation by $\EE$ preserves the
$\ord$ and $\tot$ filtrations on $D(\hr)\ast \WW$ and therefore induces a
grading on $B_{ij}$ and its associated graded modules.
 Write $\tot^{r,s} B_{ij}$ (respectively
$\ord^{r,s} B_{ij}$) to denote the elements of $\tot^rB_{ij}$
 (respectively $\ord^rB_{ij}$) with  $\EE$-degree equal to $s$.
The $\EE$-graded decompositions
$\tgr^r B_{ij} =\bigoplus_{s\in\Z} \tgr^{r,s}B_{ij}$
and  $\ogr^rB_{ij}=\bigoplus \ogr^{r,s}B_{ij}$ are defined analogously.
By construction, this bigraded structure on $\tgr B_{ij}$
is simply the one coming from the
$\TT^2$-action of \eqref{bigrad-defn}.

\begin{lem}
For each $i\geq j \geq 0$ there is an isomorphism $\tgr
B_{ij}\to \ogr B_{ij}$ which identifies $\tgr^{2r+s,s} B_{ij}$
with $\ogr^{r,s} B_{ij}$.  These isomorphisms induce a $\Z$-algebra isomorphism
$\tgr B \to \ogr B.$   \end{lem}

\begin{proof}  We fix $i\geq j\geq 0$ throughout the proof.
 By \cite[Lemma~6.10 and Corollary~6.14(1)]{GS} the
action of $\EE$ on $B_{ij}$ is diagonalisable. Therefore
$$\tgr B_{ij} = \bigoplus_{s\in \Z} \bigoplus_{r\in \Z}\,
\frac{\tot^{r,s}B_{ij}}{\tot^{r-1,s}B_{ij}}\qquad\mathrm{and}\qquad
\ogr B_{ij} =
\bigoplus_{s\in \Z} \bigoplus_{r\in \Z}\,
\frac{\ord^{r,s}B_{ij}}{\ord^{r-1,s}B_{ij}}.$$
If $b\in \ord^r B_{ij}\smallsetminus \ord^{r-1}B_{ij}$, set
$\odeg (b)=r$.
Then $ \tdeg (b) =
\Edeg (b) + 2 \odeg (b)$ and so
\begin{eqnarray*}
\ord^{r,s} B_{ij}&=& \{ b\in B_{ij}: \odeg (b) \leq r
\text{ and }\Edeg(b) = s \} \\ &
= &\{ b\in B_{ij} : \tdeg (b)
\leq 2r+s \text{ and } \Edeg(b)=s \} \\ &
 = & \tot^{2r+s,s}B_{ij}.
\end{eqnarray*} This proves the first claim of the lemma.
The second claim follows since the map $\tgr B \to \ogr B$ is
induced by the identity map on $B$.
 \end{proof}

\subsection{}\label{subsec-7.7}  Assume that $c=1/n$.
By \cite[Proposition~3.16]{GS} we can  identify
 \begin{equation}
\label{tensorfdsimple1}
eL_{c+k}(\triv) \cong eH_{c+k}\delta e_-\otimes_{U_{c+k-1}}eH_{c+k-1}\delta e
\otimes_{U_{c+k-2}}\cdots \otimes_{U_{c+1}} eH_{c+1}\delta e \otimes_{U_c}\C
\end{equation}
and
 \begin{equation}\label{eq04}
 L_{c+k}(\triv)\cong
  H_{c+k}\delta e_-\otimes_{U_{c+k-1}}eH_{c+k-1}\delta e
   \otimes_{U_{c+k-2}}\cdots \otimes_{U_{c+1}} eH_{c+1}\delta e
    \otimes_{U_c}\C.
  \end{equation}
These modules have an $\EE$-gradation and a filtration induced from the tensor
product $\tot$-filtration, which  induce a  bigraded structure on
$\tgr L_{c+k}(\triv)$ and $\tgr eL_{c+k}(\triv)$. Moreover, since the
$\tot$ filtration is $\WW$-stable, $\tgr
L_{c+k}(\triv)$ is then a bigraded $\WW$-module.

We need  a minor extra adjustment to
 the $\tot$ filtration in order to state (and prove) the conjectures
 from \cite{BEGfd}. To be precise, the conjectures are
 stated in \cite{BEGfd}
for the tensor product filtration arising from the isomorphism
\begin{equation}\label{eq05}
eL_{c+k}(\triv)  \cong  eH_{c+k}e_-
\otimes_{\UU_{c+k-1}} eH_{c+k-1}e_- \otimes_{\UU_{c+k-2}} \cdots
\otimes_{\UU_{c+1}} eH_{c+1}e_- \otimes_{\UU_c} \C.
\end{equation}
 In other words, one ignores  the copy of $\delta$ in each tensor summand
 of $eL_{c+k}$.
  We will write the $\tot$ filtration arising from \eqref{eq05} as
$\tot_{\mathrm{B}}$. Of course, this is only a book-keeping device
to account for the powers of $\delta$
and so it follows immediately that
 \begin{equation} \label{eq1}\tgr^{r,s}_{\mathrm{B}} eL_{c+k}(\triv) =
 \tgr^{r+Nk,s+Nk}  eL_{c+k}(\triv)
 \qquad\mathrm{where}\quad N=n(n-1)/2=\deg \delta.\end{equation}

Using \eqref{eq04} in place of \eqref{tensorfdsimple1},
 the same comments apply to $L_{c+k}(\triv)$ and give
 \begin{equation} \label{eq15}  \tgr^{r,s}_{\mathrm{B}} L_{c+k}(\triv) =
 \tgr^{r+Nk,s+Nk} L_{c+k}(\triv).\end{equation}

\subsection{}\label{beg-conj-sect}
We now confirm
\cite[Conjectures~7.2 and~7.3]{BEGfd}.

\begin{thm}\label{beg-conj}
Let $c=1/n$ and   $k\in \NN$.  Give $L_{c}(\triv)$ the trivial
filtration, and consider the $\TT^2$-action of \eqref{toract2}.
\begin{enumerate}
\item{}
Give $eL_{c+k}(\triv)$ the tensor product filtration arising from
\eqref{tensorfdsimple1}. Then  there is
 a bigraded isomorphism $\tgr_{\mathrm{B}} eL_{c+k}(\triv) \cong
 \mathrm{H}^0(\hio,\mathcal{L}^{k}).$
\item{}
Give $L_{c+k}(\triv)$ the tensor product filtration arising from the identity
 \eqref{eq04}.  Then there is a $\WW$-equivariant,
bigraded isomorphism
$$\tgr_{\mathrm{B}} L_{c+k}(\triv) \cong \mathrm{H}^0(\hio , \mathcal{P}\otimes
\mathcal{L}^{k-1})\otimes \sign.$$
\end{enumerate}
 \end{thm}

\begin{proof}
(1) By Lemma~\ref{totfil} and \cite[Lemma~7.2]{GS}, the tensor
product filtration on $B_{ij}$ induced from the $\tot$ filtration
on the $B_{k+1,k}$ coincides with the $\tot$ filtration on
$B_{ij}$. So we need not distinguish between them. As in the proof
of Proposition~\ref{1d-prop},   identify $L=eL_c(\triv)=\C =
U_c/\widehat{\mathfrak{m}}$ and give $L$ the trivial filtration.
Then, for any $k\geq 0$, the induced  tensor product filtration on
$L(k)\cong eL_{c+k}(\triv)$ is equal to the one coming from the
tensor product filtration on the right hand side of
\eqref{tensorfdsimple1}. By  Lemma~\ref{1d-trivia}
 there is a surjection
  $\eta^k:\tgr (B_{k0})/\tgr (B_{k0})\tgr (\widehat{\mathfrak{m}})
  \twoheadrightarrow \gr L(k).$  This  is clearly
  a graded morphism under both the
 $\tot$ and $\EE$-gradings.

 Now $\tgr \widehat{\mathfrak{m}} = \mathfrak{m} = \cxy_+^{\WW}$. Thus
 \eqref{OMT1} and Lemma~\ref{totfil} give bigraded isomorphisms
\begin{equation} \label{eq2}
\left(\frac{\tgr B_{k0}}{(\tgr B_{k0})(\tgr
\widehat{\mathfrak{m}})}\right)^{r+Nk,s+Nk} \cong
\left(\frac{A^{k}\delta^{k}}{A^{k}\delta^{k}\mathfrak{m}}\right)^{r+Nk,s+Nk}
\end{equation}
 for all $r,s\in \Z$.
 By \eqref{punctcohref99}, this final term is isomorphic to
 $\mathrm{H}^0(\hio, \LL^{k})^{r,s}$. Thus,  \eqref{eq1} and \eqref{eq2}
 combine to show that,
   for each $r$ and $s$, the map  $\eta^k$ restricts to a
     surjection
 $\eta_{r,s}: \mathrm{H}^0(\hio, \LL^{k})^{r,s}\twoheadrightarrow\tgr^{r,s}_{\mathrm{B}}
  eL_{c+k}(\triv).$
  Finally, by   Lemma~\ref{fddim} and \cite[Theorem~1.14]{BEGfd},
   we know that
  $\dim_\C \mathrm{H}^0(\hio, \LL^{k})=\dim eL_{c+k}(\triv)<\infty$.
Therefore, each $\eta_{r,s}$ must be an isomorphism  and we are done.

(2) In this case,  the filtration on each $L_{c+k}(\triv)$
is the one coming from the tensor product filtration on the right hand side
of \eqref{eq04}. However, by \cite[Lemma~7.2]{GS} this is the
same filtration as that coming from the tensor product filtration on
 $L_{c+k} \cong M\otimes_{U_c}\C$, where $M=H_{c+k}\delta e B_{k-10}$.
 Thus,  by \eqref{eq15},
\begin{equation}\label{eq24}
\tot^{r,s}_B L_{c+k}=
 \tot^{r+Nk,s+Nk}M\otimes \C.\end{equation}
By  Lemma~\ref{1d-trivia} there is  a  surjection
$\chi: \tgr M/(\tgr M)(\tgr \widetilde{\mathfrak{m}}) \twoheadrightarrow
\tgr(L_{c+k}),$ which is  bigraded under the $\tot$ and $\EE$ gradings.
By \cite[Proposition~B.1]{GS},
\begin{equation}\label{eq25}
\left(\frac{\tgr M}{\tgr M\tgr
\widehat{\mathfrak{m}}}\right)^{r+Nk,s+Nk} \cong\
\left(\frac{J^{k-1}\delta^{k}}{J^{k-1}\delta^{k}\mathfrak{m}}\right)^{r+Nk,s+Nk}
\cong \  \sign \otimes
\left(\frac{J^{k-1}\delta^{k-1}}{J^{k-1}\delta^{k-1}
\mathfrak{m}}\right)^{r+N(k-1),s+N(k-1)}
\end{equation}
 for all $r,s\in \Z$.
By \eqref{punctcohref99}, again,  the right hand side of \eqref{eq25}
 is isomorphic to $\rm{H}^0(\hio, \PP\otimes \LL^{k-1})^{r,s}\otimes \sign$. By
construction and Lemma~\ref{punctcoh},
 each of these identifications is $\WW$-equivariant.
  Therefore,  \eqref{eq25} and \eqref{eq24}
combine to give a $\WW$-equivariant surjection
$$\chi^{r,s}:  \mathrm{H}^0(\hio, \PP \otimes \LL^{k-1})^{r,s}
 \otimes \sign \twoheadrightarrow
\tgr^{r,s}_B L_{c+k} $$ for each $r$ and $s$.
By \cite[Theorem~1.11]{BEGfd} and Lemma~\ref{fddim} we know that
 $
\dim  \tgr_B L_{c+k} = (kn+1)^{n-1}=
 \dim\mathrm{H}^0(\hio, \PP \otimes \LL^{k-1})$
for all $k\geq 1$. Hence $\chi^{r,s}$ is an isomorphism, as required.
  \end{proof}

\subsection{}
\label{bigrfd} Using \eqref{bigradedactchar} and
\eqref{bigradedactchar2} one can explicitly write down the two
variable Poincar\'e series  for $L_{c+k}(\triv)$ and
$eL_{c+k}(\triv)$ arising from this theorem.  Of course, by
\eqref{toract2},  those  formul\ae\ need to be slightly modified
to apply here but, for example, it follows   from
\eqref{bigradedactchar2} that under the grading of {\it this}
theorem, $eL_{c+k}(\triv)$ has the two-variable Poincar\'e series
$p(eL_{c+k}(\triv),s,t) = C^{(k)}_n(st, st^{-1}).$
%%%%%%%%%%%%%%%%%%%%%%%%%%%%%%%%%%

\section{Characteristic cycles}\label{ses-sect}

\subsection{} In this section we expand on the results from
Sections~\ref{cohsh} and \ref{var-sect} in order to understand the
characteristic varieties $\Chr M \subseteq \hi $ and the
characteristic cycles $\Ch M $ for  category $\mathcal{O}$-modules
$M$ (the reader should recall the notation and conventions  from
Section~\ref{zalg} concerning these objects).
   Even for a standard module $\Delta=\Delta_c(\mu)$, the
answer is quite subtle and in marked contrast to the usual
associated variety $\Chro \Delta\subseteq \h\oplus\h^*/\WW$.
Indeed,  $\Chro \Delta$ is clearly independent of $\mu$  whereas
$\Chr \Delta$  varies with  and even determines $\mu$. The precise
result is given in Theorem~\ref{stand-cc} which also gives an
explicit formula for the characteristic cycle $\Ch \Delta$. Our
results are less complete for other objects of $\OO_c$, but they
are sufficient to show in Corollary~\ref{borel} that $K(\OO_c)$ is
isomorphic to the top degree Borel-Moore homology group
$\mathrm{H}_{2n-2}(Z,\C)$ for the variety $Z=\tau^{-1}(\h/\WW)$.

As usual we will assume that $c$ satisfies
Hypothesis~\ref{main-hyp} throughout the section.

\subsection{}\label{Z-grading-var}
Before proceeding, we will need a minor variant of
\cite[Corollary~4.13]{GS}. Let
$\C[\h]^{\text{co} \WW} = {\C[\h]}\big/{\C[\h]_+^{\WW}\C[\h]} $
denote the algebra of coinvariants.
 For $\mu\in\irr{\WW}$   define   the {\it fake degree of $\mu$}
 to be the polynomial
 \label{fakedegrees} $ f_{\mu}(v) = \sum_{i\geq 0}
[\C[\h]^{\text{co} \WW} _i : \mu] v^i.$

 Recall from \eqref{gradingsec} that $\cxy$ has a natural $\Z$-grading
 called the $\EE$ grading. Given an $\EE$-graded module $M=\bigoplus M_i$,
  we will write the corresponding Poincar\'e series as
  $p(M,v)=\sum v^i \dim_\C M_i$.\label{mono-poincare}

\begin{lem}\label{subsec-5.10}
 Set   $N=n(n-1)/2$ and write $K=\h\cxy \subset \cxy$ for the ideal
generated by the elements of $\h$. Then
$J^{k}\delta^{k}/J^{k}\delta^{k}K$ is $\EE$-graded  for any $k\geq
1$ and has   Poincar\'e series
\begin{equation} \label{pseriesP1}
p(J^{k}\delta^{k}/J^{k}\delta^{k}K, v)
 = v^{Nk}\frac{\sum_{\mu}
f_{\mu}(1)f_{\mu}(v) v^{k(n(\mu) - n(\mu^t))}}{\prod_{i=2}^n
(1-v^{i})}.
\end{equation}
\end{lem}

\begin{proof} As in \cite[(4.13)]{GS}, it is routine to check that  the given
modules are $\EE$-graded.

Let $C=\C[\h^{*}]^{\WW}_{+}$. Since $J^1$ is symmetric in $\h$ and
$\h^{*}$, it follows from \cite[Corollary~4.13]{GS} that
$J^{k}/J^{k}C$ has Poincar\'e series
\begin{equation} \label{pseriesP2}
p(J^{k}/J^{k}C, v)
 = \frac{\sum_{\mu}
f_{\mu}(1)f_{\mu}(v) v^{k(n(\mu) -
n(\mu^t))}[n]_{v^{-1}}!}{\prod_{i=2}^n (1-v^{i})},
\qquad\mathrm{where}\   [n]_v!=\frac{\prod_{i=1}^n(1-v^i)}{(1-v)^n}. \end{equation}

 The derivation of this formula in \cite{GS} was obtained by noting that  the
fundamental invariants, which provide the generators of $C$, form
an r-sequence in  each $J^k$. This allows one to obtain  the
Poincar\'e series for $J^{k}/J^{k}C$   from the Poincar\'e series
for $J^k$ given by \cite[Corollary~4.11]{GS}. By
\cite[Lemma~4.4(2)]{GS}, the natural generators of $K$ also form
an r-sequence in $J^k$ and so the same argument can be  used to
find the Poincar\'e series of $J^k/J^kK$.

In more detail,
 the generators of $K$ all have degree $-1$
while the fundamental invariants  have degrees $-2\geq -r\geq -n$ and
 $\delta^k$ has degree $Nk$.  In the notation of  \cite[(4.11)]{GS},
 the formula \cite[(4.13.3)]{GS} gives
$$   p(J^{k}/J^{k}C, v) =
\left(\prod_{i=2}^n(1-t^i)p(J^k,s,t)\right)_{s=v,t=v^{-1}}.$$
Using the same argument one obtains:
 $$p(J^{k}\delta^{k}/J^{k}\delta^{k}K, v)=
  v^{Nk} \Big((1-t)^{n-1}p(J^k,s,t)\Big)_{s=v,t=v^{-1}}.$$
  Comparing these two formul\ae\    gives
   $$p(J^{k}\delta^{k}/J^{k}\delta^{k}K, v)=
 \frac{v^{Nk}(1-v^{-1})^{(n-1)}}{\prod_{i=2}^{n}  (1-v^{-i})}\,
 p(J^{k}/J^{k}C, v)= \frac{v^{Nk}}{[n]_{v^{-1}}!}\,
 p(J^{k}/J^{k}C, v).$$ Combined with  \eqref{pseriesP2}  this implies
 \eqref{pseriesP1}. \end{proof}

\subsection{}  \label{subsec-7.12}
The next few results   give partial information on the
characteristic cycle $\Ch \Delta$ for a standard module $\Delta$ that will be
used in the complete description of  $\Ch \Delta$ in Theorem~\ref{stand-cc}.

 For a $\WW$-module $V$ and $\mu\in \irr \WW$   set
 $V_{\mu} = (V\otimes \mu)^{\WW}$. By construction,
  the fibres of the Procesi bundle
 $\PP$ carry the regular $\WW$-representation and so  $\PP_{\mu}$ is a vector
 bundle of rank $\dim \mu$.
  Write $K =\h\cxy$ as in \eqref{Z-grading-var}.
  It follows from \eqref{A-1-defn} that $\cxy \cong \OO(X_n)
 = \mathrm{H}^0(\hi, \PP)$. Thus $\cxy$
 acts on $\PP$ and so $\PP K$ and $(\PP K)_\mu$ are naturally defined for any
 $\mu\in\irr\WW$.

Let $c\geq 0$. We give $\Delta_c(\mu)$ the induced filtration
$\Gamma^t (\Delta_c(\mu)) = \ord^t H_c\cdot (1 \otimes \mu)$.
Since $\Delta_c(\mu)\cong \C[\h]\otimes \mu$, this is the same as
the trivial filtration $\Gamma^t (\Delta_c(\mu)) = \Delta_c(\mu)$
for all $t$, but the present formulation will be more convenient
in the sequel.

\begin{lem} For each $\mu\in\irr\WW$  there is an isomorphism
$\widehat{\Phi}_{\Gamma}(\Delta_c(\mu)) \cong  \PP_{\mu}/( \PP K)_{\mu}$
of sheaves on $\hi$.
\end{lem}

\begin{proof}
Let $X=  H_c\otimes_{\cy} \C$, filtered by
$\Gamma^t (X) = \ord^t H_c\cdot (1 \otimes \C)$, and   set
$$X(k) \ = \ S_{c+k-1}\circ \cdots \circ S_c (eX)\
 \cong\  eH_c(k) \otimes_{\cy} \C, \ = \
 B_{k0}\otimes_{U_c} eH_c\otimes_{\cy}\C\qquad \mathrm{ for}\ k\in \NN.$$
The filtration $\Gamma$ on $X$ induces the tensor product
 filtration on $X(k)$, which we also call $\Gamma$. By
\cite[Lemma~7.2]{GS}, the tensor product and $\ord$ filtrations on
$ B_{k0}\otimes_{U_c} eH_c$ coincide and so
$\Gamma^t X(k) = (\ord^t eH_{c}(k))\cdot ( 1\otimes \C)$ for all $t$ and $k$.
  Thus, by
Theorem~\ref{cohpro}, there  is a surjective mapping
$p':  J^k\delta^k \twoheadrightarrow \gr_{\Gamma} X(k).$ Moreover,
since $\cy$ acts trivially on $1\otimes 1\in X$, this factors through
a surjection
$$p: \frac{J^k\delta^k}{J^k\delta^k K} \longrightarrow \gr_{\Gamma} X(k).$$

Since the $\EE$-gradation commutes with the $\ord$-filtration and the
isomorphism of Theorem~\ref{cohpro} is $\EE$-graded,  the
 map  $p$ is $\EE$-graded.
Comparing \cite[Corollary~6.14(2)]{GS} with Lemma~\ref{Z-grading-var}
shows that $X(k)$ (which equals $\underline{N(k)}$ in the notation of
\cite{GS}) and $J^k\delta^k/J^k\delta^kK$
have the same Poincar\'e series under the $\EE$-grading.
Since $p$ is surjective
on each (finite dimensional) component, it is therefore an isomorphism.
Equivalently, $p$ induces an isomorphism $\gr_\Gamma \widetilde{X} \cong
\bigoplus_{k\geq 0} J^k\delta^k/J^k\delta^kK.$

There is a natural action of $\WW$  on $X(k)$ given by $w\cdot
(r\otimes 1) = rw^{-1}\otimes 1$ for $w\in \WW$ and $r\in H_c(k)$,
with the induced action on $\gr_{\Lambda} X(k)$. We next show that
$p$ is $\WW$-equivariant under this action.
 For $r\in H_{c}(k)$, let $\bar{r}$ denote the
image of the principal symbol $\sigma(r)$ in $J^k\delta^k/J^k\delta^kK$.
The definition of multiplication in $D(\hr)\ast \WW$ from \eqref{chered-defn}
 says that  $\omega r \omega^{-1}=\omega(r)$, for $r\in
\C[\h]\C[\h^*]\subseteq H_{c+k}$
 and $\omega\in \WW$.  By \cite[Lemma~6.11(1)]{GS} and \eqref{PBW},
 $eH_c(k)\subseteq eH_{c+k}=e\C[\h]\C[\h^*]$
and so the construction of $p'$  implies that, under the natural action of
$\WW$ on $J^k\delta^k$,
  $$p(\omega(\bar{r}))\ \equiv \ e \omega(\sigma(r)) \otimes 1
\ \equiv\ e\omega \sigma(r) \omega^{-1}\otimes 1 \ \equiv\
 e \sigma(r) \omega^{-1} \otimes 1
\ \equiv \ \omega(p(\bar{r}))\qquad \mathrm{mod} \, J^k\delta^kK$$
 for all  such $r $ and $\omega$. Thus
  $p$ and (by construction)  the  functors $S_{c+j}$ are indeed
   $\WW$-equivariant.

Next,  give
$H_c\otimes_{\cy\ast \WW} \C\WW $ the  $\WW$-action analogous to that on
$X$ so  the   isomorphism
$X\cong H_c\otimes_{\cy\ast \WW} \C\WW $ is $\WW$-equivariant. Then
\begin{equation} \label{factoronRees22} X_\mu \ \cong \
[H_c\otimes_{\cy\ast \WW}\C\WW]_{\mu} \ = \
H_c\otimes_{\cy\ast \WW}(\C\WW)_{\mu} \ = \  H_c\otimes_{\cy\ast \WW} \mu
\ = \  \Delta_c(\mu). \end{equation}
Since $p$ and the $S_j$ are $\WW$-equivariant, \eqref{factoronRees22}
induces an isomorphism
\begin{equation}\label{factoronRees23}
 (J^k\delta^k/J^k\delta^kK)_\mu \cong (\gr_\Lambda X(k))_\mu =
 \gr_\Lambda( X_\mu(k)) =
 \gr_\Lambda e\Delta_c(\mu)(k). \end{equation}
Thus  \eqref{factoronRees23} induces
 an isomorphism of graded $A$-modules
 \begin{equation} \label{standcoh}
\widetilde{e\Delta_c(\mu)}\ \cong\ \bigoplus_{k\geq
0}\left(\frac{J^k\delta^k }{J^k\delta^kK}\right)_{\mu}.
\end{equation}

Finally, we need to translate this equation to one involving sheaves.
Consider the  short exact sequence
\begin{equation} \label{factoronRees} 0\longrightarrow
\bigoplus_{k\geq 0} J^k\delta^kK\longrightarrow \bigoplus_{k\geq
0}J^k\delta^k \longrightarrow \bigoplus_{k\geq 0}
\frac{J^k\delta^k}{J^k\delta^kK}\longrightarrow 0
\end{equation} of $A$-modules.
By \eqref{A-1-defn}, the image of this sequence
in $A\lqgr\cong \coh \hi$ is isomorphic to  the short exact sequence of sheaves
$0\to \PP K \to \PP \to
\PP/\PP K\to 0$. By \eqref{good-P} this isomorphism is $\WW$-equivariant and so
   the   module
$\bigoplus_{k\geq 0} (J^k\delta^k)_{\mu}\in A\lgr$ has image
$\PP_\mu\in \coh \hi$. Combining these  facts shows that
 the image of
\eqref{standcoh}  in $\coh \hi$ is precisely   $\PP_\mu/(\PP K)_\mu$.
 \end{proof}

\subsection{}\label{Z-variety}
To describe the characteristic cycles of objects from $\OO_c$ we need to
recall a set of  Lagrangian subvarieties which were introduced by
Grojnowski in \cite[Section~3]{gro} (see also \cite[Sections~7.2 and 9.2]{Nak}).

Recall from  \eqref{hin-defn} that we  identify
$\C^{2n}/\WW$ with $S^n\C^2$, write $\widehat{\tau}: \hin \to \C^{2n}/\WW$
and identify $\h$ with a subspace of $\C^n\times\mathbf{0}\subset \C^{2n}$.
Set $\widetilde{C} = \widehat{\tau}^{\,-1}(\C^n\times\mathbf{0})$;
\label{tilde-C-defn} equivalently $\widetilde{C} =  \{
I\in \hin: \supp \C[x,y]/I \subseteq S^nC\}$ where
 $C =\{ (z,0): z\in \C\} \subseteq \C^2$.
By \cite[Proposition~3]{gro},  $\widetilde{C}$ is a pure variety
of dimension $n$. To describe its irreducible components, we use
the stratification $S^n\C^2 = \coprod_{\lambda\vdash n}
S_{\lambda}\C^2$ where, for a partition $\lambda =
(\lambda_1,\ldots ,\lambda_r)$ with exactly $r$ non-zero
parts,
$$S_{\lambda}\C^2 = \{ \sum_{i=1}^r \lambda_i [x_i] : \ x_1,
\dots ,x_r\in \C^2 \text{ are\ distinct} \}.$$
By \cite[Proposition~3]{gro} or \cite[p.111]{Nak}
the irreducible components
of $\widetilde{C}$ are the closures $\overline{C}_\lambda$  of the
subvarieties
$$\widetilde{C}_{\lambda} \ =\
 \widetilde{C}\, \cap \, \{ I \in \hin: \supp \C[x,y]/I\in S_{\lambda}\C^2\}.$$
Moreover, $\dim \overline{C}_\lambda = n$ for all $n$.

Set $Z = \tau^{-1}(\h/\WW)$ which, by  the last paragraph, is
equal to $\widetilde{C}\cap \hi$. We will occasionally write
$Z=Z(n)$ when we need to specify $n$. In order to identify the
components of $Z$, we recall from \cite[Proof of
Corollary~4.10]{GS} that the identification $\C[\C^{2n} ]^\WW =
\cxy^\WW[\mathbf{z}, \mathbf{z}^*]$ induces
 a factorisation  $\hin=\hi\times \C^2$, where
$\C^2=\mathrm{Spec}\, \C[\mathbf{z}, \mathbf{z}^*]$. It follows that
$\widetilde{C} = Z\times \C$ and so we have proved the following result.

\begin{lem}\label{Z-variety1}
The variety $Z$  is  pure of  dimension $n-1$ with
 irreducible components   $$
Z_{\lambda} = \overline{\strut \hi \cap \widetilde{C}_{\lambda}}
\qquad\mathrm{for}\  \lambda \in \irr \WW.    \eqno\qed $$
\end{lem}

In this notation, the punctual Hilbert scheme $\hio$ equals $Z_{(n)}$.

\subsection{Lemma}
\label{st-equi}
{\it If $\mu\in\irr\WW$ then $\Chr ( \Delta_c(\mu))$
is a union of subvarieties $Z_{\lambda}$, for some  $\lambda\in \irr \WW$.
 In particular, $\Chr ( \Delta_c(\mu))$ is equidimensional.
}

\begin{proof}
 By Lemma~\ref{subsec-7.12}, we need to calculate the
support variety of the sheaf $\PP_\mu/(\PP K)_\mu$.  By \eqref{factoronRees},
$$ \tau \left(\supp \PP/\PP K \right) \ = \
\mathcal{V}\left(K\cap \cxy^\WW\right) \ = \ \h/\WW, $$ as a subvariety of
$ \h\oplus \h^*/\WW.$
Therefore,  by the definition of $Z$, we have
$\supp \PP_{\mu}/(\PP K)_{\mu}\subseteq Z$.

Let $\{y_1,\ldots ,y_{n-1}\}$ be a basis of $\h$. Arguing as in
\cite[(4.9.2)]{GS} we see that $y_1,\ldots ,y_{n-1}$ is a regular
sequence in $S=\cxy[tJ^1\delta]$, the Rees ring of $X_n$, and so
we have an $\WW$-equivariant Koszul resolution
\begin{equation*} \label{kosres2} 0\rightarrow \OO_{X_n} \otimes
\bigwedge^{n-1}\h^* \rightarrow \cdots \rightarrow \OO_{X_n} \otimes
 \bigwedge^i \h^*
\rightarrow \cdots \rightarrow \OO_{X_n} \rightarrow
\frac{\OO_{X_n}}{\OO_{X_n} K}\rightarrow 0.
\end{equation*}
Tensoring this by $\mu$, applying $\rho_*$ and taking $\WW$-fixed
points gives an exact sequence
\begin{equation}
\label{kosres3}
0\rightarrow (\PP \otimes
\bigwedge^{n-1}\h^* \otimes \mu)^{\WW} \rightarrow \cdots
\rightarrow
(\PP \otimes \bigwedge^i \h^* \otimes \mu)^{\WW} \rightarrow
\cdots \rightarrow (\PP \otimes \mu)^{\WW} \rightarrow
\frac{\PP_\mu}{(\PP K)_\mu}
\rightarrow 0.
\end{equation}
As this resolution has length $n-1$, the new
intersection theorem \cite[Corollary~5.2]{BKR} implies that the irreducible
components of the support of
the cohomology sheaves   have codimension at most $n-1$ in
$\hi$. In other words, the irreducible components of
$\supp \PP_{\mu}/(\PP K)_{\mu}$ have dimension at least $n-1$.
Since $\supp \PP_{\mu}/(\PP K)_{\mu}\subseteq Z$, the result therefore
 follows from Lemma~\ref{Z-variety1}.
\end{proof}

\subsection{Lemma} \label{supp-punc} {\it The multiplicity
  $\mult_{Z_{(n)}} \left( {\PP}/{\PP K}\right) $ equals $ 1.$}

\begin{proof}
Let $\PP_1$  denote the Procesi bundle on $\hin$, as in the proof
of Lemma~\ref{punctcoh}. Then $\PP/\PP K = \PP_1 /(\PP_1 K +
\PP_1(\mathbf{z},\mathbf{z^*}))$ under the embedding $\cxy
\hookrightarrow \cxy[{\bf z}, {\bf z^*}] = \C[\C^{2n}]$ of
\eqref{hin-defn}.

For a partition  $\lambda$ of $n$, recall  that  the monomial
 ideal $I_{\lambda}$ from \eqref{I-eta-defn}
is a point of $\hio$. We claim, as  is well-known,
that $I_{\lambda}$ is also a point of $Z_{\lambda}$.
To see this, suppose that
$\sum_{i=1}^r \lambda_i [q_i] \in S_{\lambda}\C^2 \cap \h$.
Then the ideal
$$I (\lambda,\mathbf{q})= \bigcap_{i=1}^r   \left(
y^{\lambda_{i}} \C[x,y] + (x-q_{i}) \C[x,y]\right) $$
 clearly   belongs to $Z_{\lambda}$. Since $S_\lambda
\C^2$ is stable under the $\TT^2$ action of \eqref{toract1}, so is
$Z_\lambda$. Thus, $I_\lambda= \lim_{t\to 0}I(\lambda,t\mathbf{q})
$ is also in $Z_\lambda$.

Since $I_{\lambda}$ is killed by $\mathbf{z}$,
   the geometric fibres of
$\PP/\PP K$ and $\PP_1/(\PP_1 K + \PP_1\mathbf{z^*})$ at
$I_{\lambda}$ are equal.   Let
$\TT^2$ act on $\C^2$ as in \eqref{toract1}, and use the
conventions of \eqref{bi-poincare} for the corresponding bigraded
decomposition $\PP_1 (I_{\lambda})=\sum_{ij} \PP_1
(I_{\lambda})^{ij}$. These bigraded components are $\WW$-modules;
indeed, by \cite[(46) and Proposition~3.4]{hai1}
(which is proved in \cite[Section~3.9]{hai3}), the bigraded $\WW$-equivariant
structure of $\PP_1 (I_{\lambda})$ is given by $[\PP_1
(I_{\lambda})] = \sum_{\mu} s^{n(\lambda)}
K_{\mu\lambda}(t,s^{-1})[\mu]$. Factoring $\PP_1$ by $\PP_1
K + \PP_1\mathbf{z^*}$ kills precisely the elements with positive
$\h$-degree, and so the bigraded $\WW$-structure of
$\PP_1(I_{\mu})/(\PP_1 K + \PP_1\mathbf{z^*})(I_{\mu})$ is given
by $\sum_{\lambda} s^{n(\lambda)}K_{\mu \lambda}(0,s^{-1}) [\mu]$.
Specialising $s$ to $1$ and using \eqref{kostka-defn} shows that
the $\WW$-equivariant structure of these fibres is given by
\begin{equation}
\label{GHkos} \frac{\PP}{\PP K}(I_{\lambda}) =
\frac{\mathcal{P}_1}{( \mathcal{P}_1 K+
\PP_1\mathbf{z^*})}(I_{\lambda}) = \sum_{\mu} {K}_{\mu \lambda}
[\mu].
\end{equation}
 In particular \eqref{kostka-defn} implies that
$\dim \, (\PP/\PP K)(I_{(n)}) = 1.$
As the dimension of the geometric fibre at $I_{(n)}$ is an upper bound for the
multiplicity of $\PP/\PP K$ along $Z_{(n)}$, this implies that $\mult_{Z_{(n)}}
 \PP/\PP K\leq 1.$

On the other hand,  by construction,  $\PP_{(n)} = \OO_{\hi}$
 and so there  is an inclusion
$$\OO_{Z_{(n)}} = \OO_{\hi}/(K\cap \OO_{\hi}) \subseteq \PP/\PP K.$$
Since $Z_{(n)}$ is an irreducible component of $Z$ it follows that
$1 = \mult_{Z_{(n)}} \OO_Z \leq \mult_{Z_{(n)}} \PP/\PP K.$
\end{proof}

\subsection{} We can now describe the characteristic cycle of any
standard module in terms of the varieties  $Z_\lambda$ and
  the Kostka numbers $K_{\mu\lambda}$
from \eqref{kostka-defn}.

\begin{thm} \label{stand-cc} Let $\Delta_c(\mu)$ be the standard
$H_c$-module corresponding to $\mu\in\irr\WW$. Then
$$ \Ch (\Delta_c(\mu))   = \sum_{\lambda} K_{\mu
\lambda} [Z_{\lambda}].$$
\end{thm}

\noindent
{\bf Remarks.} (1)
 By \eqref{kostka-defn},  $K_{\mu \lambda} =0$ unless $\mu\geq
\lambda$ and  $K_{\mu \mu}=1$. Combined with the theorem, this therefore shows
that $\mu$ is determined by   $\Chr  \Delta_c(\mu)  $.

(2)  \label{GG-question} Very recently, Gan and Ginzburg \cite{GG}
have used $\mathcal{D}$-modules to give another way of defining
the characteristic cycle of a category $\OO_c$-module as a
subscheme of $\hi$. It would be interesting to know the
relationship between their cycle and ours.

\begin{proof}
By Lemma~\ref{st-equi} it remains to determine the multiplicity of
$\PP_\mu/(\PP K)_\mu$  along each
$Z_{\lambda}$. We will do this by reducing to the case of
Lemma~\ref{supp-punc}.

We first need to introduce some notation. Let $m\in \NN$ and  write
$\widehat{\tau}_m : \text{Hilb}^m\C^2\to \C^{2m}/\mathfrak{S}_m$ for the crepant
resolution. The Procesi bundles for $\text{Hilb}(m)$ and $\text{Hilb}^m\C^2$ will be
written as $\PP(m)$ and $\PP_1(m)$ respectively, whilst the analogue of $K$
will be written $K(m)$. For $\alpha\in \C$, set
$$\overline{\PP}_{1,\alpha}(m) =
\frac{\PP_1(m)}{\PP_1(m)K(m) + \PP_1(m)({\bf z}-\alpha, {\bf
z^*})}
 \qquad \text{and}
\qquad \overline{\PP}(m) = \overline{\PP}_{1,0}(m) = \frac{\PP(m)}{\PP(m) K(m)}.$$
As in \cite[Definition~3.2.4]{hai3} or  \cite[(4.2)]{GS}, the
{\it isospectral Hilbert scheme}\label{isospec-defn} $\mathbb{X}_m$ (which is
the $\C^{2m}$ analogue of $X_n$)  is
described as a set by
$$\mathbb{X}_m= \left\{ (I, p_1,\ldots , p_m):
I\in \text{Hilb}^m\C^2, \, p_1, \ldots, p_m\in \C^2, \,
 \supp I = \sum [p_i]\right\}.$$

 We have
$\mult_{Z_{\lambda}} \overline{\PP}(n) =
\mult_{Z_{\lambda}} \overline{\PP}_{1,0}(n)$.
If $I$ is a generic point of $Z_{\lambda}$ then $\mult_{Z_{\lambda}} \overline{\PP}_{1,0}(n) =
\mult_{I} \overline{\PP}_{1,0}(n)$, the multiplicity of
$\overline{\PP}_{1,0}(n)$ at $I$ in $\hin$. Therefore
\begin{equation} \label{kenisearly}
\mult_{Z_{\lambda}} \overline{\PP}(n) = \mult_{I} \overline{\PP}_{1,0}(n).
\end{equation}
Write $\supp I = \sum [p_i] = \lambda_1[q_1] +\cdots \lambda_r[q_r]$ for
distinct $q_1,\ldots ,q_r\in C = \C\times \{ 0\}$. By
\cite[Lemma~3.3.1]{hai3} and induction there is a neighbourhood of
$(I,p_1,\ldots , p_n) \in \mathbb{X}_n$ that is isomorphic to a neighbourhood
of the point
$$\big((I_1,q_1,\ldots ,q_1),(I_2, q_2,\ldots ,q_2),\ldots ,
 ( I_r, q_r,\ldots ,q_r)\big)\ \in \ \mathbb{X}_{\lambda_1}\times
\mathbb{X}_{\lambda_2}\times
\cdots \times \mathbb{X}_{\lambda_r}.$$
Here $I_1, \ldots, I_r$ are the primary
ideals with support $\lambda_1[q_1], \ldots,
\lambda_r[q_r]$ such that $I_1\cap \cdots \cap I_r = I$. It follows
that
\begin{eqnarray*}(\PP_1)_I = (\rho_{\ast} \OO_{\mathbb{X}_n})_I
\ = \
\bigoplus (\OO_{\mathbb{X}_n})_{(I,p_1,\ldots ,p_n)} & = &
\bigoplus (\OO_{\mathbb{X}_{\lambda_1}})_{(I_1,q_1, \ldots ,q_1)}
\otimes \cdots \otimes
(\OO_{\mathbb{X}_{\lambda_r}})_{(I_r,q_r, \ldots, q_r)} \\ & = &
\bigoplus \bigotimes (\PP_{1, q_j}(\lambda_j))_{I_j},
\end{eqnarray*}
where the
summations are over all distinct orderings of $p_1, \ldots , p_n$. As the
stabiliser of $(p_1,\ldots , p_n)\in \C^{2n}$ is
isomorphic to the Young subgroup
$S_{\lambda}=S_{\lambda_1}\times \cdots \times S_{\lambda_r}$,
 there are $n!/|S_{\lambda}|$ different orderings. Hence
 \eqref{kenisearly} implies that
 \begin{equation}\label{kenhasgone}
\mult_{Z_{\lambda}} \overline{\PP}(n) = \mult_{I} \overline{\PP}_{1,0}(n) =
\frac{n!}{|S_{\lambda}|} \prod_{j=1}^r \mult_{I_j}
\overline{\PP}_{1, q_j}(\lambda_j),
\end{equation}
where the multiplicities of the final term are taken at $I_j$ in
$\text{Hilb}^{\lambda_j}\C^2$.

In order to calculate the right hand side of \eqref{kenhasgone}, we need to
assume that each
$I_j$ is a generic point of the scheme
$V_j =\widehat{\tau}_{\lambda_j}^{\,-1}(\lambda_j[q_j])\subseteq
\text{Hilb}^{\lambda_j}\C^2$.
  We will confirm  by a dimension count that this can be achieved.
Since $\dim V_j=\lambda_j-1$, we have
$\dim (V_1\times \cdots \times V_r)
=\sum_{j=1}^r (\lambda_j -1).$
 We may choose the $r$-tuple $(q_1,\ldots ,q_r)$ from a
non-empty open subset of
$\{ (q_1,\ldots ,q_r): \sum_{j=1}^r \lambda_jq_r =0\}$.
Then as we vary the $q_j$'s in this parameter space the
spaces $ V_1 \times \cdots \times
V_r$ sweep out a variety of dimension
 $\sum_{j=1}^r (\lambda_j -1 ) +(r-1) = n-1 =
\dim Z_{\lambda}$.  We conclude that $I$ can be chosen so that the $I_j$'s are
generic.  It follows that
$\mult_{I_j} \overline{\PP}_{1,q_j}(\lambda_j) =
\mult_{V_j} \overline{\PP}_{1,q_j}(\lambda_j).$
Shifting the base from $q_j$ to $0$
and writing $Z_{(m)}=\widehat{\tau}_m^{\,-1}(\mathbf{0})$ for
 the punctual Hilbert scheme in $\mathrm{Hilb}^m\C^2$  shows that
$$\mult_{V_j} \overline{\PP}_{1,q_j}(\lambda_j) =
\mult_{Z_{(\lambda_j)}} \overline{\PP}_{1,0}(\lambda_j) =
\mult_{Z_{(\lambda_j)}} \overline{\PP}(\lambda_j).$$
By Lemma~\ref{supp-punc} this final multiplicity
equals 1. Combined with \eqref{kenhasgone} this gives
$\mult_{Z_{\lambda}} \overline{\PP}(n) = {n!}/{|S_{\lambda}|}$.

Finally we compute $\mult_{Z_{\lambda}} \overline{\PP}(n)_{\mu} $.
As the dimension of geometric fibres bounds multiplicity from
above, and $I_{\lambda}\in Z_{\lambda}$, \eqref{GHkos} shows that
\begin{equation}
\label{nice}
\mult_{Z_{\lambda}} \overline{\PP}(n)_{\mu}
\ \leq\ \dim\, \overline{\PP}(n)(I_{\lambda})_\mu =  K_{\mu \lambda}.
\end{equation}
Hence
$$ \frac{n!}{|S_{\lambda}|} \ = \ \mult_{Z_{\lambda}} \overline{\PP}(n)\  =\
\sum_{\nu} \mult_{Z_{\lambda}} \overline{\PP}(n)_{\nu} \dim \nu\  \leq\
\sum_{\nu} K_{\nu \lambda} \dim \nu .$$
By   \cite[Remark following~I(7.8)]{MacD},
$\sum_{\nu} K_{\nu \lambda} \dim \nu =
\dim \text{Ind}_{S_{\lambda}}^{\WW} \C =  {n!}/{|S_{\lambda}|}$
 and so  \eqref{nice} is an equality for each $\mu$.
 \end{proof}

\subsection{} It would be interesting to know if there is a result analogous
to  Theorem~\ref{stand-cc} that holds for the simple modules
$L_c(\mu)$, but we are unable to prove this. Indeed we do not even
know  whether the characteristic varieties $\Chr L_c(\mu)$ are
equi-dimensional; this is a special case of of
Question~\ref{pure-question}.
 The theorem does however
provide information  about the restricted characteristic varieties $\Chw
L_c(\mu)$, as defined in \eqref{cyclew-defn}.

\begin{cor}\label{simp-cc}
For all $\mu\in \irr \WW$ there exist integers
$a_{\lambda, \mu}\in \NN$ such that
\begin{equation}\label{simp-cc1}
\Chw (L_c(\mu)) \ = \ [Z_{\mu}] + \sum_{\lambda < \mu} a_{\lambda, \mu}
 [Z_{\lambda}].\end{equation}
\end{cor}

\begin{proof}
Recall from Lemma~\ref{Z-variety1} that each $Z_\lambda$ has dimension $n-1$.
As in \eqref{dominance},
the  partition $(1^n)$ corresponding to the $\sign$ representation of
$\WW$ is minimal in the dominance ordering. By \cite[Remark~3.6(2)]{GS}
$L_c(\sign) = \Delta_c(\sign)$. Hence
$\Chw L_c(\sign) =
 \sum_{\lambda} K_{(1^n)\lambda} [Z_{\lambda}] = [Z_{(1^n)}],$
  by Theorem~\ref{stand-cc}  and \eqref{kostka-defn}.
This begins an induction.

Now fix a partition $\mu$ and suppose that \eqref{simp-cc1} holds
for all partitions $\nu < \mu$. By \cite[Lemma~3.5 and
Remark~3.5(1)]{GS}, there is a short exact sequence $0\to M \to
\Delta_c(\mu) \to L_c(\mu) \to 0$ where the simple composition
factors of $M$ are of the form  $L_c(\nu)$ for $\nu < \mu$. By the
inductive hypothesis and Lemma~\ref{ses-cc}(c) $\Chw M = \sum_{\nu
< \mu} b_{\nu} [Z_{\nu}]$ for some $b_\nu\in \NN$.
 On the other hand,
by Theorem~\ref{stand-cc} and \eqref{kostka-defn},
$\Chw \Delta_c(\mu) = [Z_\mu] + \sum_{\lambda < \mu} K_{\mu \lambda} [Z_\lambda].$
Applying Lemma~\ref{ses-cc}, again, shows that
$$\Chw L_c(\mu) = \Chw \Delta_c(\mu) - \Chw M =
[Z_{\mu}] + \sum_{\lambda < \mu} a_{\lambda, \mu}[Z_\lambda],$$
for some $a_{\lambda,\mu}\in \NN$.
This completes the induction step and hence the proof of the lemma.
\end{proof}

\subsection{}\label{simp-cc-fd}
  One case where it is easy to identify $\Chw L_c(\mu)$ is when
 $\mu=\triv$ and $L_c(\triv)$ is finite dimensional;
 thus $c=r/n$ with $(r,n)=1$ by \cite[Theorem~1.2]{BEGfd}.
 In this case,  $\Chro L_c(\triv)$ is
the zero orbit $\mathbf{0}$  in $\h\oplus \h^*/\WW$ and so
 Proposition~\ref{compare-char}
 shows that $\Chr L_c(\triv)\subseteq \tau^{-1}(\mathbf{0}) =
 \hio$. By Corollary~\ref{simp-cc}, this forces
 $\Chr L_c(\triv)=\Chw L_c(\triv)=[\hio]$.

\subsection{}\label{borel}
Let $H_{2n-2}(Z, \C)$ be the top degree Borel-Moore homology
group of the variety $Z$. By Lemma~\ref{Z-variety1}
and \cite[Section~8.2]{Nak},
 the irreducible components  $Z_\lambda$ of $Z$ form a basis of
$ H_{2n-2}(Z,\C)$, so its dimension is equal to the number of partitions of
$n$.
\begin{cor}
Taking the reduced characteristic cycle $\Chw M$ of a module $M\in \OO_c$
 induces an isomorphism of vector spaces
 $\Chw : K(\OO_c)\otimes_{\Z} \C \xrightarrow{\sim} H_{2n-2}(Z, \C).$
\end{cor}
\begin{proof}
We first check that $\Chw$ is well-defined.
So for any short exact sequence
$0\to M_1 \to M_2 \to M_3 \to 0$ in
$\OO_c$, we must show that
$\Chw M_2 = \Chw M_1 + \Chw M_3$. Since every object in $\OO_c$ has a
finite composition series, it is enough by Lemma~\ref{ses-cc}(c) to
prove that $\dim \Chr L_c(\mu)=n-1$ for all $\mu\in\irr\WW$.
But this is immediate from Corollary~\ref{simp-cc}.

As $K(\OO_c)\otimes_{\Z}\C$ and $H_{2n-2}(Z,\C)$
 are vector spaces with bases $\{
[L_c(\mu)] \}$ and $\{[Z_\mu]\}$,  respectively,  Corollary~\ref{simp-cc} also
shows that the matrix which represents $\Chw$ for these choices of bases is
unitriangular when we order by the dominance ordering. As such a matrix is
invertible, $\Chw$ is an isomorphism. \end{proof}

\subsection{}\label{justify-1.15} In this subsection we give a
refinement of  Corollary~\ref{borel} that also justifies the comments
made in \eqref{intro-1.15} from the introduction.

 According to
\cite{gro} and \cite[(9.13)]{Nak} there is an isomorphism
$$\xi : \bigoplus_{n\geq 0} H_{2n-2}(Z(n),\C)\xrightarrow{\sim}
\C[p_k : k\in \NN]$$ between the cohomology of the Lagrangian
subvarieties $Z(n) \subseteq \hi$ and the space of symmetric
functions. Under this identification $[Z_{\lambda}]$ is sent to
the monomial symmetric function $m_{\lambda}$,
\cite[Corollary~9.15]{Nak}. Thus by Theorem~\ref{stand-cc} and
\eqref{kos-tran}
\begin{equation}\label{justify-less}
\xi (\Chw  \Delta_c(\mu)) = \sum_{\lambda} K_{\mu\lambda}
m_{\lambda} = s_{\mu},\end{equation} the Schur function of $\mu$.

For generic $c$, $\OO_c$ is semisimple and
so this gives a natural description of  $K(\OO_c)$  in terms of
symmetric functions. When $\OO_c$ is not semisimple this is
insufficient, and this is the case we examine here.
By \cite[Corollary~2.11]{BEGqi} and \cite{DU},  $\OO_c$ is not semisimple
 if and only if $c\in \mathcal{C}$, as defined in \eqref{main-hyp}. We
 fix such a value of $c$ and {\it  we assume that $\OO_c$
 is equivalent to  $S_q\md$} where $S_q=S_q(n,n)$ is the
{\it $q$-Schur algebra} \label{schur-defn}  and
$q= \exp(2\pi i c)$ is a primitive $e$th root of unity. As mentioned in the
introduction, this is conjectured in \cite[Remark~5.17]{GGOR} and a proof has
recently been announced by  Rouquier.

Under this assumption, and in the notation of \cite[(3.5)]{GS},
there is an equality of multiplicities
\begin{equation} \label{schurmor}
[\Delta_c(\lambda) : L_c(\tau)] = [W_q(\lambda) : F_q(\tau)].\end{equation}
The multiplicities on the right can be described as follows, where we
refer the reader to \cite{lecthi} for all definitions and notation. Let
$\mathcal{F}_v$ be the $v$-deformed Fock space representation of
$U_v(\widehat{\mathfrak{gl}}_e)$, a deformation of the space of symmetric
functions in an indeterminate $v$.  The standard basis
$\{| \lambda \rangle \}$ of $\mathcal{F}_v$ indexed by \textit{all}
partitions corresponds to the basis of symmetric functions given by the
Schur functions. By \cite[Theorem~4.1]{lecthi}, there exist upper and
lower canonical bases of $\mathcal{F}_v$, written $\{G(\lambda)\}$ and
$\{ G^{-}(\lambda)\}$ respectively. If we write
$$G(\mu) = \sum_{\lambda} d_{\lambda\mu}(q) |\lambda \rangle
\qquad\mathrm{and}\qquad  G^-(\lambda) = \sum_{\mu} e_{\lambda\mu}(q) |\mu\rangle,$$
then $d_{\lambda\mu}(q)\in \Z[q]$ and $e_{\lambda\mu}(q) \in \Z[q^{-1}]$.

If we consider only partitions of $n$, then by \cite[Section~4]{lecthi}
these bases
 are related by the matrix equality
\begin{equation}  \label{lt1} [e_{\lambda\mu}(q^{-1})] =
 [d_{\lambda^t\mu^t}(q)]^{-1}.\end{equation}
The confirmation of \cite[Conjecture~5.2]{lecthi} in \cite{vasvar}
shows that \label{lt2}$d_{\lambda^{t} \mu^t}(1) =
[W_q(\lambda) : F_q(\mu)].$
Combined with \eqref{schurmor}  and \eqref{lt1}, this implies  that
$\Chw (L_c(\lambda)) = \xi^{-1}(\sum_{\mu} e_{ \lambda \mu} (1)
s_{\mu}) = \xi^{-1}(G^-(\lambda)_{v=1}).$
In other words,

\begin{prop}\label{justify-1.15-2}   Assume that $\OO_c$ is   equivalent to
$S_q\md$ and pick $\lambda\in \irr\WW$. Then  the reduced
characteristic cycle $\Chw (L_c(\lambda))$ is exactly the
canonical basis element $\xi^{-1}(G^-(\lambda)_{v=1})$   in
$H_{2n-2}(Z,\C)$ corresponding to the lower canonical basis
element $G^-(\lambda)$ of the ring of symmetric functions from
\cite{lecthi}.
 \end{prop}

\noindent {\bf Remark.} Thanks to \cite{vasvar}, the polynomials
$d_{\lambda \mu}(q)$ are parabolic Kazhdan-Lusztig polynomials.
This hints at the possibility that the rational Cherednik algebra
connects Macdonald theory with Kazhdan-Lusztig theory by providing
a bridge between $q$-Schur algebras and Hilbert schemes. As an
example of this, Theorem~\ref{stand-cc} and \eqref{schurmor} give
a number of different positive integral factorisations of the
Kostka matrix in terms of evaluations of various Kazhdan-Lusztig
polynomials.

\subsection{An example}\label{s2-eg}
 The results of this section
are nicely illustrated by the case $n=2$ and $c=\frac{1}{2}$. (By
Remark~\ref{main-hyp}(1)  all the results of this paper do hold in
this case.) For these values of $n$ and $c$,
\cite[Proposition~8.2]{EG} implies that $U_c=\C\langle
x^2,xy,y^2\rangle \cong U(\mathfrak{sl}_2(\C))/(\Omega)$, where
$\Omega $ is the Casimir element. Moreover,  $M=e\Delta_c(\triv)$
 identifies with the non-simple Verma module over
$U(\mathfrak{sl}_2(\C))/(\Omega)$ while
 $N=e\Delta_c(\sign)$ embeds into $e\Delta_c(\triv)$ with factor
 $L_c(\triv)\cong \C$.

  First consider $M$ and note that $M\cong \C[x^2]$
as a module over $\C[x^2]=\C[\h]^\WW$.
  Give $M$ the induced $\ord$ filtration $\Lambda$: thus $M=\Lambda^0M$.
  Now, $M(1)=eH_{c+1}\delta e\otimes M \cong e\Delta_{c+1}(\triv)\cong \C[x^2]$,
   by \cite[Proposition~3.16]{GS}. However, the tensor product filtration
on $M(1)$
   is a little subtle. Indeed,  since $\delta=x$ it is clear that
   $$\Lambda^0M(1)\supseteq e\C[x^2]x\delta e \otimes 1
   =\C[x^2]x^2,$$
   but one can check that the first occurrence of $1$ in the filtration
 comes from
   $ey\delta e\otimes 1 =-2e\otimes 1.$
   Thus, under the tensor product filtration $\Lambda$, one has
   $\Lambda^0M(1)=\C[x^2]x^2$ and
   $\Lambda^1M(1)/\Lambda^0M(1)\cong \C$.
   We leave it to the reader to prove, by induction, that

   \begin{lem}\label{delta-triv} As above, take $n=2$, $c=\frac{1}{2}$ and let
   $M=e\Delta_c(\triv)\supset N=e\Delta_c(\sign)$. Pick
$k\geq 1$.
\begin{enumerate}
\item
   The tensor product filtration $\Lambda$  on $M(k)=B_{c+k,c}\otimes M$
   satisfies $\Lambda^0M(k)=\C[x^2]x^{2k}$ and
   $x^{2k-2i}\in \Lambda^iM(k)\smallsetminus \Lambda^{i-1}M(k)$
   for $0\leq i\leq k$.

\item $\gr_{\Lambda}\widetilde{M} = D_1\oplus D_2$,
   where $D_1\cong \C[x^2]x^2\oplus \C[x^2]x^4\oplus \cdots$
 and $D_2\cong \C\otimes_{A^0}A $.  Consequently,
$\Chw M=\Chr M=\V_1+\V_2$,
   where $\V_1=Z_{(1^n)}$ and $\V_2=Z_{(n)}$
   is the punctual Hilbert scheme. \qed
   \end{enumerate}
   \end{lem}

Note that part (2) of the lemma is in accordance with
   Theorem~\ref{stand-cc}. We emphasise that
   ${V}(\gr\widetilde{M}) $   has two components whereas
${V}(\gr\widetilde{N}) $ has just one. In contrast,
 under the $\ord$ filtration  $M$ and
$N$ both have associated graded module isomorphic to
 $\C[x^{2}]$ and so  $\Chro M = \Chro N$.

      \subsection{}\label{s2-eg2}
 Keep the notation and hypotheses from \eqref{s2-eg}.
  Proposition~\ref{1d} and Theorem~\ref{stand-cc} imply that
  $L=eL_c(\triv)$ satisfies $\gr_\Lambda\widetilde{L} =D_2$ with
   $\Chw L=
Z_{(n)}$, while  $N=e\Delta(\sign)$ satisfies $\gr_\Lambda
\widetilde{N} =D_1$ with
$\Chw  N = Z_{(1^n)}$. This also follows from the
computations in \eqref{s2-eg}.  Indeed, give $N$ the  filtration
induced from that on $M$. Then  $N=\C[x^2]x^2$ and
$N(k)=e\Delta_{c+k}(\sign)=B_{c+k,c}\otimes N = \C[x^2]x^{2k+2}$.
 Repeating the argument from \eqref{s2-eg} then shows that
 $N(k)=\Lambda^0N(k)$ for all $k$, whence the result.

  \subsection{}\label{s2-eg3}
  There are a couple of consequences of these computations that are
  worth noting. First, for a $U_{c}$-module $P$, even one in
  category $\mathcal{O}_c$,
   it is important to use the tensor
  product filtration on $\widetilde{P}$, even if the $\ord$ filtration
  would appear more natural. For example, suppose that one
  begins with $M=e\Delta_c(\triv)$ as above and gives each $M(k)
  \cong \C[x^2]$ the $\ord$ filtration $\Gamma$; thus $\Gamma^0M(k)
  =M(k)$ for each $k$.
  Then $\gr_\Gamma \widetilde{M}  = \C[x^2]\oplus \C[x^2]\cdots$.
But this module has a different characteristic variety to the variety
$\Chr M$   computed   in \eqref{s2-eg2}. This does not contradict the earlier
results since an easy variant of the proof of  Lemma~\ref{ch-is-wd}
  can  be used to show  that  $\gr_\Gamma \widetilde{M} $  cannot  finitely
  generated as an $A$-module.

Second, recall from \eqref{warning-subsec}  that
the  tensor product and $\ord$ filtrations on the $B_{ij}$ are
equal. The  last paragraph shows that
the analogous result fails for  $M=eL_c(\triv)$.

%%%%%%%%%%%%%%%%%%%%%%%%%%%%%%%%%%%%%%%%%%%%%
\section{Bimodules}\label{bimodule-sect}

\subsection{} In this short section we show how to extend our
results on the characteristic  varieties of one-sided modules  to
$H_c$- and $U_c$-bimodules, and we discuss the resulting
characteristic varieties of Harish-Chandra bimodules.
Since it requires no extra work,
we will actually  prove these results for $(U_c,U_d)$-bimodules
and so we will
 assume throughout
 that $c,d\in \C$ satisfy the hypothesis of
\eqref{main-hyp}. As usual,   analogues of the results proved  here do hold
for more general values of $c$ and $d$.

\subsection{}\label{bimodthm-sec}
 Let $\omega$ be the anti-automorphism of
$D(\hr)\ast \WW$ defined  by $\omega (x) = x, \ \omega(y) = -y$
and $\omega(w) = w^{-1}$ for $x\in \h^*, \ y \in \h$ and $w\in
\WW$. By \cite[(4.2)]{GGOR},  $\omega$ restricts to an involutive
anti-automorphism on both $H_c$ and $U_c$. Let
$R \boxtimes S$ denote  the usual
external tensor product of two $\C$-algebras $R$ and $S$.
 Then we can identify
the category $(U_c,U_d)\bmd$ of finitely generated
$(U_c,U_d)$-bimodules with the category of
finitely generated left $U_c\boxtimes  U_d$-modules: for any
such bimodule $M$ and  $m\in M$,  we set
$ (p\boxtimes q)\cdot m = p m \omega(q)$
for $p\in U_c$ and $q\in U_d$.

We have a straightforward generalisation of
Theorem~\ref{intro-1.2}.
\begin{thm}
\label{bimodthm} If $c,d $ satisfy Hypothesis~\ref{main-hyp}
there exists a filtered $\Z$-algebra $C$ such
that
\begin{enumerate}
\item $(U_c,U_d)\bmd$ is equivalent to $C\lqgr$;

\item  $\gr C$, the associated graded ring  of $C$, is isomorphic
to (the $\Z$-algebra associated with) the homogeneous coordinate
ring $ \bigoplus_{k \geq 0} \mathrm{H}^0 \left(\hi\times \hi, \
\pi_1^{\ast}\LL^k \otimes \pi_2^{\ast}\LL^k\right),$ where
$\pi_j$ denotes the projection map onto the $j^{\mathrm{th}}$
component of $\hi\times \hi$.
\end{enumerate}
\end{thm}
\begin{proof}
As above, we replace $(U_c,U_d)\bmd$ by $U_c\boxtimes U_d\lmod$. By
\eqref{shift-defn}, the $(U_{c+1}\boxtimes U_{d+1}, U_c\boxtimes
U_d)$-bimodule $eH_{c+1}\delta e\boxtimes eH_{d+1}\delta e$
induces a Morita equivalence between $U_c\boxtimes U_d\lmod$ and
$U_{c+1}\boxtimes U_{d+1}\lmod$. Hence,
defining $C$ to be the Morita $\Z$-algebra associated to the data
$$\{ U_{c+i}\boxtimes U_{d+i}, \   eH_{c+i+1}\delta e\boxtimes
eH_{d+i+1}\delta e \ ; \ i\in \NN\}$$ and applying
Lemma~\ref{zalgex1} gives the first part of the theorem.

In the notation of \eqref{Mij-defn},   $C_{ij} = B_{ij}(c) \boxtimes B_{ij}(d)$
for all $i\geq j\geq 0$. If we give $C$ the componentwise tensor
product filtration,
$$F^k C_{ij} = \sum_{r+s \leq k} \ord^r B_{ij} \boxtimes \ord^s
B_{ij},$$ then we deduce from \cite[Theorem~6.4(2)]{GS} that $$\gr
C =\bigoplus_{i\geq j\geq 0} \gr B_{ij} \boxtimes \gr B_{ij} \cong
\bigoplus_{i\geq j \geq 0} \mathrm{H}^0(\hi , \ \LL^{i-j})
\boxtimes \mathrm{H}^0 (\hi, \ \LL^{i-j}).$$ The
 right hand side of this
equation is the $\Z$-algebra associated to the graded algebra
$$\bigoplus_{k\geq 0} \mathrm{H}^0(\hi , \ \LL^{k}) \boxtimes
\mathrm{H}^0 (\hi, \ \LL^{k}) = \bigoplus_{k \geq 0} \mathrm{H}^0
(\hi \times \hi, \ \pi_1^{\ast} \LL^k \otimes \pi_2^{\ast}
\LL^k),$$ which proves part (2) of the theorem.
\end{proof}

\subsection{}\label{upsilon-subsec}
Write $(U_c,U_d)\fbmd$ for the category of $(U_c,U_d)$-bimodules
$M$ with a good filtration $\Lambda$;
by definition this means that
$( \ord^i U_c)(\Lambda^jM )(\ord^k U_d)\subseteq \Lambda^{i+j+k}M$, for all $i,j,k$
and that the associated graded object  $\gr_\Lambda M$ is finitely generated
as an $(\ogr U_c, \ogr U_d)$-bimodule. Under the identification of $(U_c,U_d)\bmd$ with
 $U_c\boxtimes U_d\md$, this category is simply $U_c\boxtimes U_d \fmd$.
Thanks to Theorem~\ref{bimodthm}, the requirements
(1)--(4) of \eqref{filter-discuss} are fulfilled, and so we have a
natural functor
$$\Upsilon' : (U_c, U_d) \fbmd \longrightarrow \coh (\hi \times \hi)\qquad
(M,\Lambda) \mapsto \pi \gr_\Lambda M.$$

We wish to adjust $\Upsilon'$ to take account of the twisting by
$\omega$ which occurs when we pass from $(U_c,U_d)\bmd$ to
$U_c\boxtimes U_d\lmod$. Note first that the $\ord$ filtration is
preserved by $\omega$, so we obtain an induced action of $\omega$
on $\cxy\subset \ogr D(\hr)\ast \WW$ given by the formul\ae\ in
\eqref{bimodthm-sec}. Clearly this is also induced by the
automorphism  $\bar{\omega}$  of  $\h \oplus \h^*$ which is the
identity on $\h$ and multiplies $\h^*$ by $-1$. Then
$\bar{\omega}$ induces automorphisms of $\hi$ and $\h\oplus \h^*
/\WW$ which we also denote by $\bar{\omega}$. We set $$ \Upsilon =
(1\times \bar{\omega}^{\ast} ) \circ \Upsilon ' : (U_c,U_d)\fbmd
\longrightarrow \coh (\hi \times \hi).$$

Following the discussion in Section~\ref{zalg},
 we define {\it the characteristic variety}\label{bi-charvar}
 $\Chr M\subseteq \hi \times \hi$  of a finitely generated
$(U_c,U_d)$-bimodule $M$ to be the  support variety of $\Upsilon_\Lambda
M$, for any good filtration $\Lambda$. Similarly,  {\it the
associated variety} $\Chro M\subseteq (\h\oplus \h^*)/\WW\times
(\h\oplus \h^*)/\WW$  is the support variety of $(1\times
\bar{\omega}^*) (\gr_\Lambda M)$.

\subsection{Harish-Chandra bimodules} \label{subsec-hc}
Under the canonical embeddings of $\cx^\WW$
into $U_c$ and $U_d$ any $(U_c,U_d)$-bimodule $M$ inherits
 the structure of a $\cx^\WW$-bimodule. Similarly,
$M$ is a $\C[\h^*]^\WW$-bimodule.
Following \cite[Definition~8.8]{BEGqi}  a {\it Harish-Chandra
$(U_c,U_d)$-bimodule}\label{HC-defn} can therefore be  defined to
be a finitely generated $(U_c,U_d)$-bimodule on which the induced
adjoint actions
 of    $\cx^{\WW}$ and $\cy^{\WW}$ are locally
nilpotent. We let $\hcdc$ denote the category of Harish-Chandra
$(U_c,U_d)$-bimodules, a full subcategory of $(U_c,U_d)\bmd$.

\begin{prop}\label{stein}
If $M\in \hcdc$, then $\Chr M \subseteq \hi
\times_{\h\oplus \h^*/\WW} \hi.$
\end{prop}

\noindent
{\bf Remark.} The variety $\hi \times_{\h\oplus\h^*/\WW} \hi$ is
an analogue of the Steinberg variety
$T^*\mathcal{B}\times_{\mathcal{N}} T^*\mathcal{B}$ which features
prominently in geometric representation theory. We refer the reader to
\cite[Section~3]{BB} for   results on the characteristic varieties of
Harish-Chandra bimodules in the Lie-theoretic context.

\begin{proof}
Set $\Omega= (\h\oplus \h^*)/\WW\times (\h\oplus \h^*)/\WW$ and
let $\Gamma$ be a good filtration on $M$.
We claim that the associated
variety $\Chro M$ is contained in the diagonal $\Delta$ of
$\Omega$. Given the claim,
 Proposition~\ref{compare-char}
implies that   $$\Chr M \
\subseteq\  (\tau\times \tau)^{-1}(\Chro M) \ \subseteq \  (\tau\times
\tau)^{-1}(\Delta)\  =\   \hi \times_{\h\oplus\h^*/\WW} \hi,$$
where $\tau : \hi \to \h\oplus
\h^*/\WW$ is the usual resolution of singularities. Thus the proposition
follows.

Let $I = \mathrm{ann}_{\C[\Omega]}(\gr_\Gamma M)$ and write
$\text{Graph}(\bar{\omega}) \ = \  \{ (z,
\bar{\omega}(z)) : z \in \h\oplus \h^*/\WW \} \ \subset\ \Omega.$
Since the definition of $\Chr_0$ in \eqref{upsilon-subsec}
has a twist by $\bar{\omega}^*$,
 the claim is equivalent to the assertion that
$\mathcal{V}(I) \subseteq \text{Graph}(\bar{\omega})$.
To prove this, note first that, as in \cite[Appendix~D]{EG},
 $\gr U_c\boxtimes \gr U_d$ inherits a Poisson structure
 from multiplication in $U_c\boxtimes U_d$.
  For $a,b,a',b' \in \cxy^{\WW}$ this is
given by
$$\{a\boxtimes b, a'\boxtimes b'\} = \{a,a'\}\boxtimes bb' + aa'
\boxtimes \{b,b'\},$$ where the bracket $\{ \,\,  , \,   \}$  on
$\cxy^{\WW}$  is the   one  arising from multiplication in either
$\ogr U_c$ or $\ogr U_d$. (These two  brackets are the same since they are both
equal to  the one arising from
 the  canonical symplectic structure on $\h\oplus \h^*$.)

Since $M$ is a Harish-Chandra module,
$p\boxtimes 1 - 1 \boxtimes \omega(p)$
 acts locally nilpotently on $  M$ for any
element $p\in \cx^{\WW}$ or $p\in \cy^{\WW}$. We will write
$\hat{p}=\sigma(p)$ for the principal symbol of $p\in U_c$ (or $p\in U_d$). Since
$\gr_{\Gamma} M$ is finitely generated, on passing
 to associated graded modules,
it follows that $\hat{p}\boxtimes 1 - 1 \boxtimes
\omega(\hat{p})\in \sqrt{I}$ for all homogeneous elements $\hat{p}\in \cx^{\WW}$ or
$\hat{p}\in \cy^{\WW}$. As
$\omega$ is an anti-automorphism of $U_d$ we have $\{\omega(\hat{p}),
\omega(\hat{q})\} =
\omega\{\hat{q},\hat{p}\}$ for any homogeneous elements
$\hat{p},\hat{q}\in \ogr U_c=\ogr U_d$ and so
$$\{ \hat{p} \boxtimes 1 - 1\boxtimes \omega(\hat{p}),\,\hat{q}\boxtimes 1 - 1
\boxtimes \omega(\hat{q}) \} \ =\ \{ \hat{p},\hat{q} \}\boxtimes 1 +
1\boxtimes\{ \omega(\hat{p}) ,\omega(\hat{q}) \}
\ =\ \{ \hat{p},\hat{q} \}\boxtimes 1 - 1\boxtimes \omega(\{\hat{p}
,\hat{q}\}),$$
for any such $\hat{p}, \hat{q}$. By Gabber's integrability theorem
\cite[Theorem~1]{GAB} or \cite[Theorem~1.3]{Bj2},
  $\{ \sqrt{I}, \sqrt{I} \} \subseteq
\sqrt{I}$ while, by \cite[Appendix~2]{wal}, $\cxy^{\WW}$ is generated by
the homogeneous elements in
$\cx^{\WW}$, $\cy^{\WW}$ and   their iterated Poisson brackets.
The displayed equation therefore implies
  that $\sqrt{I}$ contains $r\boxtimes 1 - 1 \boxtimes
\omega(r)$ for all $r\in \cxy^{\WW}$. In other words,
 $\mathcal{V}(I) \subseteq \text{Graph}(\bar{\omega})$.  This
completes the proof of the claim, and hence of the proposition.
\end{proof}

\subsection{} The obvious analogue of Proposition~\ref{stein}
also holds for the  category of Harish-Chandra $(H_c,H_d)$-bimodules.
 The reason is that,  thanks to the natural  analogue of
\cite[Theorem~3.3]{GS}, the idempotent functor which sends $M\in
(H_c,H_d)\bmd$ to $eMe\in (U_c,U_d)\bmd$ is an equivalence of categories.

\subsection{Questions}\label{HC-questions}
We end this section with several questions and comments about
Harish-Chandra modules.

 Let $Y$ be the variety $\hi \times_{\h\oplus \h^*/\WW} \hi$
 introduced in Proposition~\ref{stein}. It follows from
 the proof of Theorem~\ref{bkr} that the resolution of singularities
 $\tau: \hi \to \h\oplus \h^*/\WW$ is  semismall in the sense of
 \cite[Definition~6.7]{Nak}.
 Following \cite[Section~6.2]{Nak} we can make this more precise.
First,  $\h\oplus \h^*/\WW$ is stratified as in \eqref{Z-variety}
by the locally closed subvarieties $S_{\lambda} = (\h\oplus
\h^*/\WW) \cap S_{\lambda}\C^2$
 for $\lambda$ a partition of $n$. The
resolution $\tau$ then restricts to a topological fibre bundle
$\tau: \tau^{-1}(S_\lambda) \to S_\lambda$
with fibre equal to the product of punctual Hilbert schemes
$Z_{(\lambda_1)}\times \cdots \times Z_{(\lambda_r)}$, where
$\lambda=(\lambda_1,\dots,\lambda_r)$ has exactly $r$ nonzero parts. If we set
$$Y_{\lambda} = \{ (z, z') \in Y : \tau (z) = \tau(z') \in S_{\lambda}\}$$
then $Y = \coprod_{\lambda} Y_{\lambda}$ and
$\dim Y_{\lambda} = (2r - 2) + (n-r) + (n-r) =
2(n-1)$ for each $\lambda$ with  $r$ nonzero parts.
 Thus $Y$ is pure of dimension $2(n-1)$ with irreducible
components equal to the closures of the various $Y_{\lambda}$.

This raises two questions:
\begin{enumerate}
\item{} {\it Does every  Harish-Chandra bimodule have (restricted)
characteristic cycle equal to a weighted sum of the $Y_\lambda$? }

\item{} {\it If so, what are those weights  for  important
 Harish-Chandra bimodules, for example the
 bimodules $B_{ij}(c)$ from \eqref{Mij-defn}?  Further examples of
 Harish-Chandra bimodules   appear in \cite[Section~8]{BEGqi}.}
\end{enumerate}

\subsection{} Recall from \cite[Section~2.7]{CG} the convolution
product on Borel-Moore homology.
\begin{lem}\label{conv-lemma}
Under the convolution product there is an algebra isomorphism
$$\Xi: H_{4(n-1)}(Y,\C) \cong \text{Rep}(\WW)\otimes_{\Z} \C,$$ where $\text{Rep}(\WW)$
is the representation ring of $\WW$.
\end{lem}
\begin{proof}
This follows immediately from the formalism of
\cite[Section~8.9]{CG}. Indeed, let $x\in S_{\lambda}$. Since
$\tau^{-1}(x)$ is irreducible, the local system on $S_{\lambda}$
associated to $H_{2(n-r)}(\tau^{-1}(x), \C)$ is trivial. Applying
\cite[Corollary~8.9.8 and Proposition~8.9.9(b)]{CG} shows that
$H_{4(n-1)}(Y,\C)$ is a commutative semisimple $\C$-algebra with
basis labelled by partitions of $n$. Since
$\text{Rep}(\WW)\otimes_{\Z} \C$  has a non-degenerate associative
bilinear form given by the inner product on characters is
non-degenerate, it is semisimple, and so the lemma follows.
\end{proof}

A solution to the next problem would provide a significant refinement of the
lemma.
\begin{problem}   Identify the $\TT^2$-equivariant
Grothendieck group $K^{\TT^2}(Y)$ under the convolution product.
\end{problem}

\subsection{} When $c=d$, write
 $\hcdc= \mathcal{HC}_c $.
  It is known that $ \mathcal{HC}_c$ is monoidal
under tensoring over $U_c$, and that $\mathcal{HC}_c $ is
equivalent to $\C \WW\lmod$
 as a monoidal category for $c\in \NN$, \cite[Lemma~8.3 and
Theorem~8.5]{BEGqi}. It acts on $\OO_c$ by tensoring over $H_c$.

\begin{question} Is there a ring isomorphism $\Theta :
K(\mathcal{HC}_c) \longrightarrow \text{Rep} (\WW)$ such that
\begin{enumerate}
\item the following diagram  commutes.
$$
\setlength{\unitlength}{1mm}
\begin{picture}(40,15)
\put(0,13){$K(\mathcal{HC}_c)\otimes_{\Z}\C$}
\put(22,14){\vector(1,0){15}} \put(40,13){$H_{4(n-1)}(Y,\C)$}
\put(26,15){$\Chw$} \put(15,11){\vector(1,-1){7}}
\put(12,6){$\Theta$} \put(45,11){\vector(-1,-1){7}}
\put(45,6){$\Xi$} \put(20,0){$\text{Rep}(\WW)\otimes_{\Z}\C$}
\end{picture}
$$
\item $K(\OO_c)$ corresponds to the regular representation under
the map sending $[\Delta_c(\lambda)]$ to $[\lambda]$?
\end{enumerate}
\end{question}

Implicit in these questions is the assertion  that the simple objects of
$\mathcal{HC}_c$ can be labelled by partitions of $n$; once again, this is
known only for $c\in \NN$ where it follows from the equivalence
of monoidal categories
$\mathcal{HC}_c\simeq \C \WW \lmod$.
 One can ask: {\it Does   such an
equivalence exists for  all $c\notin \mathcal{C}$.}

If part (1) of the question is true, then
Question~\ref{HC-questions}(1) would follow for simple  and hence
all Harish-Chandra modules.

%%%%%%%%%%%%%%%%%%%%%%%%%%%%%%%%%%%%%%%%%%%%%%%%%%

\section{A Conze embedding}\label{conze}

\subsection{} Let $Q$ be a minimal primitive ideal in the enveloping algebra
 $U(\mathfrak{g})$ of a semisimple
complex Lie algebra $\mathfrak{g}$.
A classic result of Conze \cite{con}
shows that $U(\mathfrak{g})/Q$ embeds in a Weyl algebra
$A_m=D(\C^m)$ in such a way that the restriction to $U(\mathfrak{g})/Q$
of the canonical $D(\C^m)$-module $\mathcal{O}(\C^m)$ is   a Verma module.
Moreover, $ D(\C^m)$ is a union of Harish-Chandra bimodules $P'_i$ coming
from the
translation principle \cite[Proposition~5.5]{JS}.

The main result of this section shows that the natural analogue of
this result holds for $U_c$. Indeed, as we show,
$U_c\hookrightarrow D(\h/\WW)=\bigcup P_i$,
where the $P_i$ are
the natural analogues of the Harish-Chandra bimodules
 $P_i'$: they are the duals of the progenerative bimodules
$B_{i0}(c)$. In particular, $D(\h/\WW)$ is a flat left $U_c$-module.
As might   be expected from the earlier sections, this result also
corresponds to a basic fact about $\hi$.
In the notation of \cite{haidis} we consider the principal affine cell
 $V=\mathcal{O}_{(1^n)}$
 in  $\hi$  corresponding to the partition $(1^n)$.
 It follows from
  \cite[Corollaries~2.7 and~2.8]{haidis} that   $V\cong \C^{2n-2}$
 is the open set $\mathcal{D}(\delta) = \hi\smallsetminus \V(\delta)$
 complementary to the zero set of $\delta$.
As we show, $\ogr D(\h/\WW)=\bigcup \delta^{-i}A^i=
\mathcal{O}(V)$.

\subsection{}\label{subsec-8.2} We now turn to the proofs of these
assertions.
 Pick $c\in \C$ satisfying Hypothesis~\ref{main-hyp}.

For $i \in \NN$ let $P^i_{i+1} = \delta^{-1}e_-H_{c+i+1}
e$\label{Pij-defn}
 and notice that, by \eqref{shift-defn}, this is  the dual of the progenerative
  $(\UU_{c+i+1}, \UU_{c+i})$-bimodule  $B_{i+1,i}(c)$. Thus,
  for $i>j\geq 0$,
  $$P^j_i =  P^j_{j+1} P^{j+1}_{j+2} \cdots   P^{i-1}_i
  = \left( B_{ij}(c)\right)^* $$
  is a progenerative $(U_{c+j},\, U_{c+i})$-bimodule.
By \cite[Lemma~6.11(1)]{GS},
$B_{ij}(c)\subseteq U_{c+i}\cap U_{c+i+1}$
and so, by induction,
 $U_c\subseteq P^0_i\subseteq P^0_{i+1}$, for all $i\geq 1.$

\subsection{Lemma}\label{subsec-8.3}
 {\it
For all $\ell\geq 1$,  $P^0_{\ell}\subseteq D(\h/\WW)$. }

\begin{proof} Pick $i\geq 0$. The algebra
$H_{c+i+1}$ acts naturally on $\C[\h]$ (see the proof
of \cite[Proposition~4.5]{EG}) and so
$P^{i}_{i+1}=\delta^{-1}e_-H_{c+i+1}e$ acts naturally on $\C[\reg{\h}]$.
It is important to notice that \eqref{conj} is an equality in
$D(\hr)\ast\WW$ and so the action of $H_{c+i}$ on $\C[\reg{\h}] $
induced from that of $P^i_{i+1}$ does equal the action coming
from the inclusion $H_{c+i}\subseteq D(\hr)\ast\WW$.

In fact, the action of $P^{i}_{i+1}$ on $\C[\reg{\h}] $   restricts to an
action on $\C[\h]^\WW=\C[\h/\WW]$. To see this note that, for $p\in
H_{c+i+1}$,
the element  $e_-pe$ maps $\C[\h]^\WW$  to $e_- \C[\h] =
\C[\h]^{\sign}$. Since
$\C[\h]^{\sign}=\delta \C[\h]^\WW$,
we conclude that $\delta^{-1} e_- p e$ maps $\C[\h]^{\WW}$ to
itself.

It remains to check that elements $\delta^{-1} e_- p e\in P^i_{i+1}$
 are actually differential
operators on $\C[\h]^{\WW}$. But if  $q\in
\C[\h]^{\WW}$ then $[q,\delta^{-1}e_- p e] = \delta^{-1}e_-
[q,p] e$.   By  \cite[Lemma~3.3(v)]{BEGqi}
$\C[\h]^{\WW}$ acts ad-nilpotently  on $H_{c+i}$ and so
 $q$ also acts ad-nilpotently on
$\delta^{-1} e_- p e\in P^i_{i+1}$.
By definition, this means that  $P^{i}_{i+1}\subseteq D(\h/\WW)$. By
induction, the
same is true for the  multiples $P^j_i$ of these spaces. \end{proof}

  \subsection{}\label{conze-com} We next want to understand
  $\bigcup \delta^{-i}A^i$. To begin, consider
  $\hin$ and the corresponding union $\bigcup_i\delta^{-i}\AAA^i$,
  where $\AAA^0 = \C[\C^{2n}]^{\WW}$,
  $\AAA^1 = \C[\C^{2n}]^{\sign}$, and $\AAA^i = (\AAA^1)^i$
  in the notation of
  \cite[(4.3)]{GS}.  We bigrade $\AAA^0$ and hence the $\AAA^i$
  by the Euler and order gradations; thus
   $\C[\C^{2n}] = \C[u_1,\dots,u_{n},v_1,\dots, v_{n}]$
  is graded by $\bideg(u_i) = (1,0)$ while
  $\bideg(v_j)=(-1,1)$, and we have chosen our variables so that $\WW$
  acts by permutation of the subscripts.

  In the notation of \cite[Section~2]{haidis} we want to consider
  the affine open cell $ U=U_{(1^n)}\subset \hin$.
  By \cite[Corollary~2.7]{haidis},  $\mathcal{O}(U)$ is generated
  by elements $\Delta_D\Delta_{(1^n)}^{-1}$ where, by
  \cite[(2.22) and the paragraph after (2.32)]{haidis}, the $\Delta_D$
  are a basis of $\AAA^1$ and $\Delta_{(1^n)}=\delta$. Thus
  $\mathcal{O}(U)= \bigcup \delta^{-i}\AAA^i=\C[\delta^{-1}\AAA^1]$.
Moreover, by
   \cite[Corollary~2.8]{haidis}, $\mathcal{O}_{U} =\mathcal{O}(U)$
   is a polynomial ring
  $R=\C[e_1,\dots e_{n}, a_1,\dots, a_{n}]$ (where our indexing   differs
by one
   from Haiman's).
 By \cite{haidis}, the $e_i$'s are the elementary
  symmetric functions in the $u_j$'s and so (possibly after reordering)
  $\bideg(e_j) = (j,0)$ for each $j$. Note that $e_1=\sum u_i=
  \mathbf{z}$, in the
  notation of  \eqref{hi-defn}.
   The $a_j$ are the coefficients of the Lagrange
  interpolation polynomial $\phi_a(x) = \sum_{i=1}^n a_{i}x^{i-1}
  = \sum_j v_j\prod_{k\not=j}(x-u_k)(u_j-u_k)^{-1}.$ Since this polynomial
   is homogeneous as a polynomial in either the $v$'s or $u$'s,
  with
 total  degrees $1$ in the $v$'s and   $0$ in the $u$'s,  it follows that
  $\bideg(a_j) = (-j,1)$ for each $j$.
In this case it is not easy to express     the element
 $\mathbf{z}^*=\sum v_i$  from
  \eqref{hi-defn} in terms of these variables, but  at least
 $\bideg(\mathbf{z}^*) = (-1,1)  = \bideg(a_1)$.

  Now consider $\hi$. By \cite[Lemma~4.9(1)]{GS} and the comments of the last
  paragraph,
  $A^0=\AAA^0/(e_1,  \mathbf{z}^*)$. Therefore,
  $$\bigcup_i \delta^{-i}A^i=\C[\delta^{-1}A^1]=
  \C[\delta^{-1}\AAA^1]/(e_1, \mathbf{z}^*) = R/(e_1, \mathbf{z}^*).
  $$
  Since $S=R/(e_1)$ is a domain, the element $\mathbf{z}^*$ is a regular
element in
  $S$ of the same bidegree as the element $a_1$.
  Thus the factor rings
  $\C[\delta^{-1}A^1]   =S/(\mathbf{z}^*)$ and
  $S/(a_1)$ have the same bigraded  Poincar\'e series.
  On the other hand, the restriction of $\mathcal{O}_{U}$ to
  $\hi$ is clearly just the open subvariety $\mathcal{D} (\delta)$
  complementary to $\V(\delta)$ in $\hi$. Combining these
  observations gives:

  \begin{lem} $\C[\delta^{-1}A^1]=\bigcup \delta^{-i}A^i= \mathcal{O}(V)$,
where
  $V=\mathcal{D} (\delta)$. Also, $ \C[\delta^{-1}A^1]$
has the same bigraded Poincar\'e series as the polynomial ring
$  \C[e_2,\dots e_{n}, a_2,\dots, a_{n}]$, where
  $\bideg(e_i)=(i,0)$ and $\bideg(a_i)=(-i,1)$.
  \qed
  \end{lem}

\subsection{}\label{subsec-8.4}
We can now prove our analogue of  the Conze embedding.

\begin{prop} Assume that $c\in \C$ satisfies
Hypothesis~\ref{main-hyp}. Then $D(\h/\WW)= \bigcup_{i\geq 0}
P_i^0 $. Moreover, $$\ogr D(\h/\WW)=\bigcup \ogr P_i^0=
\C[\delta^{-1}A^1]= \mathcal{O}(V),$$ where $V\cong \C^{2n-2}$ is
the open affine subvariety
 $\mathcal{D}(\delta)$ of $\hi$.
\end{prop}

\begin{proof} By \eqref{shift-defn} the $P_i^0$ are progenerators for
$U_c$.
 Give $P_i^0$ the induced filtration $\ord$ arising
from the order
filtration $\ord$ on $D(\hr)\ast \WW$. There is a potential ambiguity,
here,
 since $P^0_i$ can also be given an    order filtration
induced from  the order filtration
on $D(\h/\WW)$. However, thanks to \cite[Lemma~3.2.1 and
Proposition~3.5]{lev} the restriction map $D(\hr )^{\WW} \to
D(\hr/\WW)$ is a filtered isomorphism and so the two filtrations agree.
In particular, $\ogr  P_i^0 \subseteq \gr D(\h/\WW).$

If $j\geq 0$, then
\cite[Lemma~6.8]{GS} shows that
$$\ogr P_{j+1}^j  =  \delta^{-1}\ogr(e_-H_{c+j+1}e)
= \delta^{-1}e_-(\C[\h\oplus \h^*\ast \WW)e
= \Theta(\delta^{-1}A^1),$$
where $\Theta : \C[\reg{\h}\oplus \h^*]^\WW \to
e(\C[\reg{\h}\oplus \h^*]^\WW)e$ is the isomorphism given by $\Theta(f)=efe$.
Therefore, by induction and \cite[Lemma~6.7]{GS}
$$\Theta(\delta^{-i}A^i) \ = \  \Theta((\delta^{-1}A^1)^i)\ = \
(\ogr P_1^0)(\ogr P_2^1)\cdots (\ogr P_i^{i-1}) \ \subseteq\ \ogr P_i^0$$
and so, by Lemma~\ref{subsec-8.3},
\begin{equation}\label{conze3}  \C[\delta^{-1}A^1]\ = \
\bigcup_{i\geq 0} \delta^{-i} A^i \text{ is isomorphic to a subspace of }
 \bigcup_{i\geq 0}
    \ogr P_i^0\  \subseteq\  \ogr D(\h/\WW).
    \end{equation}

Clearly each $P^0_i$ is graded under the $\EE$-grading and so
  $\gr P_i^0$
has an induced $(\EE, \ord)$
bigrading, as in \eqref{conze-com}.
Now $\C[\h]^{\WW}= \C[e_2,\ldots ,e_n]$ is a polynomial ring whose
generators have degree and hence $\EE$-degree
$\Edeg(e_i)=i$.
Since the equation $\partial_{e_i}e_i-e_i\partial_{e_i}=1$
in $D(\h/\WW)$ is $\EE$-homogeneous,
$\Edeg \partial_{e_i} = -i$ for each $i$. Therefore
 $\ogr D(\h/\WW) =
\C[e_2,\ldots ,e_{n}, d_2,\ldots ,d_n]$ is a polynomial ring whose
generators satisfy $\bideg e_i = (i,0) $ and $ \bideg d_i =
(-i,1).$

By \eqref{conze3} and  Lemma~\ref{conze-com},
$\C[\delta^{-1}A^1]$ is therefore a bigraded algebra with exactly the same
bigraded Poincar\'e series as that of its overring $\ogr D(\h/\WW)$. Since
the bigraded components of these algebras are finite dimensional, they
must be equal.     In particular,  \eqref{conze3} implies that
     $\bigcup_{i\geq 0} \ogr P_i^0 = \ogr D(\h /\WW)$  and
    hence that $ \bigcup_{i\geq 0 } P_i^0 = D(\h/\WW)$. Finally,
    as $\ogr D(\h/\WW)$ is a polynomial ring, this also implies that
    $V= \spc(\C[\delta^{-1}A^1])$ is isomorphic to $\C^{2n-2}$.
\end{proof}

  \subsection{}\label{homdim} Since the $P_i^0$ are projective left
  $U_c$-modules we obtain:

\begin{cor} Keep the hypotheses  and notation of
Proposition~\ref{subsec-8.4}.
 Then $D(\h/\WW)$ is a flat left $U_c$-module. \qed
 \end{cor}

We remark that $D(\h/\WW)$ will typically not be
 a flat right module over $U_c$. Indeed this even fails for $n=2$, when
 $U_c$ is a primitive factor ring of $U(\mathfrak{sl}_2)$ (see
 \cite[(5.7)]{JS}).

%%%%%%%%%%%%%%%%%%%%%%%

%\clearpage

\section*{Index of Notation}\label{index}
\begin{multicols}{2}
{\small  \baselineskip 14pt

$A^1$, alternating polynomials, $A=\bigoplus A^k\delta^k$,
\hfill\eqref{A-1-defn}

$\widehat{A} = \bigoplus_{i\geq j \geq 0}
A^{i-j}\delta^{i-j}$,  \hfill\eqref{Aij-defn}

$\BB$, the tautological rank $n$  bundle, \hfill\eqref{tauto-defn}

$B=\bigoplus B_{ij}$,   \hfill\eqref{B-ring-defn}

$\mathcal{C}$, the set of singular values for $c$,
\hfill\eqref{curlyC-defn}

 $\widetilde{C}$,
${\widetilde{C}_{\lambda}}$, \hfill\eqref{tilde-C-defn}

$\Chr$,   characteristic varieties,
\hfill\eqref{charvar-defn},\eqref{bi-charvar}

 $\Chro$, associated  varieties,
  \hfill\eqref{charvar-defn},\eqref{bi-charvar}

$\Ch$, characteristic cycle, \hfill\eqref{charcycle-defn}

$\Chw$, restricted characteristic cycle, \hfill\eqref{cyclew-defn}

$d_\mu= \{ (i,j)\in \NN\times \NN : j <
\mu_{i+1}\}$, \hfill\eqref{d-mu-defn}

$\delta=\prod_{s\in \mathcal{S}} \alpha_s$, \hfill\eqref{delta-defn}

 $\Delta_c(\tau) $, the standard module, \hfill\eqref{standard-defn}

dominance ordering on $\irr\WW$,  \hfill\eqref{dominance-defn}

%Dunkl-Cherednik representation $\theta_c$, \hfill\eqref{theta-defn}

Euler operator $\EE=\sum x_i\delta_i$,  $\Edeg$, \hfill\eqref{Euler-defn}

$e, e_-$, trivial and sign idempotents, \hfill\eqref{e-defn}

fake degrees $f_\tau$, \hfill\eqref{fakedegrees}

good filtrations, \hfill \eqref{good-defn}

 $H_c$, the rational Cherednik algebra, $\h,\h^*$,  \hfill\eqref{hc-defn}

% Harish-Chandra bimodule, \hfill\eqref{HC-defn}

Hilbert schemes $\hin$, $\hi$, \hfill\eqref{hin-defn}

$I_\mu$, monomial ideal for partition $\mu$,  \hfill\eqref{I-eta-defn}

$\irr\WW$, the irred. reps. of $\WW$, \hfill\eqref{irred-defn}

$J^1= \cxy A^1$, \hfill\eqref{A-1-defn}

 $K_{\mu\lambda }(t,s)$, Kostka-Macdonald
coefficients, \hfill\eqref{kostka-defn}

$L_c(\tau)$, simple $\mathcal{O}_c$-module, \hfill\eqref{L-defn}

$\LL=\bigwedge^n \BB \cong\OO_{\hi}(1)$, \hfill\eqref{LLL-defn}

 $\mathcal{O}_c$, category $\mathcal{O}$ for $H_c$,
  \hfill\eqref{cat-O-defn}

 $\ord,\ogr$, order filtration and gradation,
 \hfill\eqref{order-filt-defn}

% $P_i^j =B_{ij}^*$,  \hfill\eqref{Pij-defn}

 $\PP=\rho_*\OO_{X_n}$,  the  Procesi bundle,
\hfill\eqref{PP-defn}

 $\PP_\mu$, $V_\mu$, \hfill\eqref{subsec-7.12}

$\Phi$,  $\widehat{\Phi}$,\hfill\eqref{Phi-defn}

Poincar\'{e} series $p(M,s,t)$, $p(M,v)$,
\hfill\eqref{bi-poincare},\eqref{mono-poincare}

$\qgr$, $\QGr$, \hfill\eqref{qgr-defn}

%$R(n,l)= \mathrm{H}^0(\hi, \PP\otimes \BB^{l})$, \hfill\eqref{RR-defn}

$\rho: X_n\to \hin$, \hfill\eqref{Xn-defn}

% $S_q = S_q(n,n)$, $q$-Schur algebra, \hfill\eqref{schur-defn}

% $\mathcal{S}$ the reflections in $\WW$,
% \hfill\eqref{involution-defn}

$\sigma(m)$, the principal symbol of $m$, \hfill\eqref{princ-symbol}

$\sign$, the sign representation of $\WW$, \hfill\eqref{sign-defn}

$\TT^2$-action \hfill\eqref{bigrad-defn0}, \eqref{bigrad-defn}

  $\widehat{\tau} : \hin\to  \C^{2n}/\WW$,  \hfill\eqref{tau-defn}

  $\tau : \hi\to \h\oplus\h^*/\WW$, \hfill\eqref{tau-defn}

$\tot$, $\tgr$,  total degree filtration/gradation \hfill\eqref{total-filt-defn}

$\triv$, the trivial $\WW$-representation, \hfill\eqref{triv-defn}

$U_c=eH_{c}e$, the spherical subalgebra, \hfill\eqref{spherical-defn}

%$W_q(\tau)$, $q$-Weyl module, \hfill\eqref{q-weyl-defn}

$\WW=\mathfrak{S}_n$, the symmetric group, \hfill\eqref{symmetric-defn}

%$X_n$, the isospectral Hilbert scheme, \hfill\eqref{Xn-defn}

$Z=Z(n)=\sigma^{-1}(\h/\WW)$, $Z_\lambda$, \hfill\eqref{Z-variety}

$\hio=Z_{(n)}$, the punctual Hilbert scheme,   \hfill\eqref{punctual-sect}

}
\end{multicols}

%\clearpage

\end{document}